\documentclass[twoside,
a4paper]{amsart}
\providecommand{\cosp}{\loglike{cosp}}
\providecommand{\sinp}{\loglike{sinp}}
\providecommand{\alli}{\iota}
\providecommand{\uir}[1]{\rho_{#1}}
\usepackage{xr-hyper}
\input{cyr.fd}
\usepackage{graphicx}

\PassOptionsToPackage{pdfauthor={Vladimir V. Kisil},%
    pdftitle={Erlangen Program at Large---2: Inventing a wheel. The parabolic one},%
    pdfsubject={mathematics},%
    backref=page,%
  pdfkeywords={symmetries, dual numbers, double numbers, Moebius map}
}{hyperref}
\IfFileExists{eulervm.sty}{\usepackage{eulervm}
}{}

\usepackage[printedin,extnum]{myamsart}
\ppnum{LEEDS--MATH--PURE--2007-07
}{arXiv: \href{http://arXiv.org/abs/0707.4024}{arXiv:0707.4024}}
{2007}

\input{mydef}
\IfFileExists{noweb.sty}{\usepackage{noweb}}{}
\makeatletter
\nwopt@smallcode
\let\tab=&
\def\idxexample#1{\nwix@id@uses#1}
\makeatother

\def\nwlbrace{\textbf{\texttt{\char123}}}
\def\nwrbrace{\textbf{\texttt{\char125}}}
\newcommand{\CPP}{\texttt{C++}}
\newcommand{\NoWEB}{\texttt{noweb}}

\providecommand{\GiNaC}{\textsf{GiNaC}}

\providecommand{\rs}{\mathring{\sigma}}
\providecommand{\clifford}[2][]{\ifcase #1 #2\or \tilde{#2} \or \breve{#2} \fi}

\begin{document}
\title[EPAL2: Inventing a parabolic wheel]
{Erlangen Program at Large---2:\\
Inventing a wheel. The parabolic one}

\author[Vladimir V. Kisil]%
{\href{http://maths.leeds.ac.uk/~kisilv/}{Vladimir V. Kisil}}
\thanks{On  leave from Odessa University.}

\address{%
School of Mathematics\\
University of Leeds\\
Leeds LS2\,9JT\\
UK
}

\email{\href{mailto:kisilv@maths.leeds.ac.uk}{kisilv@maths.leeds.ac.uk}}

\urladdr{\href{http://maths.leeds.ac.uk/~kisilv/}%
{http://maths.leeds.ac.uk/\~{}kisilv/}}

\dedicatory{Dedicated to 300\(^\text{th}\) anniversary of Leonhard Euler's birth}
\begin{abstract}
  We discuss parabolic versions of Euler's identity \[\rme^{\rmi
    t}=\cos t +\rmi \sin t.\] A purely algebraic approach based on
  dual numbers is known to produce a very trivial relation
  \(\rme^{\rmp t} = 1+\rmp t\). Therefore we use a geometric setup of
  parabolic rotations to recover the corresponding non-trivial
  algebraic framework.  Our main tool is M\"obius transformations which
  turn out to be closely related to induced representations of the
  group \(\SL\).
\end{abstract}
\AMSMSC{13A50}{08A99, 15A04, 20H05, 22D30, 51M10}
\maketitle

\tableofcontents

\section{Introduction: a Parabolic Wheel---an Algebraic Approach}
\label{sec:introduction}

A mathematical picture of a wheel, which uniformly rotates around its axis, is
given by the following ``model'':
\begin{equation}
  \label{eq:trigonometric}
  x=\cos t,\qquad y=\sin t, 
\end{equation}
where \(x\) and \(y\) denote the coordinates of a point on the unit
distance from the axis of rotation. The principal ingredients are
\emph{sine} and \emph{cosine} functions, which are known for more than
two thousand years.

\subsection{Complex Numbers}
\label{sec:complex-numbers}
\begin{figure}[htbp]
  \centering
  \includegraphics{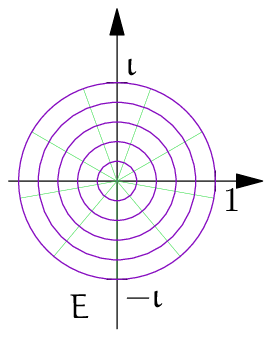}\qquad
  \includegraphics{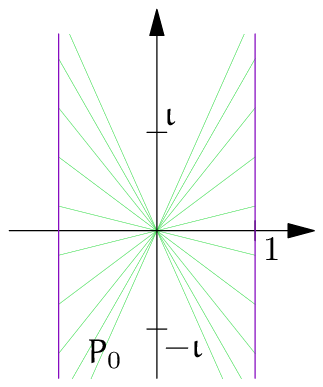}\qquad
  \includegraphics{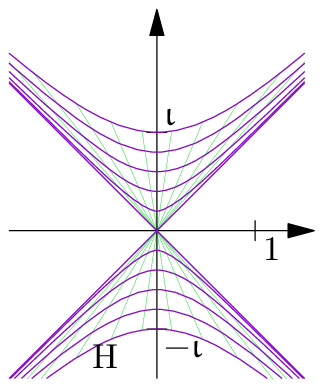}
  \caption[Rotations of wheels---Algebraic]{Rotations of algebraic
    wheels: elliptic (\(E\)), trivial parabolic (\(P_0\)) and
    hyperbolic (\(H\)). All blue rims of wheels are defined by the identity
    \(x^2-\alli^2y^2=1\). Green ``spokes'' (straight lines from the origin
    to a point on the rims) are ``rotated'' by multiplication by
    \(e^{\alli t}\).}
  \label{fig:rotations}
\end{figure}

A later invention of complex numbers \(z=x+\rmi y\), \(\rmi^2=-1\)
allows to write down two identities~\eqref{eq:trigonometric} as a
single one:
\begin{equation}
  \label{eq:complex}
  z=\cos t +\rmi \sin t.
\end{equation}
The next big advance is known as
\href{http://en.wikipedia.org/wiki/Euler's_formula}{Euler's formula},
which expresses trigonometric functions through the exponent of
an imaginary number:
\begin{equation}
  \label{eq:Euler}
  \rme^{\rmi t}=\cos t +\rmi \sin t.
\end{equation}
Thus the geometrical meaning of multiplication by \(\rme^{\rmi t}\) is an
isometric rotation of the plane \(\Space{R}{2}\), see
Fig.~\ref{fig:rotations}(E) with the (elliptic) metric given by:
\begin{equation}
  \label{eq:ell-metric}
  x^2+y^2=(x+\rmi y)(x-\rmi y).
\end{equation}
What are possible extensions of this results?

\subsection{Double Numbers}
\label{sec:double-numbers}

Complex numbers is not the only possible extension of the reals.
There are other variants of imaginary units, for example\footnote{Here
  \(\rmh\) is not a real number---a clarification which may be omitted
  in the case of \(\rmi^2=-1\).} \(\rmh^2=1\). Replacing \(\rmi\) by
\(\rmh\) in~\eqref{eq:Euler} we get 
a key to hyperbolic trigonometry:
\begin{equation}
  \label{eq:hyp-exp}
  \rme^{\rmh t}=\cosh t +\rmh \sinh t.
\end{equation}
Here expressions \(x+\rmh y\) form the algebra of \emph{double
  numbers}---the simplest case of
\href{http://en.wikipedia.org/wiki/Hypernumber}{hypernumbers}.
Multiplication by \(\rme^{\rmh t}\) is a map of double numbers into
itself which preserves the hyperbolic metric, cf.~\eqref{eq:ell-metric}:
\begin{equation}
  \label{eq:hyp-metric}
  x^2-y^2=(x+\rmh y)(x-\rmh y).
\end{equation}
Geometrically this may be viewed as \emph{hyperbolic rotation}, see
Fig.~\ref{fig:rotations}(H), in contrast to the elliptic
case~\eqref{eq:Euler}. 

\subsection{Dual number}
\label{sec:dual-number}

To make the picture complete we may wish to add the parabolic case
through the imaginary unit\footnote{We use different scripts of the
  epsilon: \(\rmh\) denote hyperbolic imaginary unit and
  \(\rmp\)---parabolic one.} of \emph{dual numbers} defined by
\(\rmp^2=0\). Since \(\rmp^n=0\) for any integer \(n>1\) we get from
Taylor's series of the exponent function the following identity:
\begin{equation}
  \label{eq:par-exp}
  \rme^{\rmp t}=1+\rmp t.
\end{equation}
Then parabolic rotations associated with \(\rme^{\rmp t}\) acts on dual
numbers as follows:
\begin{displaymath}
  \rme^{\rmp x}: a+\rmp b \mapsto a+\rmp (a x+b).
\end{displaymath}
This links the parabolic case with the Galilean
group~\cite{Yaglom79}.

Should we conclude, cf.~\cites{HerranzOrtegaSantander99a,Yaglom79}, from
here that: 
\begin{itemize}
\item the parabolic trigonometric functions are trivial:
  \begin{equation}
    \label{eq:par-trig-0}
    \cosp t =1, \qquad \sinp t=t\text{?}
  \end{equation}
\item the parabolic distance is independent from \(y\),
  cf.~\eqref{eq:ell-metric} and~\eqref{eq:hyp-metric}:
  \begin{equation}
    \label{eq:par-metr-0}
    x^2=(x+\rmh y)(x-\rmh y)\text{?}
  \end{equation}
\item the polar decomposition of a dual number is defined by~\cite{Yaglom79}*{App.~C(30')}:
  \begin{equation}
    \label{eq:p-polar-yaglom}
    u+\rmp v = u(1+\rmp \frac{v}{u}), \quad \text{ thus }
    \quad \modulus{u+\rmp v}=u, \quad \arg(u+\rmp v)=\frac{v}{u}\text{?}
  \end{equation}
\item the parabolic wheel looks rectangular, see  Fig.~\ref{fig:rotations}(\(P_0\))?
\end{itemize}

The analogies \eqref{eq:Euler}--\eqref{eq:hyp-exp}--\eqref{eq:par-exp}
and
\eqref{eq:ell-metric}--\eqref{eq:hyp-metric}--\eqref{eq:par-metr-0}
are quite explicit and widely accepted as an ultimate source for
parabolic
trigonometry~\cites{LavrentShabat77,HerranzOrtegaSantander99a,Yaglom79}. However
we will see shortly that there exists a less trivial form as well.

\begin{rem}
  The parabolic imaginary unit \(\rmp\) is a close relative to the
  infinitesimal number \(\varepsilon\) from non-standard
  analysis~\cites{Devis77,Uspenskii88}. The former has the property that its
  square is \emph{exactly} zero, meanwhile the square of the later is
  \emph{almost} zero at its own scale.
\end{rem}

\begin{rem}
  Introduction of double and dual numbers is not as artificial as it
  may looks from the traditional viewpoint, see
  Rem.~\ref{re:double-natural}. 
\end{rem}

\begin{rem}
  In cases when we need to consider simultaneously several imaginary
  units we use \(\alli\) to denote any of \(\rmi\), \(\rmp\),
  \(\rmh\).
\end{rem}

\section{A Parabolic Wheel---a Geometrical Viewpoint}
\label{sec:second-attempt}

We make a second attempt to describe parabolic rotations. If
multiplication (linear transformation) is not sophisticated enough for
this we would advance to the next level of complexity:
linear-fractional.

\subsection{Matrices}
\label{sec:matrices}
Imaginary units do not need to be seen as abstract quantities. 
We may realise them through zero-trace \(2\times 2\) matrices as
follows:
\begin{equation}
  \label{eq:imag-unit-matrices}
  \rmi=  \begin{pmatrix}
    0&1\\-1&0
  \end{pmatrix},   \qquad 
  \rmp=
  \begin{pmatrix}
    0&1\\0&0
  \end{pmatrix},   \qquad 
  \rmh=
  \begin{pmatrix}
    0&1\\1&0
  \end{pmatrix},
\end{equation}
with the parabolic \(\rmp\) nicely siting between the elliptic \(\rmi\)
and hyperbolic \(\rmh\).  Then the matrix multiplication implies
\(\rmi^2=-I\), \(\rmp^2=0\cdot I\), \(\rmh^2=I\), where \(I\) is the
\(2\times 2\) identity matrix. Correspondingly we have a matrix form of
the identities \eqref{eq:Euler}--\eqref{eq:hyp-exp}:
\begin{equation}
  \label{eq:matrix-exp-eh}
  \exp\begin{pmatrix}
    0&t\\-t&0
  \end{pmatrix}= 
  \begin{pmatrix}
    \cos t & \sin t \\ -\sin t & \cos t
  \end{pmatrix}, \qquad 
  \exp\begin{pmatrix}
    0&t\\t&0
  \end{pmatrix}= 
  \begin{pmatrix}
    \cosh t & \sinh t \\ \sinh t & \cosh t
  \end{pmatrix}.
\end{equation}
However the above pattern is only partially reproduced in the matrix
form of \eqref{eq:par-exp}: 
\begin{equation}
  \label{eq:parab-exp-matrix}
  \exp
  \begin{pmatrix}
    0&t\\0&0
  \end{pmatrix}=
  \begin{pmatrix}
    1&t\\0&1
  \end{pmatrix}.
\end{equation}
There is also some arbitrariness  in our choice of a matrix
representation for \(\rmp\), it may be equally well given by
the lower-triangular form:
\begin{equation}
  \label{eq:parab-prime-exp-matrix}
  \rmp'=
  \begin{pmatrix}
    0&0\\1&0
  \end{pmatrix}\quad 
  \text{implying}
  \quad
  \exp
  \begin{pmatrix}
    0&0\\t&0
  \end{pmatrix}=
  \begin{pmatrix}
    1&0\\t&1
  \end{pmatrix}.
\end{equation}
On the first glance this is not a radical difference, however, it
does have some implications for the M\"obius transform from
Section~\ref{sec:mobius-maps}.

\subsection{Cayley Transform}
\label{sec:cayley-transform}

Another matrix form of the identity~\eqref{eq:Euler} is
provided by the Cayley transform: 
\begin{equation}
  \label{eq:ell-cayley}
  \frac{1}{2}
  \begin{pmatrix}
    1 & -\rmi \\ -\rmi &1
  \end{pmatrix}
  \begin{pmatrix}
    \cos t & -\sin t \\ \sin t & \cos t
  \end{pmatrix}
  \begin{pmatrix}
    1 & \rmi \\ \rmi &1
  \end{pmatrix} =
  \begin{pmatrix}
    \rme^{\rmi t} & 0 \\ 0 &  \rme^{-\rmi t}
  \end{pmatrix}, 
\end{equation}
where the matrix
\begin{equation}
  \label{eq:cayley-matr}
  C_\rmi= \frac{1}{\sqrt{2}} 
  \begin{pmatrix}
    1 & -\rmi \\ -\rmi &1
  \end{pmatrix}
\end{equation}
is the Cayley transform  from the upper-half plane to the unit
disk. It have its hyperbolic cousin 
\begin{equation}
  \label{eq:hyp-cayley-matr}
  C_\rmh=
  \frac{1}{\sqrt{2}} 
  \begin{pmatrix} 
    1 & \rmh \\ -\rmh &1
  \end{pmatrix},
\end{equation}
which produces a matrix form of~\eqref{eq:hyp-exp}:
\begin{equation}
  \label{eq:hyp-cayley}
  \frac{1}{2}
  \begin{pmatrix}
    1 & \rmh \\ -\rmh &1
  \end{pmatrix}
  \begin{pmatrix}
    \cosh t & \sinh t \\ \sinh t & \cosh t
  \end{pmatrix}
  \begin{pmatrix}
    1 & -\rmh \\ \rmh &1
  \end{pmatrix} =
  \begin{pmatrix}
    \rme^{\rmh t} & 0 \\ 0 &  \rme^{-\rmh t}
  \end{pmatrix}.
\end{equation}
In the parabolic case we use the same pattern as
in~\eqref{eq:cayley-matr} and~\eqref{eq:hyp-cayley-matr}:
\begin{displaymath}
  C_\rmp=
  \begin{pmatrix} 
    1 & -\rmp \\ -\rmp &1  
  \end{pmatrix}
\end{displaymath}
The Cayley transform of matrix~\eqref{eq:parab-exp-matrix} is: 
\begin{equation}
  \label{eq:par-cayley}
  \begin{pmatrix}
    1 & -\rmp\\ -\rmp &1
  \end{pmatrix}
  \begin{pmatrix}
    1 & t \\ 0 & 1
  \end{pmatrix}
  \begin{pmatrix}
    1 & \rmp \\ \rmp &1
  \end{pmatrix} =
  \begin{pmatrix}
    1+\rmp t & t \\ 0 & 1-\rmp t
  \end{pmatrix}=
    \begin{pmatrix}
    \rme^{\rmp t} & t \\ 0 & \rme^{-\rmp t}
  \end{pmatrix}
\end{equation}
This is again not far from the previous
identities~\eqref{eq:ell-cayley} and~\eqref{eq:hyp-cayley}, however, the
off-diagonal \((1,2)\)-term destroys harmony. 

\begin{rem}
  \label{re:sl2-intro}
  It is not senseless to consider three
  matrices~\eqref{eq:imag-unit-matrices}, which materialise the
  imaginary units, together. In fact those trace-less matrices
  form a basis of the \(\algebra{sl}_2\) Lie algebra of the group
  \(\SL\)~\cite{Lang85}. Moreover they are generators of the
  one-parameter subgroups \(K\), \(N\), \(A\) correspondingly, which
  form the Iwasawa decomposition \(\SL=ANK\) of the group \(\SL\), see
  Section~\ref{sec:induc-repr}.
\end{rem}

\subsection[Moebius Maps]{M\"obius maps}
\label{sec:mobius-maps}
\begin{figure}[htbp]
  \centering
  \includegraphics{parab-rot-k.eps}\hfill
  \includegraphics{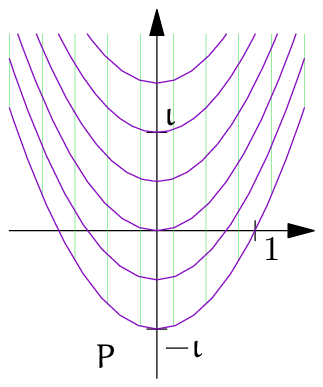}\hfill
  \includegraphics{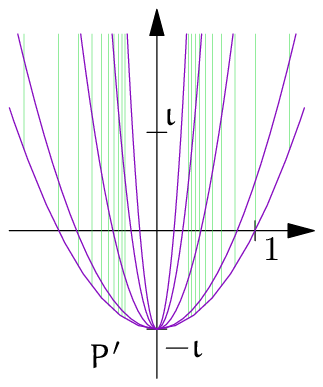}\hfill
  \includegraphics{parab-rot-a.eps}
  \caption{Rotation of geometric wheels: elliptic (\(E\)), two parabolic (\(P\)
    and \(P'\)) and hyperbolic (\(H\)). Blue orbits are level lines for the
    respective moduli. Green straight lines join points with the same value of
    argument and are drawn with the constant ``angular
    step'' in each case.}
  \label{fig:p-rotations}
\end{figure}

The matrix version of Euler's identity from the previous section can be
folded back to numbers through the linear-fractional (or M\"obius)
transformations. Indeed any \(2\times 2\) matrix define a map:
\begin{equation}
  \label{eq:moebius-def}
  \begin{pmatrix}
    a&b\\c&d
  \end{pmatrix}: z\mapsto \frac{az+b}{cz+d}, \text{ where } z=u+\alli v.
\end{equation}
Notably this is a group homomorphism of invertible matrices under
multiplication (or the group \(\SL\)) into transformations of conformally
completed plane~\cites{HerranzSantander02b,Kisil06b}. More precisely,
see~\cite{Kisil05a} or~\cite{Kisil06a} (an easy-reading), real-valued
matrices~\eqref{eq:matrix-exp-eh} and~\eqref{eq:parab-exp-matrix} act
as transformations of the ``upper half-plane'';
matrices~\eqref{eq:ell-cayley}, \eqref{eq:hyp-cayley}
and~\eqref{eq:par-cayley} act as transformations of the respective
``unit disk'', see Fig.~\ref{fig:p-rotations}. Those unit disks are
images of the upper half-planes under respective Cayley
transforms~\cite{Kisil05a}*{\S~\ref{E-sec:unit-circles}}.

The actions of diagonal matrices from the right-hand side of
identities~\eqref{eq:ell-cayley} and~\eqref{eq:hyp-cayley} are
straightforward: they are multiplications by \(\rme^{-2\rmi t}\) and
\(\rme^{-2\rmh t}\) correspondingly. Notably the images of the point
\(-\alli\) are:
\begin{equation}
  \label{eq:sin-cos-from-Moebius}
  \begin{pmatrix}
    \rme^{\rmi t}&0\\0&\rme^{-\rmi t}
  \end{pmatrix}: -\rmi \mapsto \sin 2t - \rmi \cos 2t; 
  \quad
  \begin{pmatrix}
    \rme^{\rmh t}&0\\0&\rme^{-\rmh t}
  \end{pmatrix}: -\rmh \mapsto - \sinh 2t- \rmh\cosh 2t .
\end{equation}
However the parabolic action of matrix~\eqref{eq:par-cayley} in
formula~\eqref{eq:moebius-def} is not such a simple one. 

\subsubsection{The Upper-Triangular Subgroup}
\label{sec:upper-triang-subgr}
The parabolic version of the relations~\eqref{eq:sin-cos-from-Moebius} 
with the upper-triangular matrices from \(N\) becomes:
\begin{equation}
  \label{eq:par-moebius}
  \begin{pmatrix}
    \rme^{\rmp t}&t\\0&\rme^{-\rmp t}
  \end{pmatrix}: -\rmp \mapsto t +\rmp (t^2-1).
\end{equation}
This coincides with the \emph{cyclic rotations} defined
in~\cite{Yaglom79}*{\S~8}.  A comparison of this result
with~\eqref{eq:sin-cos-from-Moebius} seemingly confirms that \(\sinp
t=t\) but suggest a new expression for \(\cosp t\):
\begin{displaymath}
  \cosp t = 1-t^2, \qquad \sinp t= t.
\end{displaymath}
Therefore the parabolic Pythagoras' identity would be:
\begin{equation}
  \label{eq:pyhagoras-p}
  \sinp^2 t +  \cosp t =1,
\end{equation}
which nicely fits in between of the elliptic and hyperbolic versions:
\begin{displaymath}
  \sin^2 t+\cos^2 t  =1, \qquad \sinh^2 t -  \cosh^2 t =-1.
\end{displaymath}
The identity~\eqref{eq:pyhagoras-p} is also less trivial than
the version \( \cosp^2 t =1\)
from~\cite{HerranzOrtegaSantander99a} (see
also~\eqref{eq:par-trig-0}, \eqref{eq:par-metr-0}). 

\subsubsection{The Lower-Triangular Subgroup}
\label{sec:lower-triang-subgr}
There is also the second option to define parabolic rotations, it is
generated by the lower-triangular variant of the above construction,
cf.~\eqref{eq:parab-prime-exp-matrix}. The important difference now
is: the reference point cannot be \(-\rmp\) since it is a fixed
point---as well as any point on the vertical axis. Instead we take
\(\rmp^{-1}\), which is an ideal element (a point at
infinity~\cite{Yaglom79}*{App.~C}) since \(\rmp\) is a divisor of
zero. The proper compactifications by ideal elements for all three
cases were discussed in~\cites{Kisil05a,Kisil06b}.

We denote the subgroup of lower-triangular matrices by \(N'\). We obtain
with it:
\begin{equation}
  \label{eq:par-moebius-prime}
  \begin{pmatrix}
    \rme^{\rmp t}&0\\t&\rme^{-\rmp t}
  \end{pmatrix}: \frac{1}{\rmp} \mapsto \frac{1}{t}+ \rmp \left(\frac{1}{t^2}-1\right). 
\end{equation}

A comparison with~\eqref{eq:par-moebius} shows that this form is obtained
by the change \(t\mapsto t^{-1}\). The same transformation gives new
expressions for parabolic trigonometric functions. The parabolic
``unit circle'' (or \emph{cycle}~\cites{Kisil05a,Yaglom79}) is
defined by the equation \(x^2-y=1\) in both cases . However other
orbits are different and we will give their description in the next
Section.

\section{Rebuilding Algebraic Structures from Geometry}
\label{sec:rebu-algebr-struct}

Rotations in elliptic and hyperbolic cases are given by products of
complex or double numbers correspondingly, however the multiplication
of dual numbers produces only the trivial parabolic rotation from
Fig.~\ref{fig:rotations}(\(P_0\)) rather than more interesting ones
from Fig.~\ref{fig:p-rotations}(\(P\)) or
Fig.~\ref{fig:p-rotations}(\(P'\)).  Also the coordinate-wise addition of vectors on
the plane is invariant under elliptic and hyperbolic rotations but is not
under the parabolic one. Can we find such algebraic operations for
vectors which will be compatible with parabolic rotations?

It is common in mathematics to ``revert a theorem into a definition''
and we will use it systematically in this section to recover a
compatible algebraic structure.

\subsection{Modulus and Argument}
\label{sec:modulus-argument}

In the elliptic and hyperbolic cases orbits of rotations are points
with the constant norm (modulus): either \(x^2+y^2\) or
\(x^2-y^2\). In the parabolic case we employ this point of view
as well:
\begin{defn}
  \label{de:norm}
  Orbits of actions~\eqref{eq:par-moebius}
  and~\eqref{eq:par-moebius-prime} are contour line for the following
  functions which we call respective moduli (norms):
  \begin{equation}
    \label{eq:parab-norm}
    \text{ for } N:\ \modulus{u+\rmp v}=u^2-v, \qquad \text{ for } N':\
    \modulus{u+\rmp v}'=\frac{u^2}{v+1}.
  \end{equation}
\end{defn}
\begin{rem}
  \begin{enumerate}
  \item The expression \( \modulus{(u,v)}=u^2-v\) represents a parabolic distance
    from \((0,\frac{1}{2})\) to \((u,v)\),
    see~\cite{Kisil05a}*{Lem.~\ref{E-le:n-orbits-concentric}} and is in
    line with the ``parabolic Pythagoras' identity''~\eqref{eq:pyhagoras-p}.
  \item Modulus for \(N'\) expresses the parabolic focal length from \((0,-1)\)
    to  \((u,v)\) as described in~\cite{Kisil05a}*{Lem.~\ref{E-le:np-orbits-p-confocal}}.
  \end{enumerate}
\end{rem}

The only straight lines preserved by the both parabolic rotations
\(N\) and \(N'\) are vertical lines, thus we will treat them as
``spokes'' for parabolic wheels. Elliptic spokes in mathematical
terms are ``points on the complex plane with the same argument'', thus
we again use it for the parabolic definition:
\begin{defn}
  \label{de:arg}
  Parabolic arguments are defined as follows:
  \begin{equation}
    \label{eq:parab-arg}
    \text{ for } N:\ \arg(u+\rmp v)=u, \qquad \text{ for } N':\
    \arg'(u+\rmp v)=\frac{1}{u}.
  \end{equation}
\end{defn}
Both Definitions~\ref{de:norm} and~\ref{de:arg} possess natural
properties with respect to parabolic rotations:
\begin{prop}
  \label{pr:rot-norm-arg}
  Let \(w_s\) is a parabolic rotation of \(w\) by angle \(s=t\)
  in~\eqref{eq:par-moebius} or \(s=t^{-1}\)
  in~\eqref{eq:par-moebius-prime}. Then:
  \begin{displaymath}
    \modulus{w_s}^{(\prime)}=\modulus{w},\qquad \arg^{(\prime)} w_s=
    \arg^{(\prime)} w+s,
  \end{displaymath}
  where primed versions are used for subgroup \(N'\).
\end{prop}
\begin{rem}
  Note that in the commonly accepted
  approach~\cite{Yaglom79}*{App.~C(30')} parabolic modulus and
  argument are given by expressions~\eqref{eq:p-polar-yaglom}, which
  are in a sense opposite to our conclusions.
\end{rem}

\subsection{Rotation as Multiplication}
\label{sec:rotat-as-mult}

We revert again theorems into definitions to assign multiplication.
In fact we require an extended version of properties stated in
Proposition~\ref{pr:rot-norm-arg}:
\begin{defn}
  \label{de:product}
  The product of vectors \(w_1\) and \(w_2\) is defined by the following
  two conditions:  
  \begin{enumerate}
  \item \(\arg(w_1 w_2)=\arg w_1 + \arg w_2\);
  \item \(\modulus{w_1 w_2} =\modulus{w_1}\cdot \modulus{w_2}\).
  \end{enumerate}
\end{defn}
We also need a special form of parabolic conjugation.
\begin{defn}
  Parabolic conjugation is given by \(\overline{u+\rmp v}=-u+\rmp v\).
\end{defn}
Combination of Definitions~\ref{de:norm}, \ref{de:arg}
and~\ref{de:product} uniquely determine expressions for products.
\begin{prop}
  The parabolic product of vectors is defined by formulae:
  \begin{align}
    \label{eq:parab-prod}
    \text{for } &N:& (u,v)*(u',v') & =
    (u+u',(u+u')^2-(v-u^2)(v'-u'^2)); \\
    \label{eq:parab-prime-prod}
    \text{for } &N':&(u,v)*(u',v') & = \left(\frac{uu'}{u+u'},\frac{(v+1)(v'+1)}{(u+u')^2}-1\right).
  \end{align}
\end{prop}
Although both expressions looks unusual they have many familiar properties:
\begin{prop}
  Both products~\eqref{eq:parab-prod} and~\eqref{eq:parab-prime-prod}
  satisfy to the following conditions:
  \begin{enumerate}
  \item\label{item:prod-comm-ass} They are commutative and
    associative;
  \item\label{it:rot-as-mult} The respective
    rotations~\eqref{eq:par-moebius} and~\eqref{eq:par-moebius-prime}
    are given by multiplications with a dual number with the unit norm.
  \item\label{item:prod-inv} The product \(w_1\bar{w}_2\) is invariant
    under respective rotations~\eqref{eq:par-moebius}
    and~\eqref{eq:par-moebius-prime} .
  \item\label{item:prod-norm-sq} The second component of the product
    \(w\bar{w}\) is \(\modulus{w}^2\).
  \end{enumerate}
\end{prop}

\subsection{Invariant Linear Algebra}
\label{sec:invar-line-algebra}

Now we wish to define a linear structure on \(\Space{R}{2}\) which would
be invariant under point multiplication from the previous Subsection
(and thus under the parabolic rotations,
cf.~Prop.\ref{it:rot-as-mult}). Multiplication by a scalar is
straightforward (at list for a positive scalar): it should preserve
the argument and scale the norm of vectors. Thus we have formulae for
\(a>0\):
\begin{eqnarray}
  \label{eq:scalar-prod}
  a\cdot (u,v)&=&(u,a v+u^2(1-a))\quad\text{for } N,\\
  \label{eq:scalar-prod-prime}
  a\cdot (u,v)&=&\left(u,\frac{v+1}{a}-1\right)\quad\text{for } N'.
\end{eqnarray}

On the other hand addition of vectors can be done in several different
ways. We present two solutions: one is tropical and another---exotic.

\subsubsection{Tropical form}
\label{sec:tropical-form}
Let us introduce the lexicographic order on \(\Space{R}{2}\):
\begin{displaymath}
  (u,v)\prec(u',v') \quad \text{if and only if} \quad
  \left\{\begin{array}{ll}
      \text{either}& u<u'; \\
      \text{or}& u=u',\  v<v'.
  \end{array}\right.
\end{displaymath}
One can define functions \(\min\) and \(\max\) of a pairs of points
from on \(\Space{R}{2}\) correspondingly. An addition of two vectors
can be defined either as their minimum or maximum. A similar
definition is used in \emph{tropical mathematics}, also known as
Maslov dequantisation or \(\Space[\min]{R}{}\) and
\(\Space[\max]{R}{}\) algebras, see~\cite{Litvinov05} for a
comprehensive survey. It is easy to check that such an addition is
distributive with respect to vector
multiplications~\eqref{eq:scalar-prod}---\eqref{eq:scalar-prod-prime}
and consequently is invariant under 
parabolic rotations. Although it looks promising to investigate this
framework we do not study it further for now.

\subsubsection{Exotic form}
\label{sec:exsotic-form}
Addition of vectors for both subgroups \(N\) and \(N'\) can be defined
by the common rules, where subtle differences are hidden within
corresponding Definitions~\ref{de:norm} (norms) and~\ref{de:arg}
(arguments).
\begin{defn}
  \label{de:p-add}
  Parabolic addition of vectors is defined by the following formulae:
  \begin{eqnarray}
    \label{eq:p-add-arg-exotic}
    \arg^{(\prime)}(w_1+w_2)&=&\frac{\arg^{(\prime)} w_1\cdot
      \modulus{w_1}^{(\prime)} 
      +\arg^{(\prime)} w_2\cdot\modulus{w_2}^{(\prime)}}{\modulus{w_1+w_2}^{(\prime)}},\\
    \label{eq:p-add-norm-exotic}
    \modulus{w_1+w_2}^{(\prime)}&=&\modulus{w_1}^{(\prime)}+\modulus{w_2}^{(\prime)},
  \end{eqnarray}
  primed versions are used for the subgroup \(N'\).
\end{defn}
The rule for the norm of sum~\eqref{eq:p-add-norm-exotic} may looks
too trivial at a first glance. We should say in its defence that it
nicely sits in between of the elliptic  \(\modulus{w+w'}\leq
\modulus{w}+\modulus{w'}\) and hyperbolic \(\modulus{w+w'}\geq
\modulus{w}+\modulus{w'}\) inequalities for norms.

Both
formulae~\eqref{eq:p-add-arg-exotic}--\eqref{eq:p-add-norm-exotic}
together uniquely define explicit expressions for additions of
vectors.  Although those expressions are rather cumbersome and not
really much needed. Instead we list properties of this operations:
\begin{prop}
  Vector additions for subgroups \(N\) and \(N'\) defined
  by~\eqref{eq:p-add-arg-exotic}--\eqref{eq:p-add-norm-exotic} satisfy
  to the following conditions:
  \begin{enumerate}
  \item\label{item:add-is-comm-ass} They are commutative and
    associative.
  \item\label{item:distrib} They are distributive for
    multiplications~\eqref{eq:parab-prod} and~\eqref{eq:parab-prime-prod};
    consequently:
  \item\label{item:add-rot-inv} They are parabolic rotationally invariant;
  \item\label{item:distrib-scalar} They are distributive in both ways for the scalar
    multiplications~\eqref{eq:scalar-prod}
    and~\eqref{eq:scalar-prod-prime} respectively:
    \begin{displaymath}
      a\cdot(w_1+w_2)=a\cdot w_1+a\cdot w_2,\qquad
      (a+b)\cdot w=a\cdot w+b\cdot w.
    \end{displaymath}
  \end{enumerate}
\end{prop}
To complete the construction we need to define the zero vector and
inverse.
\begin{prop}
  \begin{itemize}
  \item[(\(N\))] The zero vector is \((0,0)\) and consequently the inverse
    of \((u,v)\) is \((u,2u^2-v)\).
  \item[(\(N'\))] The zero vector is \((\infty,-1)\) and consequently the inverse
    of \((u,v)\) is \((u,-v-2)\).
  \end{itemize}
\end{prop}
Consequently we can check that scalar
multiplications by negative reals are given by the same identities
~\eqref{eq:scalar-prod} and \eqref{eq:scalar-prod-prime} as for
positive ones. 

\subsubsection{The Real and Imaginary Parts}
\label{sec:real-imaginary-parts}
Having the vector addition at our hands we may wish to define the real
and imaginary parts compatible with it.
We can start from the familiar formulae
\(\frac{1}{2}(w+\bar{w})\) and \(\frac{1}{2}(w-\bar{w})\). While such
a real part has a reasonable value \((0,\modulus{w})\) (the subgroup
\(N\) case), the imaginary part suffers from the malformed denominator
in~\eqref{eq:p-add-arg-exotic} due to the fact
\(\modulus{w}=\modulus{\bar{w}}\). 

This should not be seen as a defect of the exotic addition. A moment
of reflection reveals that ``purely real'' dual numbers are naturally
defined as having zero argument, e.g. the vertical axis for the
subgroup \(N\). This seems to agree well with the above
real part.

However what are ``purely imaginary'' dual numbers? The horizontal
axis can be the first suggestion, but all its points have different
arguments. One may still prefer to choose that the number \((1,0)\)
would be purely imaginary (since \((0,1)\) is purely real) and all
imaginary numbers have the same argument. So purely imaginary numbers
\begin{defn}
  \label{de:imaginary_nums}
  For both subgroups \(N\) and \(N'\):
  \begin{enumerate}
  \item purely real numbers are defined by the condition \(\arg
    w= 0\).
  \item purely imaginary numbers are defined by the condition \(\arg
    w= 1\).
  \end{enumerate}
\end{defn}
  Since arguments of \(\Re w\) and \(\Im  w\) are fixed by this
  Definition we need only to find their moduli.
\begin{prop}
  The natural condition \(w=\Re w+ \Im w\) together with
  Defn.~\ref{de:imaginary_nums} uniquely define \(\Re w\) and \(\Im
  w\). It is also determined by the identities:
  \begin{displaymath}
    \modulus{\Re w} = (1-\arg w)\modulus{w}, \qquad     \modulus{\Im
      w} = \arg w\modulus{w}.
  \end{displaymath}
\end{prop}
The explicit formulae for the real and imaginary parts in the cases of
both subgroup \(N\) and \(N'\) can be found in
App.~\ref{sec:outp-symb-calc}. 

\subsubsection{Linearisation of the exotic form}
\label{sec:line-exot-form}

Some useful information can be obtained from the transformation
between the parabolic unit disk and its linearised model. In such
linearised coordinates \((a,b)\) the
addition~\eqref{eq:p-add-arg-exotic}--\eqref{eq:p-add-norm-exotic} is
done in the usual coor\-dinate\--wise manner:
\((a,b)+(a',b')=(a+a',b+b')\).  

To this end we calculate the value of \((u,v)=a\cdot(u_1,v_1)+b\cdot(u_2,v_2)\) for
\((u_1,v_1)=(1,0)\) and \((u_2,v_2)=(-1,0)\).  For the subgroups \(N\) the
transform is given by:
\begin{align}
  u&=\frac{a-b}{a+b},& v&=\frac{(a-b)^2}{(a+b)^2}-(a+b), &
  a&=\frac{u^2-v}{2}(1+u),& b&=\frac{u^2-v}{2}(1-u).
\end{align}
For the subgroup \(N'\) such a transformation is:
\begin{align}
  u&=\frac{a+b}{a-b},& v&=\frac{(a+b)}{(a-b)^2}-1,&
  a&=\frac{u(u+1)}{2(v+1)},& b&=\frac{u(u-1)}{2(v+1)}.
\end{align}
We also note that both norms~\eqref{eq:parab-norm} have exactly
the same value \(a+b\) in the respective \((a,b)\)-coordinates.

\begin{rem}
  \label{re:conformality}
  The irrelevance of the standard linear structure for parabolic
  rotations manifests itself in many different ways, e.g. in an
  apparent ``non-conformality'' of lengths from parabolic foci, e.g.
  with the parameter \(\rs=0\)
  in~\cite{Kisil05a}*{Prop.~\ref{E-it:conformity-length-foci}}.  An
  adjustment of notions to the proper framework restores the clear
  picture. 

  The initial definition of
  conformality~\cite{Kisil05a}*{Defn.~\ref{E-de:conformal}} considered
  the usual limit \(y'\rightarrow y\) along a straight line, i.e.
  ``spoke'' in terms of Fig.~\ref{fig:rotations}. This is justified in
  the elliptic and hyperbolic cases. However in the parabolic setting
  the right ``spokes'' are vertical lines, see
  Fig.~\ref{fig:p-rotations}, so the limit should be taken along
  them~\cite{Kisil05a}*{Prop.~\ref{E-pr:parab-conf}}.
\end{rem}

\section{Induced Representations as a Source of Imaginary Units}
\label{sec:induc-repr}

As we already mentioned in Rem.~\ref{re:sl2-intro} all three matrix
exponents \eqref{eq:matrix-exp-eh}--\eqref{eq:parab-exp-matrix} are
one-parameter subgroups of the group \(\SL\)---the group of \(2\times
2\) matrices with unit determinant. Moreover any one-parameter
subgroup of \(\SL\) is a conjugate to either of subgroup \(A\), \(N\) or
\(K\), see Rem.~\ref{re:sl2-intro} for their descriptions.
Thus our consideration may be applied for construction of
\emph{induced representations} of \(\SL\) in an extended meaning.

The general scheme of induced representations is as follows,
see~\cite{Kirillov76}*{\S~13.2}, \cite{Kisil97c}*{\S~3.1}. We
denote \(\SL\) by \(G\) and let \(H\) be its subgroup.  Let \(
\Omega=G / H\) be the corresponding homogeneous space and \(s: \Omega
\rightarrow G\) be a continuous function~\cite{Kirillov76}*{\S~13.2}
which is a left inverse to the natural projection \(G\rightarrow
G/H\). In our case we choose:
\begin{equation}
  \label{eq:s-map}
  s: (u,v) \mapsto
  \frac{1}{\sqrt{v}}
  \begin{pmatrix}
    v & u \\ 0 & 1
  \end{pmatrix}, \qquad (u,v)\in\Space{R}{2},\  v>0.
\end{equation}
Then any \(g\in G\) has a unique decomposition of the form
\(g=s(\omega)h\) where \(\omega\in \Omega\) and \(h\in H\). We will
write:
\begin{equation}
  \label{eq:r-map}
  \omega=s^{-1}(g),\qquad h=r(g):={(s^{-1}(g))}^{-1}g.
\end{equation}
Note that \(\Omega \) is a left homogeneous space with the
\(G\)-action defined in terms of \(s\) as follows:
\begin{equation}
  \label{eq:g-action}
  g: \omega  \mapsto g\cdot \omega=s^{-1}(g^{-1}* s(\omega)),
\end{equation}
where \(*\) is the multiplication on \(G\).

Let \(\chi: H \rightarrow \Space{R}{2} \) be a ``unitary character'' of
\(H\) in some generalised sense illustrated bellow.  Then
it induces a  ``unitary'' representation of \(G\), which is very
close to induced representations in the sense of
Mackey~\cite{Kirillov76}*{\S~13.2}.  This representation has the
canonical realisation \(\uir{}\) in the space \( \FSpace{L}{2}(\Omega)\)
of square integrable \(\Space{R}{2}\)-valued functions. It is given by
the formula (cf.~\cite{Kirillov76}*{\S~13.2.(7)--(9)}): 
\begin{equation} 
  \label{eq:def-ind}
  [\uir{\chi}(g) f](\omega)= \chi_0(r(g^{-1} * s(\omega)))  f(g\cdot \omega),
  \qquad 
  \chi_0(h)=\chi(h)\left( 
    \frac{d\mu(h\cdot \omega)}{d\mu(\omega)} \right)^{ \frac{1}{2} },
\end{equation}
where \(g\in G\), \(\omega\in\Omega\), \(h\in H\) and \(r: G
\rightarrow H\), \(s: \Omega \rightarrow G\) are maps defined
above; \(*\)~denotes multiplication on \(G\) and \(\cdot\) denotes the
action~\eqref{eq:g-action} of \(G\) on \(\Omega\) from the left.

\subsection{Induction from $K$}
\label{sec:induction-from-k}

This is the most traditional case in the representation theory. The
action~\eqref{eq:g-action} takes the form:
\begin{displaymath}
  \begin{pmatrix}
    a&b\\c&d
  \end{pmatrix}: (u,v)\mapsto
  \left(\frac{(au+b)(c u+d) +cav^2}{( c u+d)^2 +(cv)^2}, \frac{v}{( c u+d)^2 +(cv)^2}\right)
\end{displaymath}
Obviously it preserves the upper-half plane \(v>0\). Moreover with the
help of the imaginary unit \(\rmi^2=-1\) it can be naturally
represented as a M\"obius transformation:
\begin{displaymath}
  \begin{pmatrix}
    a&b\\c&d
  \end{pmatrix}: w\mapsto  \frac{aw+b}{c w+d}, \quad w=u+\rmi v.
\end{displaymath}

Thus it is justified \emph{in this particular case} to look for
complex-valued characters of \(K\). They are parametrised by an
integer \(n\in\Space{Z}{}\):
\begin{displaymath}
  \uir{n}(h_t)=e^{\rmi nt}=(\cos t+\rmi\sin t)^n, \qquad \text{ where } h_t=
  \begin{pmatrix}
    \cos t&-\sin t\\ \sin t & \cos t
  \end{pmatrix}.
\end{displaymath}
We can also calculate that:
\begin{displaymath}
  r(g^{-1} * s(\omega))=\frac{1}{\sqrt{( c u+d)^2 +(cv)^2}}
  \begin{pmatrix}
    cu+d&-cv\\cv&cu+d
  \end{pmatrix}.
\end{displaymath}
Taking into account the identity
\(\frac{\modulus{a}}{a}=\left(\frac{\bar{a}}{a}\right)^{\frac{1}{2}}\)
we obtain such a realisation of~\eqref{eq:def-ind}:
\begin{displaymath}
  [\uir{n}(g) f](w)=\left(\frac{{c}\bar{w}+{d}}{cw+d}\right)^{\frac{n}{2}}\ 
  f\left(\frac{aw+b}{c w+d}\right), \quad \text{ where }
  g^{-1}=\begin{pmatrix}
    a&b\\c&d
  \end{pmatrix}, \ 
  w=u+\rmi v.
\end{displaymath}

\subsection{Induction from $A$}
\label{sec:induction-from-a}

In this case the action~\eqref{eq:g-action} takes the form:
\begin{displaymath}
  \begin{pmatrix}
    a&b\\c&d
  \end{pmatrix}: (u,v)\mapsto
  \left(\frac{(au+b)(c u+d) -cav^2}{( c u+d)^2 -(cv)^2}, \frac{v}{( c u+d)^2 -(cv)^2}\right)
\end{displaymath}
This time the map \textbf{does not} preserve the upper-half plane
\(v>0\): the sign of \(( c u+d)^2 -(cv)^2\) is not determined. To
express this map as a M\"obius transformation we require the double
numbers imaginary unit \(\rmh ^2=1\):
\begin{displaymath}
  \begin{pmatrix}
    a&b\\c&d
  \end{pmatrix}: w\mapsto  \frac{aw+b}{c w+d}, \quad w=u+\rmh v.
\end{displaymath}
\begin{rem}
  \label{re:double-natural}
  As we can see now the double numbers naturally appear in relation with the
  group \(\SL\) and thus their introduction in
  \S~\ref{sec:double-numbers} was not ``a purely generalistic
  attempt'', cf.~\cite{Pontryagin86a}*{p.~4}. The same is true for dual
  numbers as can be seen in the next subsection.
\end{rem}
Under such conditions it does not make much sense to look for a
complex-valued characters of \(A\). Instead we will take double number
valued characters which are parametrised by a real number \(\sigma\):
\begin{displaymath}
  \uir{\sigma }(h_t)=e^{\rmh \sigma t}=(\cosh t+\rmh\sinh t)^\sigma,
  \qquad \text{ where } h_t= 
  \begin{pmatrix}
    \cosh t&\sinh t\\ \sinh t & \cosh t
  \end{pmatrix}.
\end{displaymath}
We can also calculate that
\begin{displaymath}
  r(g^{-1} * s(\omega))=\frac{1}{\sqrt{( c u+d)^2 -(cv)^2}}
  \begin{pmatrix}
    cu+d&cv\\cv&cu+d
  \end{pmatrix}.
\end{displaymath}
Thus the formula~\eqref{eq:def-ind} becomes:
\begin{displaymath}
  [\uir{\sigma}(g) f](w)=\left(\frac{{c}\bar{w}+{d}}{cw+d}\right)^{\frac{\sigma}{2}}\ 
  f\left(\frac{aw+b}{c w+d}\right), \quad \text{ where }
  g^{-1}=\begin{pmatrix}
    a&b\\c&d
  \end{pmatrix}
  , \ w=u+\rmh v.
\end{displaymath}

\subsection{Induction from $N$}
\label{sec:induction-from-n}

We consider here the lower-triangular matrices forming the subgroup
\(N'\). The action~\eqref{eq:g-action} takes now the form:
\begin{displaymath}
  \begin{pmatrix}
    a&b\\c&d
  \end{pmatrix}: (u,v)\mapsto
  \left(\frac{a u+b}{c u+d}, \frac{v}{( c u+d)^2}\right)
\end{displaymath}
This map preserves the upper-half plane \(v>0\) as the elliptic case
of \(K\). To express this map as a M\"obius transformation we require
the dual numbers imaginary unit \(\rmp ^2=0\):
\begin{displaymath}
  \begin{pmatrix}
    a&b\\c&d
  \end{pmatrix}: w\mapsto  \frac{aw+b}{c w+d}, \quad w=u+\rmp v.
\end{displaymath}
Similarly to the previous hyperbolic case we would look not for
complex-valued characters but rather use parabolic rotations described
in Section~\ref{sec:mobius-maps}. A distinction from the hyperbolic
case is that they are not given by multiplication of double
numbers. Such a ``character'' parametrised by a real number \(\kappa\)
is defined by:
\begin{displaymath}
  \uir{\kappa }(h_t): w  \mapsto \frac{(1+\rmp \kappa t) w +\kappa t}{(1-\rmp
    \kappa t)}
  =(1+2\rmp \kappa t) w +(\kappa t+\rmp \kappa^2 t^2)
  , \quad \text{ where } h_t=
  \begin{pmatrix}
    1&0\\ t & 1
  \end{pmatrix}
  .
\end{displaymath}
Furthermore we calculate that
\begin{displaymath}
  r(g^{-1} * s(\omega))=
  \begin{pmatrix}
    1 & 0 \\ \frac{v c}{ c u+d} &1
  \end{pmatrix}.
\end{displaymath}
Thus the formula~\eqref{eq:def-ind} has the following realisation:
\begin{displaymath}
  [\uir{\kappa }(g) f](w)=
\left(1-\rmp\frac{2\kappa v c}{ c u+d}\right)\,f\left(\frac{aw+b}{c
    w+d}\right)-\frac{\kappa v c}{ c u+d}
+\rmp \frac{(\kappa v c)^2}{( c u+d)^2}
\end{displaymath}
\begin{displaymath}
   \quad \text{ where }
  g^{-1}=\begin{pmatrix}
    a&b\\c&d
  \end{pmatrix}, \quad 
  w=u+\rmp v.
\end{displaymath}
The vector space of functions, where this representation acts, should
be also considered with linear operations defined in
\S~\ref{sec:exsotic-form}.

These three examples will be used to build the corresponding versions of
the Cauchy integral formula along the Erlangen Program at Large
outlined in~\cites{Kisil97c,Kisil05a}.

\section*{Acknowledgments}
\label{sec:acknowledgments}
I am grateful to Prof.~S.L.~Blyumin for pointing out the
paper~\cite{Litvinov05} to my attention. Anastasia Kisil read the
manuscript and made useful suggestions.

{\small
\bibliography{abbrevmr,akisil,analyse,arare,aclifford,aphysics}
}

\appendix

\IfFileExists{parab-rotation.d}{
\section{Output of Symbolic Calculations}
\label{sec:outp-symb-calc} 
Here are the results of our symbolic calculations. The source code can
be obtained from this paper~\cite{Kisil07a} source at
\url{http://arXiv.org}. It uses Clifford algebra
facilities~\cite{Kisil05b} of the \GiNaC\ library~\cite{GiNaC}. The
source code is written in \NoWEB~\cite{NoWEB} literature programming
environment.

\noindent Cayley of the matrix x: \(\left(\begin{array}{cc}- (\clifford[0]{e}^{{1} }*\clifford[0]{e}^{{0} }) x+\clifford[0]{\mathbf{1}}& \clifford[0]{e}^{{0} } x\\0&- (\clifford[0]{e}^{{1} }*\clifford[0]{e}^{{0} }) x+\clifford[0]{\mathbf{1}}\end{array}\right)\)\\
Rotation by x: \(\left(\begin{array}{cc}u+x&x^{2}+2  u x+v\end{array}\right)\)\\
Rotation of \((u_0, v_0)\) by \(x\): \(\left(\begin{array}{cc}x&-1+x^{2}\end{array}\right)\)\\
Parabolic norm: \(-v+u^{2}\)\\
Real number x as a dual number: \(\left(\begin{array}{cc}0&- x\end{array}\right)\)\\
Product: \(\left(\begin{array}{cc}u+u'&- u^{2} u'^{2}+u'^{2}+ u'^{2} v- v' v+2  u u'+ u^{2} v'+u^{2}\end{array}\right)\)\\
Product by a scalar: \(\left(\begin{array}{cc}u& a v+u^{2}- u^{2} a\end{array}\right)\)\\
Real part: \(\left(\begin{array}{cc}0&- u v+v-u^{2}+u^{3}\end{array}\right)\)\\
Imag part: \(\left(\begin{array}{cc}1&1+ u v-u^{3}\end{array}\right)\)\\
Lin comb of two vectors a*(1, 0)+b*(-1, 0): \(\left(\begin{array}{cc}\frac{a-b}{a+b}&-\frac{3  a b^{2}+2  a b+b^{3}-b^{2}-a^{2}+3  a^{2} b+a^{3}}{{(a+b)}^{2}}\end{array}\right)\)\\
P is the sum Re(P) and Im(P): \textbf{true}\\
The real part of a real dual number is itself: \textbf{true}\\
norm is invariant: \textbf{true}\\
Product is invariant: \textbf{true}\\
Product is norm squared: \textbf{true}\\
Product \((u, v)*(u_0, v_0)\) is \((u, v)\): \textbf{true}\\
Add is commutative: \textbf{true}\\
Add is associative: \textbf{true}\\
S-mult commutative: \textbf{true}\\
S-mult associative: \textbf{true}\\
S-mult distributive 1: \textbf{true}\\
S-mult distributive 2: \textbf{true}\\
Product is symmetric (commutative): \textbf{true}\\
Prod is associative: \textbf{true}\\
Product is distributive: \textbf{true}\\
\vspace{2mm}\hrule
Cayley of the matrix x: \(\left(\begin{array}{cc}- (\clifford[0]{e}^{{1} }*\clifford[0]{e}^{{0} }) x+\clifford[0]{\mathbf{1}}&0\\ \clifford[0]{e}^{{0} } x&- (\clifford[0]{e}^{{1} }*\clifford[0]{e}^{{0} }) x+\clifford[0]{\mathbf{1}}\end{array}\right)\)\\
Rotation by x: \(\left(\begin{array}{cc}-\frac{u}{-1+ u x}&-\frac{ u^{2} x^{2}-2  u x-v}{1+ u^{2} x^{2}-2  u x}\end{array}\right)\)\\
Rotation of \((u_0, v_0)\) by \(x\): \(\left(\begin{array}{cc}-\frac{1}{x}&-\frac{-1+x^{2}}{x^{2}}\end{array}\right)\)\\
Parabolic norm: \(\frac{u^{2}}{1+v}\)\\
Real number x as a dual number: \(\left(\begin{array}{cc}\infty&\frac{\infty^{2}-x}{x}\end{array}\right)\)\\
Product: \(\left(\begin{array}{cc}\frac{ u u'}{u+u'}&-\frac{-1+u'^{2}-v'- v' v+2  u u'-v+u^{2}}{{(u+u')}^{2}}\end{array}\right)\)\\
Product by a scalar: \(\left(\begin{array}{cc}u&-\frac{-1+a-v}{a}\end{array}\right)\)\\
Real part: \(\left(\begin{array}{cc}\infty&\frac{u+\infty^{2}+ \infty^{2} v-u^{2}}{ {(-1+u)} u}\end{array}\right)\)\\
Imag part: \(\left(\begin{array}{cc}1&-\frac{-1+u-v}{u}\end{array}\right)\)\\
Lin comb of two vectors a*(1, 0)+b*(-1, 0): \(\left(\begin{array}{cc}\frac{a+b}{a-b}&\frac{a+2  a b-b^{2}-a^{2}+b}{{(a-b)}^{2}}\end{array}\right)\)\\
P is the sum Re(P) and Im(P): \textbf{true}\\
The real part of a real dual number is itself: \textbf{true}\\
norm is invariant: \textbf{true}\\
Product is invariant: \textbf{true}\\
Product is norm squared: \textbf{true}\\
Product \((u, v)*(u_0, v_0)\) is \((u, v)\): \textbf{true}\\
Add is commutative: \textbf{true}\\
Add is associative: \textbf{true}\\
S-mult commutative: \textbf{true}\\
S-mult associative: \textbf{true}\\
S-mult distributive 1: \textbf{true}\\
S-mult distributive 2: \textbf{true}\\
Product is symmetric (commutative): \textbf{true}\\
Prod is associative: \textbf{true}\\
Product is distributive: \textbf{true}\\
\vspace{2mm}\hrule
Elliptic case of induced representations\\
map \(r(M)\): \(\left(\begin{array}{cc}\frac{d^{2}}{d^{2}+c^{2}}&-\frac{ d c}{d^{2}+c^{2}}\\\frac{ d c}{d^{2}+c^{2}}&\frac{d^{2}}{d^{2}+c^{2}}\end{array}\right)\)\\
map \(s^{-1}(M)\): \(\left(\begin{array}{cc}\frac{1}{d}&\frac{ c^{2} b+c+ d^{2} b}{d^{2}}\\0&\frac{d^{2}+c^{2}}{d}\end{array}\right)\)\\
character: \(\left(\begin{array}{cc}\frac{2  u d c+d^{2}+ u^{2} c^{2}}{2  u d c+ c^{2} v^{2}+d^{2}+ u^{2} c^{2}}&-\frac{ c v {(d+ u c)}}{2  u d c+ c^{2} v^{2}+d^{2}+ u^{2} c^{2}}\\\frac{ c v {(d+ u c)}}{2  u d c+ c^{2} v^{2}+d^{2}+ u^{2} c^{2}}&\frac{{(d+ u c)}^{2}}{2  u d c+ c^{2} v^{2}+d^{2}+ u^{2} c^{2}}\end{array}\right)\)\\
Moebius map: \(\left(\begin{array}{cc}\frac{ u c b+ d b+ u a d+ a c v^{2}+ u^{2} a c}{2  u d c+ c^{2} v^{2}+d^{2}+ u^{2} c^{2}}&\frac{ a d v- c v b}{2  u d c+ c^{2} v^{2}+d^{2}+ u^{2} c^{2}}\end{array}\right)\)\\
Moebius map is given by the imaginary unit: \textbf{true}\\
\vspace{2mm}\hrule
Parabolic (\(N^\prime\)) case of induced representations\\
map \(r(M)\): \(\left(\begin{array}{cc}1&0\\\frac{c}{d}&1\end{array}\right)\)\\
map \(s^{-1}(M)\): \(\left(\begin{array}{cc}\frac{1}{d}&b\\0&d\end{array}\right)\)\\
character: \(\left(\begin{array}{cc}1&0\\\frac{ c v}{d+ u c}&1\end{array}\right)\)\\
Moebius map: \(\left(\begin{array}{cc}\frac{ u a+b}{d+ u c}&\frac{ a d v- c v b}{{(d+ u c)}^{2}}\end{array}\right)\)\\
Moebius map is given by the imaginary unit: \textbf{true}\\
\vspace{2mm}\hrule
Hyperbolic case of induced representations\\
map \(r(M)\): \(\left(\begin{array}{cc}\frac{d^{2}}{d^{2}-c^{2}}&\frac{ d c}{d^{2}-c^{2}}\\\frac{ d c}{d^{2}-c^{2}}&\frac{d^{2}}{d^{2}-c^{2}}\end{array}\right)\)\\
map \(s^{-1}(M)\): \(\left(\begin{array}{cc}\frac{1}{d}&-\frac{ c^{2} b+c- d^{2} b}{d^{2}}\\0&\frac{d^{2}-c^{2}}{d}\end{array}\right)\)\\
character: \(\left(\begin{array}{cc}\frac{2  u d c+d^{2}+ u^{2} c^{2}}{2  u d c- c^{2} v^{2}+d^{2}+ u^{2} c^{2}}&\frac{ c v {(d+ u c)}}{2  u d c- c^{2} v^{2}+d^{2}+ u^{2} c^{2}}\\\frac{ c v {(d+ u c)}}{2  u d c- c^{2} v^{2}+d^{2}+ u^{2} c^{2}}&\frac{{(d+ u c)}^{2}}{2  u d c- c^{2} v^{2}+d^{2}+ u^{2} c^{2}}\end{array}\right)\)\\
Moebius map: \(\left(\begin{array}{cc}\frac{ u c b+ d b+ u a d- a c v^{2}+ u^{2} a c}{2  u d c- c^{2} v^{2}+d^{2}+ u^{2} c^{2}}&\frac{ a d v- c v b}{2  u d c- c^{2} v^{2}+d^{2}+ u^{2} c^{2}}\end{array}\right)\)\\
Moebius map is given by the imaginary unit: \textbf{true}\\
\vspace{2mm}\hrule
}{}

\IfFileExists{parab-rotation.tex}{
\section{Program for Symbolic Calculations} 
\label{sec:progr-symb-calc}
This is a documentation for our symbolic calculations. You can obtain
the program itself from the
\href{http://arXiv.org/abs/math/yymmnnn}{source files} of this
paper~\cite{Kisil07a} at \url{arXiv.org}; \LaTeX\ compilation of it
will produces the file \texttt{parab-rotation.nw} in the current directory.
This is a \NoWEB~\cite{NoWEB} code of the program. It uses Clifford
algebra facilities~\cite{Kisil05b} of the \GiNaC\
library~\cite{GiNaC}.

\def\LA{\begingroup\maybehbox\bgroup\setupmodname\Rm$\langle$}\def\RA{$\rangle$\egroup\endgroup}\providecommand{\MM}{\kern.5pt\raisebox{.4ex}{\begin{math}\scriptscriptstyle-\kern-1pt-\end{math}}\kern.5pt}\providecommand{\PP}{\kern.5pt\raisebox{.4ex}{\begin{math}\scriptscriptstyle+\kern-1pt+\end{math}}\kern.5pt}\def\commopen{/\begin{math}\ast\,\end{math}}\def\commclose{\,\begin{math}\ast\end{math}\kern-.5pt/}\def\begcomm{\begingroup\maybehbox\bgroup\setupmodname}\def\endcomm{\egroup\endgroup}\nwfilename{parab-rotation.nw}\nwbegindocs{1}\nwdocspar
\nwenddocs{}\nwbegindocs{2}%
\nwenddocs{}\nwbegindocs{3}%
\nwenddocs{}\nwbegindocs{4}%
\nwenddocs{}\nwbegindocs{5}%
\nwenddocs{}\nwbegindocs{6}%
\nwenddocs{}\nwbegindocs{7}%
\nwenddocs{}\nwbegindocs{8}%
\nwenddocs{}\nwbegindocs{9}%
\nwenddocs{}\nwbegindocs{10}%
\nwenddocs{}\nwbegindocs{11}%
\nwenddocs{}\nwbegindocs{12}%
\nwenddocs{}\nwbegindocs{13}%
\nwenddocs{}\nwbegindocs{14}%
\nwenddocs{}\nwbegindocs{15}%
\nwenddocs{}\nwbegindocs{16}%
\nwenddocs{}\nwbegindocs{17}%
\nwenddocs{}\nwbegindocs{18}%
\nwdocspar
This piece of software is licensed under
\href{http://www.gnu.org/licenses/gpl.html}{GNU General Public
  License} (Version 3, 29 June 2007)~\cite{GNUGPL}.

\nwenddocs{}\nwbegindocs{19}\nwdocspar
\subsection{Class {\Tt{}\Rm{}{\bf{}dual\_number}\nwendquote}}
\label{sec:class-dual_number}

\nwenddocs{}\nwbegindocs{20}\nwdocspar
\subsubsection{Public Methods}
\label{sec:public-methods}

\nwenddocs{}\nwbegindocs{21}\nwdocspar
A dual number can be created simply by listing its two components.
\nwenddocs{}\nwbegincode{22}\sublabel{NWgqRer-1eKCCy-1}\nwmargintag{{\nwtagstyle{}\subpageref{NWgqRer-1eKCCy-1}}}\moddef{Public methods~{\nwtagstyle{}\subpageref{NWgqRer-1eKCCy-1}}}\endmoddef\Rm{}\nwstartdeflinemarkup\nwusesondefline{\\{NWgqRer-4GrAIY-2}}\nwprevnextdefs{\relax}{NWgqRer-1eKCCy-2}\nwenddeflinemarkup
        {\bf{}dual\_number}({\bf{}const} {\bf{}ex} & {\it{}a}, {\bf{}const} {\bf{}ex} & {\it{}b});

\nwalsodefined{\\{NWgqRer-1eKCCy-2}\\{NWgqRer-1eKCCy-3}\\{NWgqRer-1eKCCy-4}\\{NWgqRer-1eKCCy-5}\\{NWgqRer-1eKCCy-6}\\{NWgqRer-1eKCCy-7}\\{NWgqRer-1eKCCy-8}\\{NWgqRer-1eKCCy-9}}\nwused{\\{NWgqRer-4GrAIY-2}}\nwidentuses{\\{{\nwixident{dual{\_}number}}{dual:unnumber}}}\nwindexuse{\nwixident{dual{\_}number}}{dual:unnumber}{NWgqRer-1eKCCy-1}\nwendcode{}\nwbegindocs{23}Alternatively you can provide a \(1\times 2\) or \(2\times 1\)
matrix,  a list, another dual number {\Tt{}\Rm{}{\it{}P}\nwendquote} or a complex expression
with a non-zero imaginary part to give two components. 
If {\Tt{}\Rm{}{\it{}P}\nwendquote} does not have two components and is a real-valued
expression, it will be embedded into dual numbers with zero argument
and norm equal to {\Tt{}\Rm{}{\it{}P}\nwendquote}.
\nwenddocs{}\nwbegincode{24}\sublabel{NWgqRer-1eKCCy-2}\nwmargintag{{\nwtagstyle{}\subpageref{NWgqRer-1eKCCy-2}}}\moddef{Public methods~{\nwtagstyle{}\subpageref{NWgqRer-1eKCCy-1}}}\plusendmoddef\Rm{}\nwstartdeflinemarkup\nwusesondefline{\\{NWgqRer-4GrAIY-2}}\nwprevnextdefs{NWgqRer-1eKCCy-1}{NWgqRer-1eKCCy-3}\nwenddeflinemarkup
        {\bf{}dual\_number}({\bf{}const} {\bf{}ex} & {\it{}P});

\nwused{\\{NWgqRer-4GrAIY-2}}\nwidentuses{\\{{\nwixident{dual{\_}number}}{dual:unnumber}}}\nwindexuse{\nwixident{dual{\_}number}}{dual:unnumber}{NWgqRer-1eKCCy-2}\nwendcode{}\nwbegindocs{25} We can also obtain the module and argument of a {\Tt{}\Rm{}{\bf{}dual\_number}\nwendquote}.
\nwenddocs{}\nwbegincode{26}\sublabel{NWgqRer-1eKCCy-3}\nwmargintag{{\nwtagstyle{}\subpageref{NWgqRer-1eKCCy-3}}}\moddef{Public methods~{\nwtagstyle{}\subpageref{NWgqRer-1eKCCy-1}}}\plusendmoddef\Rm{}\nwstartdeflinemarkup\nwusesondefline{\\{NWgqRer-4GrAIY-2}}\nwprevnextdefs{NWgqRer-1eKCCy-2}{NWgqRer-1eKCCy-4}\nwenddeflinemarkup
        {\bf{}ex} {\it{}arg}() {\bf{}const};
        {\bf{}ex} {\it{}norm}() {\bf{}const};

\nwused{\\{NWgqRer-4GrAIY-2}}\nwendcode{}\nwbegindocs{27}We define the conjugate of a {\Tt{}\Rm{}{\bf{}dual\_number}\nwendquote} by \(\overline{u+\rmp
  v}=-u+\rmp v\).
\nwenddocs{}\nwbegincode{28}\sublabel{NWgqRer-1eKCCy-4}\nwmargintag{{\nwtagstyle{}\subpageref{NWgqRer-1eKCCy-4}}}\moddef{Public methods~{\nwtagstyle{}\subpageref{NWgqRer-1eKCCy-1}}}\plusendmoddef\Rm{}\nwstartdeflinemarkup\nwusesondefline{\\{NWgqRer-4GrAIY-2}}\nwprevnextdefs{NWgqRer-1eKCCy-3}{NWgqRer-1eKCCy-5}\nwenddeflinemarkup
        {\bf{}ex} {\it{}conjugate}() {\bf{}const} {\nwlbrace} {\bf{}return} {\bf{}dual\_number}(-{\it{}u\_comp}, {\it{}v\_comp}); {\nwrbrace}

\nwused{\\{NWgqRer-4GrAIY-2}}\nwidentuses{\\{{\nwixident{dual{\_}number}}{dual:unnumber}}}\nwindexuse{\nwixident{dual{\_}number}}{dual:unnumber}{NWgqRer-1eKCCy-4}\nwendcode{}\nwbegindocs{29} Real part may be defined out of the formula
\(\Re(u,v)=\frac{1}{2}((u,v)+\overline{(u,v)})\), but it simply
reduces to the value of norm for the {\Tt{}\Rm{}{\bf{}dual\_number}\nwendquote}.
\nwenddocs{}\nwbegincode{30}\sublabel{NWgqRer-1eKCCy-5}\nwmargintag{{\nwtagstyle{}\subpageref{NWgqRer-1eKCCy-5}}}\moddef{Public methods~{\nwtagstyle{}\subpageref{NWgqRer-1eKCCy-1}}}\plusendmoddef\Rm{}\nwstartdeflinemarkup\nwusesondefline{\\{NWgqRer-4GrAIY-2}}\nwprevnextdefs{NWgqRer-1eKCCy-4}{NWgqRer-1eKCCy-6}\nwenddeflinemarkup
        {\bf{}ex} {\it{}real\_part}() {\bf{}const};
        {\bf{}ex} {\it{}imag\_part}() {\bf{}const}; 

\nwused{\\{NWgqRer-4GrAIY-2}}\nwendcode{}\nwbegindocs{31}Negative of a {\Tt{}\Rm{}{\bf{}dual\_number}\nwendquote} and its power.
\nwenddocs{}\nwbegincode{32}\sublabel{NWgqRer-1eKCCy-6}\nwmargintag{{\nwtagstyle{}\subpageref{NWgqRer-1eKCCy-6}}}\moddef{Public methods~{\nwtagstyle{}\subpageref{NWgqRer-1eKCCy-1}}}\plusendmoddef\Rm{}\nwstartdeflinemarkup\nwusesondefline{\\{NWgqRer-4GrAIY-2}}\nwprevnextdefs{NWgqRer-1eKCCy-5}{NWgqRer-1eKCCy-7}\nwenddeflinemarkup
        {\bf{}dual\_number} {\it{}neg}() {\bf{}const} 
                {\nwlbrace} {\bf{}return} {\bf{}dual\_number}({\it{}u\_comp}, -{\it{}v\_comp} + ({\it{}is\_subgroup\_N}? 2\begin{math}\ast\end{math}{\it{}pow}({\it{}u\_comp}, 2) : -2)); {\nwrbrace}
        {\bf{}dual\_number} {\it{}power}({\bf{}const} {\bf{}ex} & {\it{}e}) {\bf{}const};

\nwused{\\{NWgqRer-4GrAIY-2}}\nwidentuses{\\{{\nwixident{dual{\_}number}}{dual:unnumber}}}\nwindexuse{\nwixident{dual{\_}number}}{dual:unnumber}{NWgqRer-1eKCCy-6}\nwendcode{}\nwbegindocs{33} We can also convert a {\Tt{}\Rm{}{\bf{}dual\_number}\nwendquote} to a {\Tt{}\Rm{}{\bf{}matrix}\nwendquote}.
\nwenddocs{}\nwbegincode{34}\sublabel{NWgqRer-1eKCCy-7}\nwmargintag{{\nwtagstyle{}\subpageref{NWgqRer-1eKCCy-7}}}\moddef{Public methods~{\nwtagstyle{}\subpageref{NWgqRer-1eKCCy-1}}}\plusendmoddef\Rm{}\nwstartdeflinemarkup\nwusesondefline{\\{NWgqRer-4GrAIY-2}}\nwprevnextdefs{NWgqRer-1eKCCy-6}{NWgqRer-1eKCCy-8}\nwenddeflinemarkup
        {\bf{}matrix} {\it{}to\_matrix}() {\bf{}const} {\nwlbrace} {\bf{}return} {\bf{}matrix} (1, 2, {\bf{}lst}({\it{}u\_comp}, {\it{}v\_comp})); {\nwrbrace}

\nwused{\\{NWgqRer-4GrAIY-2}}\nwendcode{}\nwbegindocs{35}We define the rule for parabolic norm of a sum, see~\eqref{eq:p-add-norm-exotic}.
\nwenddocs{}\nwbegincode{36}\sublabel{NWgqRer-1eKCCy-8}\nwmargintag{{\nwtagstyle{}\subpageref{NWgqRer-1eKCCy-8}}}\moddef{Public methods~{\nwtagstyle{}\subpageref{NWgqRer-1eKCCy-1}}}\plusendmoddef\Rm{}\nwstartdeflinemarkup\nwusesondefline{\\{NWgqRer-4GrAIY-2}}\nwprevnextdefs{NWgqRer-1eKCCy-7}{NWgqRer-1eKCCy-9}\nwenddeflinemarkup
        {\bf{}ex} {\it{}add\_norms}({\bf{}const} {\bf{}dual\_number} & {\it{}P}) {\bf{}const} {\nwlbrace} {\bf{}return} ({\it{}norm}()+{\it{}P}.{\it{}norm}()).{\it{}normal}(); {\nwrbrace}

\nwused{\\{NWgqRer-4GrAIY-2}}\nwidentuses{\\{{\nwixident{dual{\_}number}}{dual:unnumber}}}\nwindexuse{\nwixident{dual{\_}number}}{dual:unnumber}{NWgqRer-1eKCCy-8}\nwendcode{}\nwbegindocs{37}Algebraic operations are defined for {\Tt{}\Rm{}{\bf{}dual\_number}\nwendquote}s in a way
described in \S~\ref{sec:exsotic-form}. The standard \CPP\ operators
{\Tt{}\Rm{}+\nwendquote}, {\Tt{}\Rm{}-\nwendquote}, {\Tt{}\Rm{}\begin{math}\ast\end{math}\nwendquote}, {\Tt{}\Rm{}\begin{math}\div\end{math}\nwendquote} will be overloaded later in order to permit
natural expressions with {\Tt{}\Rm{}{\bf{}dual\_number}\nwendquote}s.
\nwenddocs{}\nwbegincode{38}\sublabel{NWgqRer-1eKCCy-9}\nwmargintag{{\nwtagstyle{}\subpageref{NWgqRer-1eKCCy-9}}}\moddef{Public methods~{\nwtagstyle{}\subpageref{NWgqRer-1eKCCy-1}}}\plusendmoddef\Rm{}\nwstartdeflinemarkup\nwusesondefline{\\{NWgqRer-4GrAIY-2}}\nwprevnextdefs{NWgqRer-1eKCCy-8}{\relax}\nwenddeflinemarkup
        {\bf{}dual\_number} {\it{}add}({\bf{}const} {\bf{}dual\_number} & {\it{}a}) {\bf{}const};
        {\bf{}dual\_number} {\it{}sub}({\bf{}const} {\bf{}dual\_number} & {\it{}a}) {\bf{}const} {\nwlbrace} {\bf{}return} {\it{}add}({\it{}a}.{\it{}neg}()); {\nwrbrace}
        {\bf{}dual\_number} {\it{}mul}({\bf{}const} {\bf{}dual\_number} & {\it{}a}) {\bf{}const};
        {\bf{}dual\_number} {\it{}mul}({\bf{}const} {\bf{}ex} & {\it{}a}) {\bf{}const} {\nwlbrace} {\bf{}return} {\it{}mul}({\bf{}dual\_number}({\it{}a})); {\nwrbrace}

\nwused{\\{NWgqRer-4GrAIY-2}}\nwidentuses{\\{{\nwixident{dual{\_}number}}{dual:unnumber}}}\nwindexuse{\nwixident{dual{\_}number}}{dual:unnumber}{NWgqRer-1eKCCy-9}\nwendcode{}\nwbegindocs{39}\nwdocspar
\subsection{Algebraic Subroutines}
\label{sec:algebr-subr}
We need a couple of global variables which help to write uniformly
algebraic rules for both cases of subgroups \(N\) and \(N^\prime\).

\nwenddocs{}\nwbegindocs{40}Firstly, we need to consider separately cases of subgroup \(N\) and
\(N^\prime\), the following global variable keeps track on it.
\nwenddocs{}\nwbegincode{41}\sublabel{NWgqRer-4EWpKB-1}\nwmargintag{{\nwtagstyle{}\subpageref{NWgqRer-4EWpKB-1}}}\moddef{N-Nprime separation~{\nwtagstyle{}\subpageref{NWgqRer-4EWpKB-1}}}\endmoddef\Rm{}\nwstartdeflinemarkup\nwusesondefline{\\{NWgqRer-36Ytqo-1}}\nwprevnextdefs{\relax}{NWgqRer-4EWpKB-2}\nwenddeflinemarkup
{\bf{}bool} {\it{}is\_subgroup\_N};

\nwalsodefined{\\{NWgqRer-4EWpKB-2}\\{NWgqRer-4EWpKB-3}}\nwused{\\{NWgqRer-36Ytqo-1}}\nwendcode{}\nwbegindocs{42}In the case of the subgroup \(N^\prime\) the reference point lies at infinity, the
following {\Tt{}\Rm{}{\bf{}realsymbol}\nwendquote} variable represents it in the symbolic
calculations.
\nwenddocs{}\nwbegincode{43}\sublabel{NWgqRer-4EWpKB-2}\nwmargintag{{\nwtagstyle{}\subpageref{NWgqRer-4EWpKB-2}}}\moddef{N-Nprime separation~{\nwtagstyle{}\subpageref{NWgqRer-4EWpKB-1}}}\plusendmoddef\Rm{}\nwstartdeflinemarkup\nwusesondefline{\\{NWgqRer-36Ytqo-1}}\nwprevnextdefs{NWgqRer-4EWpKB-1}{NWgqRer-4EWpKB-3}\nwenddeflinemarkup
{\bf{}realsymbol} {\it{}Inf}({\tt{}"Inf"}, {\tt{}"{\char92}{\char92}infty"});

\nwused{\\{NWgqRer-36Ytqo-1}}\nwendcode{}\nwbegindocs{44} We define the ``zero angle'': for the subgroup \(N\) it is \(0\), for
 \(N^\prime\)---\(\infty\).
\nwenddocs{}\nwbegincode{45}\sublabel{NWgqRer-4EWpKB-3}\nwmargintag{{\nwtagstyle{}\subpageref{NWgqRer-4EWpKB-3}}}\moddef{N-Nprime separation~{\nwtagstyle{}\subpageref{NWgqRer-4EWpKB-1}}}\plusendmoddef\Rm{}\nwstartdeflinemarkup\nwusesondefline{\\{NWgqRer-36Ytqo-1}}\nwprevnextdefs{NWgqRer-4EWpKB-2}{\relax}\nwenddeflinemarkup
{\bf{}\char35{}define}{\tt{} Arg0 (is\_subgroup\_N ? ex(0) : ex(Inf))}\nwindexdefn{\nwixident{Arg0}}{Arg0}{NWgqRer-4EWpKB-3}

\nwused{\\{NWgqRer-36Ytqo-1}}\nwidentdefs{\\{{\nwixident{Arg0}}{Arg0}}}\nwendcode{}\nwbegindocs{46}Here is the set of algebraic procedures representing definitions
made in this paper.

\nwenddocs{}\nwbegindocs{47}\nwdocspar
\subsubsection[Argument and u]{Argument and $u$}
\label{sec:argument-u}
In the case of \(N\) the value of \(u\) is simply the argument.
\nwenddocs{}\nwbegincode{48}\sublabel{NWgqRer-3jzoaG-1}\nwmargintag{{\nwtagstyle{}\subpageref{NWgqRer-3jzoaG-1}}}\moddef{Algebraic procedures~{\nwtagstyle{}\subpageref{NWgqRer-3jzoaG-1}}}\endmoddef\Rm{}\nwstartdeflinemarkup\nwusesondefline{\\{NWgqRer-3HhVai-1}}\nwprevnextdefs{\relax}{NWgqRer-3jzoaG-2}\nwenddeflinemarkup
{\bf{}ex} {\it{}u\_from\_arg}({\bf{}const} {\bf{}ex} & {\it{}a}) {\nwlbrace}
        {\bf{}if} ({\it{}is\_subgroup\_N})
                {\bf{}return} {\it{}a};

\nwalsodefined{\\{NWgqRer-3jzoaG-2}\\{NWgqRer-3jzoaG-3}\\{NWgqRer-3jzoaG-4}\\{NWgqRer-3jzoaG-5}\\{NWgqRer-3jzoaG-6}\\{NWgqRer-3jzoaG-7}}\nwused{\\{NWgqRer-3HhVai-1}}\nwendcode{}\nwbegindocs{49}In the case of \(N^\prime\) the value of \(u\) is the inverse to the argument,
and we need to treat properly the case of zero\ldots
\nwenddocs{}\nwbegincode{50}\sublabel{NWgqRer-3jzoaG-2}\nwmargintag{{\nwtagstyle{}\subpageref{NWgqRer-3jzoaG-2}}}\moddef{Algebraic procedures~{\nwtagstyle{}\subpageref{NWgqRer-3jzoaG-1}}}\plusendmoddef\Rm{}\nwstartdeflinemarkup\nwusesondefline{\\{NWgqRer-3HhVai-1}}\nwprevnextdefs{NWgqRer-3jzoaG-1}{NWgqRer-3jzoaG-3}\nwenddeflinemarkup
        {\bf{}else} {\nwlbrace}
                {\bf{}if} ({\it{}a}.{\it{}normal}().{\it{}is\_zero}())
                        {\bf{}return} {\it{}Inf};

\nwused{\\{NWgqRer-3HhVai-1}}\nwendcode{}\nwbegindocs{51}\ldots and \(\infty\).  We try to replace \(\frac{1}{\infty}\) by \(0\).
\nwenddocs{}\nwbegincode{52}\sublabel{NWgqRer-3jzoaG-3}\nwmargintag{{\nwtagstyle{}\subpageref{NWgqRer-3jzoaG-3}}}\moddef{Algebraic procedures~{\nwtagstyle{}\subpageref{NWgqRer-3jzoaG-1}}}\plusendmoddef\Rm{}\nwstartdeflinemarkup\nwusesondefline{\\{NWgqRer-3HhVai-1}}\nwprevnextdefs{NWgqRer-3jzoaG-2}{NWgqRer-3jzoaG-4}\nwenddeflinemarkup
                {\bf{}else}
                        {\bf{}try} {\nwlbrace}
                                {\bf{}realsymbol} {\it{}t}({\tt{}"t"});
                                {\bf{}return} {\it{}pow}({\it{}a}.{\it{}subs}({\it{}Inf} \begin{math}\equiv\end{math} {\it{}pow}({\it{}t}, -1)), -1).{\it{}normal}().{\it{}subs}({\it{}t} \begin{math}\equiv\end{math} 0).{\it{}normal}();
                        {\nwrbrace} {\bf{}catch} ({\it{}std}::{\it{}exception} &{\it{}p}) {\nwlbrace}
                                {\bf{}return} {\it{}pow}({\it{}a}, -1);
                        {\nwrbrace}
        {\nwrbrace}
{\nwrbrace}

\nwused{\\{NWgqRer-3HhVai-1}}\nwendcode{}\nwbegindocs{53}\nwdocspar
\subsubsection{Argument of a point}
\label{sec:argument-point}
The opposite task (finding argument of a point) is solved similarly.
\nwenddocs{}\nwbegincode{54}\sublabel{NWgqRer-4dxCds-1}\nwmargintag{{\nwtagstyle{}\subpageref{NWgqRer-4dxCds-1}}}\moddef{Dual number class further implementation~{\nwtagstyle{}\subpageref{NWgqRer-4dxCds-1}}}\endmoddef\Rm{}\nwstartdeflinemarkup\nwusesondefline{\\{NWgqRer-3HhVai-1}}\nwprevnextdefs{\relax}{NWgqRer-4dxCds-2}\nwenddeflinemarkup
{\bf{}ex} {\bf{}dual\_number}::{\it{}arg}() {\bf{}const} {\nwlbrace}
        {\bf{}if} ({\it{}is\_subgroup\_N})
                {\bf{}return} {\it{}u\_comp};

\nwalsodefined{\\{NWgqRer-4dxCds-2}\\{NWgqRer-4dxCds-3}\\{NWgqRer-4dxCds-4}\\{NWgqRer-4dxCds-5}\\{NWgqRer-4dxCds-6}\\{NWgqRer-4dxCds-7}\\{NWgqRer-4dxCds-8}\\{NWgqRer-4dxCds-9}}\nwused{\\{NWgqRer-3HhVai-1}}\nwidentuses{\\{{\nwixident{dual{\_}number}}{dual:unnumber}}}\nwindexuse{\nwixident{dual{\_}number}}{dual:unnumber}{NWgqRer-4dxCds-1}\nwendcode{}\nwbegindocs{55}Again in the case of \(N^\prime\) we need to consider cases of \(0\)\ldots
\nwenddocs{}\nwbegincode{56}\sublabel{NWgqRer-4dxCds-2}\nwmargintag{{\nwtagstyle{}\subpageref{NWgqRer-4dxCds-2}}}\moddef{Dual number class further implementation~{\nwtagstyle{}\subpageref{NWgqRer-4dxCds-1}}}\plusendmoddef\Rm{}\nwstartdeflinemarkup\nwusesondefline{\\{NWgqRer-3HhVai-1}}\nwprevnextdefs{NWgqRer-4dxCds-1}{NWgqRer-4dxCds-3}\nwenddeflinemarkup
        {\bf{}else} {\nwlbrace}
                {\bf{}if} ({\it{}u\_comp}.{\it{}normal}().{\it{}is\_zero}())
                        {\bf{}return} {\it{}Inf};

\nwused{\\{NWgqRer-3HhVai-1}}\nwendcode{}\nwbegindocs{57}\ldots and \(\infty\). We try to replace \(\frac{1}{\infty}\) by \(0\).
\nwenddocs{}\nwbegincode{58}\sublabel{NWgqRer-4dxCds-3}\nwmargintag{{\nwtagstyle{}\subpageref{NWgqRer-4dxCds-3}}}\moddef{Dual number class further implementation~{\nwtagstyle{}\subpageref{NWgqRer-4dxCds-1}}}\plusendmoddef\Rm{}\nwstartdeflinemarkup\nwusesondefline{\\{NWgqRer-3HhVai-1}}\nwprevnextdefs{NWgqRer-4dxCds-2}{NWgqRer-4dxCds-4}\nwenddeflinemarkup
                {\bf{}else} 
                        {\bf{}try} {\nwlbrace}
                                {\bf{}realsymbol} {\it{}t}({\tt{}"t"});
                                {\bf{}return} {\it{}pow}({\it{}u\_comp}.{\it{}subs}({\it{}Inf} \begin{math}\equiv\end{math} {\it{}pow}({\it{}t}, -1)), -1).{\it{}normal}().{\it{}subs}({\it{}t} \begin{math}\equiv\end{math} 0).{\it{}normal}();
                        {\nwrbrace} {\bf{}catch} ({\it{}std}::{\it{}exception} &{\it{}p}) {\nwlbrace}
                                {\bf{}return} {\it{}pow}({\it{}u\_comp}, -1);
                        {\nwrbrace}
        {\nwrbrace}
{\nwrbrace}

\nwused{\\{NWgqRer-3HhVai-1}}\nwendcode{}\nwbegindocs{59}\nwdocspar

\subsubsection{Norm}
\label{sec:norm}

The corresponding value of the parabolic norm is calculated by the
formulae~\eqref{eq:parab-norm}.
\nwenddocs{}\nwbegincode{60}\sublabel{NWgqRer-4dxCds-4}\nwmargintag{{\nwtagstyle{}\subpageref{NWgqRer-4dxCds-4}}}\moddef{Dual number class further implementation~{\nwtagstyle{}\subpageref{NWgqRer-4dxCds-1}}}\plusendmoddef\Rm{}\nwstartdeflinemarkup\nwusesondefline{\\{NWgqRer-3HhVai-1}}\nwprevnextdefs{NWgqRer-4dxCds-3}{NWgqRer-4dxCds-5}\nwenddeflinemarkup
{\bf{}ex} {\bf{}dual\_number}::{\it{}norm}() {\bf{}const} {\nwlbrace}
        {\bf{}if} ({\it{}is\_subgroup\_N})
                {\bf{}return} {\it{}pow}({\it{}u\_comp}, 2)-{\it{}v\_comp};

\nwused{\\{NWgqRer-3HhVai-1}}\nwidentuses{\\{{\nwixident{dual{\_}number}}{dual:unnumber}}}\nwindexuse{\nwixident{dual{\_}number}}{dual:unnumber}{NWgqRer-4dxCds-4}\nwendcode{}\nwbegindocs{61}The case of subgroup \(N^\prime\) require treatment of infinity.
\nwenddocs{}\nwbegincode{62}\sublabel{NWgqRer-4dxCds-5}\nwmargintag{{\nwtagstyle{}\subpageref{NWgqRer-4dxCds-5}}}\moddef{Dual number class further implementation~{\nwtagstyle{}\subpageref{NWgqRer-4dxCds-1}}}\plusendmoddef\Rm{}\nwstartdeflinemarkup\nwusesondefline{\\{NWgqRer-3HhVai-1}}\nwprevnextdefs{NWgqRer-4dxCds-4}{NWgqRer-4dxCds-6}\nwenddeflinemarkup
        {\bf{}else} {\nwlbrace}
                {\bf{}if} ({\it{}u\_comp}.{\it{}is\_zero}()) {\nwlbrace}
                        {\bf{}if} (({\it{}v\_comp}+1).{\it{}is\_zero}())
                                {\bf{}return} 1;
                        {\bf{}else}
                                {\bf{}return} {\it{}pow}({\it{}Inf}, 2)\begin{math}\div\end{math}({\it{}v\_comp}+1);
                {\nwrbrace} {\bf{}else}
                        {\bf{}return} ({\it{}pow}({\it{}u\_comp}, 2)\begin{math}\div\end{math}({\it{}v\_comp}+1)).{\it{}normal}();
        {\nwrbrace}
{\nwrbrace}

\nwused{\\{NWgqRer-3HhVai-1}}\nwendcode{}\nwbegindocs{63}\nwdocspar
\subsubsection[The value of v from the argument and norm]{The value of
  \(v\) from the argument and norm} 
\label{sec:value-v-from}
We oftenly need to find values of \(v\) such that for a given
value of argument \(A\) point \((A, v)\) will have a given norm.
\nwenddocs{}\nwbegincode{64}\sublabel{NWgqRer-3jzoaG-4}\nwmargintag{{\nwtagstyle{}\subpageref{NWgqRer-3jzoaG-4}}}\moddef{Algebraic procedures~{\nwtagstyle{}\subpageref{NWgqRer-3jzoaG-1}}}\plusendmoddef\Rm{}\nwstartdeflinemarkup\nwusesondefline{\\{NWgqRer-3HhVai-1}}\nwprevnextdefs{NWgqRer-3jzoaG-3}{NWgqRer-3jzoaG-5}\nwenddeflinemarkup
{\bf{}ex} {\it{}v\_from\_norm}({\bf{}const} {\bf{}ex} & {\it{}u}, {\bf{}const} {\bf{}ex} & {\it{}n}) {\nwlbrace}
        {\bf{}realsymbol} {\it{}l}({\tt{}"l"});
        {\bf{}if} ({\it{}is\_subgroup\_N})
                {\bf{}return} {\it{}lsolve}({\bf{}dual\_number}({\it{}u}, {\it{}l}).{\it{}norm}() \begin{math}\equiv\end{math} {\it{}n}, {\it{}l}).{\it{}normal}();
        {\bf{}else}
                {\bf{}return} {\it{}lsolve}({\it{}pow}({\bf{}dual\_number}({\it{}u}, {\it{}l}).{\it{}norm}(), -1) \begin{math}\equiv\end{math} {\it{}pow}({\it{}n}, -1), {\it{}l}).{\it{}normal}();
{\nwrbrace}

\nwused{\\{NWgqRer-3HhVai-1}}\nwidentuses{\\{{\nwixident{dual{\_}number}}{dual:unnumber}}}\nwindexuse{\nwixident{dual{\_}number}}{dual:unnumber}{NWgqRer-3jzoaG-4}\nwendcode{}\nwbegindocs{65}\nwdocspar
\nwenddocs{}\nwbegincode{66}\sublabel{NWgqRer-3jzoaG-5}\nwmargintag{{\nwtagstyle{}\subpageref{NWgqRer-3jzoaG-5}}}\moddef{Algebraic procedures~{\nwtagstyle{}\subpageref{NWgqRer-3jzoaG-1}}}\plusendmoddef\Rm{}\nwstartdeflinemarkup\nwusesondefline{\\{NWgqRer-3HhVai-1}}\nwprevnextdefs{NWgqRer-3jzoaG-4}{NWgqRer-3jzoaG-6}\nwenddeflinemarkup
{\bf{}dual\_number} {\it{}zero\_dual\_number}() {\nwlbrace}
        {\bf{}return} ({\it{}is\_subgroup\_N} ? {\bf{}dual\_number}(0, 0) : {\bf{}dual\_number}({\it{}Inf}, -1));
{\nwrbrace}

\nwused{\\{NWgqRer-3HhVai-1}}\nwidentuses{\\{{\nwixident{dual{\_}number}}{dual:unnumber}}}\nwindexuse{\nwixident{dual{\_}number}}{dual:unnumber}{NWgqRer-3jzoaG-5}\nwendcode{}\nwbegindocs{67}\nwdocspar
\nwenddocs{}\nwbegincode{68}\sublabel{NWgqRer-4dxCds-6}\nwmargintag{{\nwtagstyle{}\subpageref{NWgqRer-4dxCds-6}}}\moddef{Dual number class further implementation~{\nwtagstyle{}\subpageref{NWgqRer-4dxCds-1}}}\plusendmoddef\Rm{}\nwstartdeflinemarkup\nwusesondefline{\\{NWgqRer-3HhVai-1}}\nwprevnextdefs{NWgqRer-4dxCds-5}{NWgqRer-4dxCds-7}\nwenddeflinemarkup
{\bf{}bool} {\bf{}dual\_number}::{\it{}is\_zero}() {\bf{}const} {\nwlbrace}
        {\bf{}return} {\it{}is\_equal}({\it{}zero\_dual\_number}());
{\nwrbrace}

\nwused{\\{NWgqRer-3HhVai-1}}\nwidentuses{\\{{\nwixident{dual{\_}number}}{dual:unnumber}}}\nwindexuse{\nwixident{dual{\_}number}}{dual:unnumber}{NWgqRer-4dxCds-6}\nwendcode{}\nwbegindocs{69}\nwdocspar
\subsubsection{Real and Imaginary Parts}
\label{sec:real-imag-parts}

See \S~\ref{sec:real-imaginary-parts} for a discussion of the real and
imaginary parts of dual numbers.
\nwenddocs{}\nwbegincode{70}\sublabel{NWgqRer-3jzoaG-6}\nwmargintag{{\nwtagstyle{}\subpageref{NWgqRer-3jzoaG-6}}}\moddef{Algebraic procedures~{\nwtagstyle{}\subpageref{NWgqRer-3jzoaG-1}}}\plusendmoddef\Rm{}\nwstartdeflinemarkup\nwusesondefline{\\{NWgqRer-3HhVai-1}}\nwprevnextdefs{NWgqRer-3jzoaG-5}{NWgqRer-3jzoaG-7}\nwenddeflinemarkup
{\bf{}ex} {\bf{}dual\_number}::{\it{}real\_part}() {\bf{}const} {\nwlbrace}
        {\bf{}return} {\it{}dn\_from\_arg\_mod}(0, (1-{\it{}arg}())\begin{math}\ast\end{math}{\it{}norm}());
{\nwrbrace}

{\bf{}ex} {\bf{}dual\_number}::{\it{}imag\_part}() {\bf{}const} {\nwlbrace}
        {\bf{}return} {\it{}dn\_from\_arg\_mod}(1, {\it{}arg}()\begin{math}\ast\end{math}{\it{}norm}());
{\nwrbrace}

\nwused{\\{NWgqRer-3HhVai-1}}\nwidentuses{\\{{\nwixident{dual{\_}number}}{dual:unnumber}}}\nwindexuse{\nwixident{dual{\_}number}}{dual:unnumber}{NWgqRer-3jzoaG-6}\nwendcode{}\nwbegindocs{71}\nwdocspar
\subsubsection{Product of Two Points}
\label{sec:product-two-points}
We define now the product of two points according to the
Definition~\ref{de:product}. We also include a multiplication by a
scalar: if a factor is a scalar it is replaced by a vector with the
zero argument and norm equal to the scalar.
\nwenddocs{}\nwbegincode{72}\sublabel{NWgqRer-4dxCds-7}\nwmargintag{{\nwtagstyle{}\subpageref{NWgqRer-4dxCds-7}}}\moddef{Dual number class further implementation~{\nwtagstyle{}\subpageref{NWgqRer-4dxCds-1}}}\plusendmoddef\Rm{}\nwstartdeflinemarkup\nwusesondefline{\\{NWgqRer-3HhVai-1}}\nwprevnextdefs{NWgqRer-4dxCds-6}{NWgqRer-4dxCds-8}\nwenddeflinemarkup
{\bf{}dual\_number} {\bf{}dual\_number}::{\it{}mul}({\bf{}const} {\bf{}dual\_number} & {\it{}P}) {\bf{}const} {\nwlbrace}
        {\bf{}ex} {\it{}u}={\it{}u\_from\_arg}({\it{}arg}()+{\it{}P}.{\it{}arg}()).{\it{}normal}();
        {\bf{}return} {\bf{}dual\_number}({\it{}u}, {\it{}v\_from\_norm}({\it{}u}, {\it{}norm}()\begin{math}\ast\end{math}{\it{}P}.{\it{}norm}()));
{\nwrbrace}

\nwused{\\{NWgqRer-3HhVai-1}}\nwidentuses{\\{{\nwixident{dual{\_}number}}{dual:unnumber}}}\nwindexuse{\nwixident{dual{\_}number}}{dual:unnumber}{NWgqRer-4dxCds-7}\nwendcode{}\nwbegindocs{73}\nwdocspar
\subsubsection{Vector Addition of Two Points}
\label{sec:vector-addition-two}

\nwenddocs{}\nwbegindocs{74} The sum is calculated from the expression~\eqref{eq:p-add-arg-exotic}. 
\nwenddocs{}\nwbegincode{75}\sublabel{NWgqRer-4dxCds-8}\nwmargintag{{\nwtagstyle{}\subpageref{NWgqRer-4dxCds-8}}}\moddef{Dual number class further implementation~{\nwtagstyle{}\subpageref{NWgqRer-4dxCds-1}}}\plusendmoddef\Rm{}\nwstartdeflinemarkup\nwusesondefline{\\{NWgqRer-3HhVai-1}}\nwprevnextdefs{NWgqRer-4dxCds-7}{NWgqRer-4dxCds-9}\nwenddeflinemarkup
{\bf{}dual\_number} {\bf{}dual\_number}::{\it{}add}({\bf{}const} {\bf{}dual\_number} & {\it{}a}) {\bf{}const} {\nwlbrace}
        {\bf{}ex} {\it{}norms} = {\it{}add\_norms}({\it{}a});
        {\bf{}if} ({\it{}norms}.{\it{}normal}().{\it{}is\_zero}())
                {\bf{}return} {\it{}zero\_dual\_number}();
        {\bf{}else} {\nwlbrace}
                {\bf{}ex} {\it{}us}={\it{}u\_from\_arg}(({\it{}arg}()\begin{math}\ast\end{math}{\it{}norm}()+{\it{}a}.{\it{}arg}()\begin{math}\ast\end{math}{\it{}a}.{\it{}norm}()) \begin{math}\div\end{math} {\it{}norms}).{\it{}normal}();
                {\bf{}return} {\bf{}dual\_number}({\it{}us}, {\it{}v\_from\_norm}({\it{}us}, {\it{}norms}));
        {\nwrbrace}
{\nwrbrace}

\nwused{\\{NWgqRer-3HhVai-1}}\nwidentuses{\\{{\nwixident{dual{\_}number}}{dual:unnumber}}}\nwindexuse{\nwixident{dual{\_}number}}{dual:unnumber}{NWgqRer-4dxCds-8}\nwendcode{}\nwbegindocs{76}\nwdocspar
\nwenddocs{}\nwbegincode{77}\sublabel{NWgqRer-3jzoaG-7}\nwmargintag{{\nwtagstyle{}\subpageref{NWgqRer-3jzoaG-7}}}\moddef{Algebraic procedures~{\nwtagstyle{}\subpageref{NWgqRer-3jzoaG-1}}}\plusendmoddef\Rm{}\nwstartdeflinemarkup\nwusesondefline{\\{NWgqRer-3HhVai-1}}\nwprevnextdefs{NWgqRer-3jzoaG-6}{\relax}\nwenddeflinemarkup
{\bf{}dual\_number} {\it{}dn\_from\_arg\_mod}({\bf{}const} {\bf{}ex} & {\it{}a}, {\bf{}const} {\bf{}ex} & {\it{}n}) {\nwlbrace}
        {\bf{}ex} {\it{}us}={\it{}u\_from\_arg}({\it{}a}).{\it{}normal}();
        {\bf{}return} {\bf{}dual\_number}({\it{}us}, {\it{}v\_from\_norm}({\it{}us}, {\it{}n}));
{\nwrbrace}

\nwused{\\{NWgqRer-3HhVai-1}}\nwidentuses{\\{{\nwixident{dual{\_}number}}{dual:unnumber}}}\nwindexuse{\nwixident{dual{\_}number}}{dual:unnumber}{NWgqRer-3jzoaG-7}\nwendcode{}\nwbegindocs{78} De Moivre's Identity:
\nwenddocs{}\nwbegincode{79}\sublabel{NWgqRer-4dxCds-9}\nwmargintag{{\nwtagstyle{}\subpageref{NWgqRer-4dxCds-9}}}\moddef{Dual number class further implementation~{\nwtagstyle{}\subpageref{NWgqRer-4dxCds-1}}}\plusendmoddef\Rm{}\nwstartdeflinemarkup\nwusesondefline{\\{NWgqRer-3HhVai-1}}\nwprevnextdefs{NWgqRer-4dxCds-8}{\relax}\nwenddeflinemarkup
{\bf{}dual\_number} {\bf{}dual\_number}::{\it{}power}({\bf{}const} {\bf{}ex} & {\it{}e}) {\bf{}const} {\nwlbrace}
        {\bf{}return} {\it{}dn\_from\_arg\_mod}({\it{}arg}()\begin{math}\ast\end{math}{\it{}e}, {\it{}pow}({\it{}norm}(), {\it{}e}));
{\nwrbrace}

\nwused{\\{NWgqRer-3HhVai-1}}\nwidentuses{\\{{\nwixident{dual{\_}number}}{dual:unnumber}}}\nwindexuse{\nwixident{dual{\_}number}}{dual:unnumber}{NWgqRer-4dxCds-9}\nwendcode{}\nwbegindocs{80}All algebraic routines are defined now.

\nwenddocs{}\nwbegindocs{81}\nwdocspar
\subsection{Calculation and Tests}
This Subsection contains code for calculation of various
expression. See~\cite{Kisil04c} or \GiNaC info for usage of Clifford
algebra functions.

\nwenddocs{}\nwbegindocs{82}\subsubsection{Calculation of Expressions}
\label{sec:calc-expr}
Firstly, we output the expression of the Cayley transform for a generic
element from subgroups \(N\) and \(N^\prime\).
\nwenddocs{}\nwbegincode{83}\sublabel{NWgqRer-1Pc9Jw-1}\nwmargintag{{\nwtagstyle{}\subpageref{NWgqRer-1Pc9Jw-1}}}\moddef{Show expressions~{\nwtagstyle{}\subpageref{NWgqRer-1Pc9Jw-1}}}\endmoddef\Rm{}\nwstartdeflinemarkup\nwusesondefline{\\{NWgqRer-nXe8t-4}}\nwprevnextdefs{\relax}{NWgqRer-1Pc9Jw-2}\nwenddeflinemarkup
{\bf{}ex} {\it{}XC}={\it{}canonicalize\_clifford}(({\it{}TC}\begin{math}\ast\end{math}{\it{}X}\begin{math}\ast\end{math}{\it{}TCI}).{\it{}evalm}());
{\it{}formula\_out}({\tt{}"Cayley of the matrix x: "}, {\it{}XC}.{\it{}subs}({\it{}sign}\begin{math}\equiv\end{math}0).{\it{}normal}());

\nwalsodefined{\\{NWgqRer-1Pc9Jw-2}\\{NWgqRer-1Pc9Jw-3}\\{NWgqRer-1Pc9Jw-4}\\{NWgqRer-1Pc9Jw-5}\\{NWgqRer-1Pc9Jw-6}\\{NWgqRer-1Pc9Jw-7}\\{NWgqRer-1Pc9Jw-8}\\{NWgqRer-1Pc9Jw-9}\\{NWgqRer-1Pc9Jw-A}}\nwused{\\{NWgqRer-nXe8t-4}}\nwidentuses{\\{{\nwixident{formula{\_}out}}{formula:unout}}}\nwindexuse{\nwixident{formula{\_}out}}{formula:unout}{NWgqRer-1Pc9Jw-1}\nwendcode{}\nwbegindocs{84}Then we calculate M\"obius action of those matrix on a point.
\nwenddocs{}\nwbegincode{85}\sublabel{NWgqRer-1Pc9Jw-2}\nwmargintag{{\nwtagstyle{}\subpageref{NWgqRer-1Pc9Jw-2}}}\moddef{Show expressions~{\nwtagstyle{}\subpageref{NWgqRer-1Pc9Jw-1}}}\plusendmoddef\Rm{}\nwstartdeflinemarkup\nwusesondefline{\\{NWgqRer-nXe8t-4}}\nwprevnextdefs{NWgqRer-1Pc9Jw-1}{NWgqRer-1Pc9Jw-3}\nwenddeflinemarkup
{\bf{}dual\_number} {\it{}W}({\it{}clifford\_moebius\_map}({\it{}XC}, {\it{}P}.{\it{}to\_matrix}(), {\it{}e}).{\it{}subs}({\it{}sign}\begin{math}\equiv\end{math}0).{\it{}normal}()),
        {\it{}W1} = {\it{}W}.{\it{}subs}({\bf{}lst}({\it{}u}\begin{math}\equiv\end{math}{\it{}u1}, {\it{}v}\begin{math}\equiv\end{math}{\it{}v1}));
{\it{}formula\_out}({\tt{}"Rotation by x: "}, {\it{}W});

\nwused{\\{NWgqRer-nXe8t-4}}\nwidentuses{\\{{\nwixident{dual{\_}number}}{dual:unnumber}}\\{{\nwixident{formula{\_}out}}{formula:unout}}}\nwindexuse{\nwixident{dual{\_}number}}{dual:unnumber}{NWgqRer-1Pc9Jw-2}\nwindexuse{\nwixident{formula{\_}out}}{formula:unout}{NWgqRer-1Pc9Jw-2}\nwendcode{}\nwbegindocs{86}Next we specialise the above result to the reference point.
\nwenddocs{}\nwbegincode{87}\sublabel{NWgqRer-1Pc9Jw-3}\nwmargintag{{\nwtagstyle{}\subpageref{NWgqRer-1Pc9Jw-3}}}\moddef{Show expressions~{\nwtagstyle{}\subpageref{NWgqRer-1Pc9Jw-1}}}\plusendmoddef\Rm{}\nwstartdeflinemarkup\nwusesondefline{\\{NWgqRer-nXe8t-4}}\nwprevnextdefs{NWgqRer-1Pc9Jw-2}{NWgqRer-1Pc9Jw-4}\nwenddeflinemarkup
{\it{}formula\_out}({\tt{}"Rotation of {\char92}{\char92}((u\_0, v\_0){\char92}{\char92}) by {\char92}{\char92}(x{\char92}{\char92}): "}, \nwindexdefn{\nwixident{formula{\_}out}}{formula:unout}{NWgqRer-1Pc9Jw-3}
        {\it{}W}.{\it{}subs}({\bf{}lst}({\it{}u} \begin{math}\equiv\end{math} {\it{}u0}, {\it{}v} \begin{math}\equiv\end{math} {\it{}v0})).{\it{}subs}({\it{}Inf} \begin{math}\equiv\end{math} {\it{}pow}({\it{}y}, -1)).{\it{}normal}().{\it{}subs}({\it{}y} \begin{math}\equiv\end{math} 0).{\it{}normal}());

\nwused{\\{NWgqRer-nXe8t-4}}\nwidentdefs{\\{{\nwixident{formula{\_}out}}{formula:unout}}}\nwendcode{}\nwbegindocs{88}The expression for the parabolic norm.
\nwenddocs{}\nwbegincode{89}\sublabel{NWgqRer-1Pc9Jw-4}\nwmargintag{{\nwtagstyle{}\subpageref{NWgqRer-1Pc9Jw-4}}}\moddef{Show expressions~{\nwtagstyle{}\subpageref{NWgqRer-1Pc9Jw-1}}}\plusendmoddef\Rm{}\nwstartdeflinemarkup\nwusesondefline{\\{NWgqRer-nXe8t-4}}\nwprevnextdefs{NWgqRer-1Pc9Jw-3}{NWgqRer-1Pc9Jw-5}\nwenddeflinemarkup
{\it{}formula\_out}({\tt{}"Parabolic norm: "}, {\it{}P}.{\it{}norm}());

\nwused{\\{NWgqRer-nXe8t-4}}\nwidentuses{\\{{\nwixident{formula{\_}out}}{formula:unout}}}\nwindexuse{\nwixident{formula{\_}out}}{formula:unout}{NWgqRer-1Pc9Jw-4}\nwendcode{}\nwbegindocs{90}Embedding of reals into dual numbers.
\nwenddocs{}\nwbegincode{91}\sublabel{NWgqRer-1Pc9Jw-5}\nwmargintag{{\nwtagstyle{}\subpageref{NWgqRer-1Pc9Jw-5}}}\moddef{Show expressions~{\nwtagstyle{}\subpageref{NWgqRer-1Pc9Jw-1}}}\plusendmoddef\Rm{}\nwstartdeflinemarkup\nwusesondefline{\\{NWgqRer-nXe8t-4}}\nwprevnextdefs{NWgqRer-1Pc9Jw-4}{NWgqRer-1Pc9Jw-6}\nwenddeflinemarkup
{\it{}formula\_out}({\tt{}"Real number x as a dual number: "}, {\bf{}dual\_number}({\it{}x}));

\nwused{\\{NWgqRer-nXe8t-4}}\nwidentuses{\\{{\nwixident{dual{\_}number}}{dual:unnumber}}\\{{\nwixident{formula{\_}out}}{formula:unout}}}\nwindexuse{\nwixident{dual{\_}number}}{dual:unnumber}{NWgqRer-1Pc9Jw-5}\nwindexuse{\nwixident{formula{\_}out}}{formula:unout}{NWgqRer-1Pc9Jw-5}\nwendcode{}\nwbegindocs{92}The expression for the product of two points.
\nwenddocs{}\nwbegincode{93}\sublabel{NWgqRer-1Pc9Jw-6}\nwmargintag{{\nwtagstyle{}\subpageref{NWgqRer-1Pc9Jw-6}}}\moddef{Show expressions~{\nwtagstyle{}\subpageref{NWgqRer-1Pc9Jw-1}}}\plusendmoddef\Rm{}\nwstartdeflinemarkup\nwusesondefline{\\{NWgqRer-nXe8t-4}}\nwprevnextdefs{NWgqRer-1Pc9Jw-5}{NWgqRer-1Pc9Jw-7}\nwenddeflinemarkup
{\it{}formula\_out}({\tt{}"Product: "}, {\it{}P}\begin{math}\ast\end{math}{\it{}P1});

\nwused{\\{NWgqRer-nXe8t-4}}\nwidentuses{\\{{\nwixident{formula{\_}out}}{formula:unout}}}\nwindexuse{\nwixident{formula{\_}out}}{formula:unout}{NWgqRer-1Pc9Jw-6}\nwendcode{}\nwbegindocs{94}The expression of the product of a point and a scalar.
\nwenddocs{}\nwbegincode{95}\sublabel{NWgqRer-1Pc9Jw-7}\nwmargintag{{\nwtagstyle{}\subpageref{NWgqRer-1Pc9Jw-7}}}\moddef{Show expressions~{\nwtagstyle{}\subpageref{NWgqRer-1Pc9Jw-1}}}\plusendmoddef\Rm{}\nwstartdeflinemarkup\nwusesondefline{\\{NWgqRer-nXe8t-4}}\nwprevnextdefs{NWgqRer-1Pc9Jw-6}{NWgqRer-1Pc9Jw-8}\nwenddeflinemarkup
{\it{}formula\_out}({\tt{}"Product by a scalar: "}, {\it{}a}\begin{math}\ast\end{math}{\it{}P});

\nwused{\\{NWgqRer-nXe8t-4}}\nwidentuses{\\{{\nwixident{formula{\_}out}}{formula:unout}}}\nwindexuse{\nwixident{formula{\_}out}}{formula:unout}{NWgqRer-1Pc9Jw-7}\nwendcode{}\nwbegindocs{96}Expressions for the real and imaginary parts.
\nwenddocs{}\nwbegincode{97}\sublabel{NWgqRer-1Pc9Jw-8}\nwmargintag{{\nwtagstyle{}\subpageref{NWgqRer-1Pc9Jw-8}}}\moddef{Show expressions~{\nwtagstyle{}\subpageref{NWgqRer-1Pc9Jw-1}}}\plusendmoddef\Rm{}\nwstartdeflinemarkup\nwusesondefline{\\{NWgqRer-nXe8t-4}}\nwprevnextdefs{NWgqRer-1Pc9Jw-7}{NWgqRer-1Pc9Jw-9}\nwenddeflinemarkup
{\it{}formula\_out}({\tt{}"Real part: "}, {\it{}P}.{\it{}real\_part}());
{\it{}formula\_out}({\tt{}"Imag part: "}, {\it{}P}.{\it{}imag\_part}());

\nwused{\\{NWgqRer-nXe8t-4}}\nwidentuses{\\{{\nwixident{formula{\_}out}}{formula:unout}}}\nwindexuse{\nwixident{formula{\_}out}}{formula:unout}{NWgqRer-1Pc9Jw-8}\nwendcode{}\nwbegindocs{98}The expression for a sum of two points is too cumbersome to be
printed. 
\nwenddocs{}\nwbegincode{99}\sublabel{NWgqRer-1Pc9Jw-9}\nwmargintag{{\nwtagstyle{}\subpageref{NWgqRer-1Pc9Jw-9}}}\moddef{Show expressions~{\nwtagstyle{}\subpageref{NWgqRer-1Pc9Jw-1}}}\plusendmoddef\Rm{}\nwstartdeflinemarkup\nwusesondefline{\\{NWgqRer-nXe8t-4}}\nwprevnextdefs{NWgqRer-1Pc9Jw-8}{NWgqRer-1Pc9Jw-A}\nwenddeflinemarkup
//formula\_out("Add is: ", (P+P1).normal());

\nwused{\\{NWgqRer-nXe8t-4}}\nwidentuses{\\{{\nwixident{formula{\_}out}}{formula:unout}}}\nwindexuse{\nwixident{formula{\_}out}}{formula:unout}{NWgqRer-1Pc9Jw-9}\nwendcode{}\nwbegindocs{100}Linear combination of points \((1,0)\) and \((-1,0)\) with
coefficients \(a\) and \(b\), for the linearisation presented in
\S~\ref{sec:line-exot-form}.
\nwenddocs{}\nwbegincode{101}\sublabel{NWgqRer-1Pc9Jw-A}\nwmargintag{{\nwtagstyle{}\subpageref{NWgqRer-1Pc9Jw-A}}}\moddef{Show expressions~{\nwtagstyle{}\subpageref{NWgqRer-1Pc9Jw-1}}}\plusendmoddef\Rm{}\nwstartdeflinemarkup\nwusesondefline{\\{NWgqRer-nXe8t-4}}\nwprevnextdefs{NWgqRer-1Pc9Jw-9}{\relax}\nwenddeflinemarkup
{\it{}formula\_out}({\tt{}"Lin comb of two vectors a*(1, 0)+b*(-1, 0): "}, \nwindexdefn{\nwixident{formula{\_}out}}{formula:unout}{NWgqRer-1Pc9Jw-A}
                ({\it{}a}\begin{math}\ast\end{math}{\it{}P}+{\it{}b}\begin{math}\ast\end{math}{\it{}P1}).{\it{}subs}({\bf{}lst}({\it{}u}\begin{math}\equiv\end{math}1, {\it{}v}\begin{math}\equiv\end{math}0, {\it{}u1}\begin{math}\equiv\end{math}-1, {\it{}v1}\begin{math}\equiv\end{math}0)).{\it{}normal}());

\nwused{\\{NWgqRer-nXe8t-4}}\nwidentdefs{\\{{\nwixident{formula{\_}out}}{formula:unout}}}\nwendcode{}\nwbegindocs{102}\nwdocspar
\subsubsection{Checking Algebraic Identities }
\label{sec:check-algebr-ident}
In this Subsection we verify basic algebraic properties of the defined operations.

\nwenddocs{}\nwbegindocs{103}A dual number is the sum of its real and imaginary parts.
\nwenddocs{}\nwbegincode{104}\sublabel{NWgqRer-22TNkn-1}\nwmargintag{{\nwtagstyle{}\subpageref{NWgqRer-22TNkn-1}}}\moddef{Check identities~{\nwtagstyle{}\subpageref{NWgqRer-22TNkn-1}}}\endmoddef\Rm{}\nwstartdeflinemarkup\nwusesondefline{\\{NWgqRer-nXe8t-4}}\nwprevnextdefs{\relax}{NWgqRer-22TNkn-2}\nwenddeflinemarkup
{\it{}test\_out}({\tt{}"P is the sum Re(P) and Im(P): "}, \nwindexdefn{\nwixident{test{\_}out}}{test:unout}{NWgqRer-22TNkn-1}
        {\it{}P}-({\it{}ex\_to}\begin{math}<\end{math}{\bf{}dual\_number}\begin{math}>\end{math}({\it{}P}.{\it{}real\_part}())+{\it{}ex\_to}\begin{math}<\end{math}{\bf{}dual\_number}\begin{math}>\end{math}({\it{}P}.{\it{}imag\_part}())));

\nwalsodefined{\\{NWgqRer-22TNkn-2}\\{NWgqRer-22TNkn-3}\\{NWgqRer-22TNkn-4}\\{NWgqRer-22TNkn-5}\\{NWgqRer-22TNkn-6}\\{NWgqRer-22TNkn-7}\\{NWgqRer-22TNkn-8}\\{NWgqRer-22TNkn-9}\\{NWgqRer-22TNkn-A}\\{NWgqRer-22TNkn-B}\\{NWgqRer-22TNkn-C}\\{NWgqRer-22TNkn-D}\\{NWgqRer-22TNkn-E}\\{NWgqRer-22TNkn-F}}\nwused{\\{NWgqRer-nXe8t-4}}\nwidentdefs{\\{{\nwixident{test{\_}out}}{test:unout}}}\nwidentuses{\\{{\nwixident{dual{\_}number}}{dual:unnumber}}}\nwindexuse{\nwixident{dual{\_}number}}{dual:unnumber}{NWgqRer-22TNkn-1}\nwendcode{}\nwbegindocs{105}A dual number maid out of a real \(a\) has the norm of real part
equal to \(a\).
\nwenddocs{}\nwbegincode{106}\sublabel{NWgqRer-22TNkn-2}\nwmargintag{{\nwtagstyle{}\subpageref{NWgqRer-22TNkn-2}}}\moddef{Check identities~{\nwtagstyle{}\subpageref{NWgqRer-22TNkn-1}}}\plusendmoddef\Rm{}\nwstartdeflinemarkup\nwusesondefline{\\{NWgqRer-nXe8t-4}}\nwprevnextdefs{NWgqRer-22TNkn-1}{NWgqRer-22TNkn-3}\nwenddeflinemarkup
{\it{}test\_out}({\tt{}"The real part of a real dual number is itself: "}, \nwindexdefn{\nwixident{test{\_}out}}{test:unout}{NWgqRer-22TNkn-2}
        {\it{}ex\_to}\begin{math}<\end{math}{\bf{}dual\_number}\begin{math}>\end{math}({\bf{}dual\_number}({\it{}a}).{\it{}real\_part}()).{\it{}norm}()-{\it{}a});

\nwused{\\{NWgqRer-nXe8t-4}}\nwidentdefs{\\{{\nwixident{test{\_}out}}{test:unout}}}\nwidentuses{\\{{\nwixident{dual{\_}number}}{dual:unnumber}}}\nwindexuse{\nwixident{dual{\_}number}}{dual:unnumber}{NWgqRer-22TNkn-2}\nwendcode{}\nwbegindocs{107}The norm is invariant under parabolic rotations, i.e. they are in
agreement with Defn.~\ref{de:norm}.
\nwenddocs{}\nwbegincode{108}\sublabel{NWgqRer-22TNkn-3}\nwmargintag{{\nwtagstyle{}\subpageref{NWgqRer-22TNkn-3}}}\moddef{Check identities~{\nwtagstyle{}\subpageref{NWgqRer-22TNkn-1}}}\plusendmoddef\Rm{}\nwstartdeflinemarkup\nwusesondefline{\\{NWgqRer-nXe8t-4}}\nwprevnextdefs{NWgqRer-22TNkn-2}{NWgqRer-22TNkn-4}\nwenddeflinemarkup
{\it{}test\_out}({\tt{}"norm is invariant: "}, {\it{}P}.{\it{}norm}()-{\it{}W}.{\it{}norm}());

\nwused{\\{NWgqRer-nXe8t-4}}\nwidentuses{\\{{\nwixident{test{\_}out}}{test:unout}}}\nwindexuse{\nwixident{test{\_}out}}{test:unout}{NWgqRer-22TNkn-3}\nwendcode{}\nwbegindocs{109}The product \(w_1\bar{w}_2\) is invariant under rotations,
Prop.~\ref{item:prod-inv}. 
\nwenddocs{}\nwbegincode{110}\sublabel{NWgqRer-22TNkn-4}\nwmargintag{{\nwtagstyle{}\subpageref{NWgqRer-22TNkn-4}}}\moddef{Check identities~{\nwtagstyle{}\subpageref{NWgqRer-22TNkn-1}}}\plusendmoddef\Rm{}\nwstartdeflinemarkup\nwusesondefline{\\{NWgqRer-nXe8t-4}}\nwprevnextdefs{NWgqRer-22TNkn-3}{NWgqRer-22TNkn-5}\nwenddeflinemarkup
{\it{}test\_out}({\tt{}"Product is invariant: "}, {\it{}P}\begin{math}\ast\end{math}{\it{}P1}.{\it{}conjugate}()-{\it{}W}\begin{math}\ast\end{math}{\it{}W1}.{\it{}conjugate}());

\nwused{\\{NWgqRer-nXe8t-4}}\nwidentuses{\\{{\nwixident{test{\_}out}}{test:unout}}}\nwindexuse{\nwixident{test{\_}out}}{test:unout}{NWgqRer-22TNkn-4}\nwendcode{}\nwbegindocs{111}Product \(w\bar{w}\) is \((0, \modulus{w}^2)\), Prop.~\ref{item:prod-norm-sq}.
\nwenddocs{}\nwbegincode{112}\sublabel{NWgqRer-22TNkn-5}\nwmargintag{{\nwtagstyle{}\subpageref{NWgqRer-22TNkn-5}}}\moddef{Check identities~{\nwtagstyle{}\subpageref{NWgqRer-22TNkn-1}}}\plusendmoddef\Rm{}\nwstartdeflinemarkup\nwusesondefline{\\{NWgqRer-nXe8t-4}}\nwprevnextdefs{NWgqRer-22TNkn-4}{NWgqRer-22TNkn-6}\nwenddeflinemarkup
{\it{}test\_out}({\tt{}"Product is norm squared: "},\nwindexdefn{\nwixident{test{\_}out}}{test:unout}{NWgqRer-22TNkn-5}
                 ({\it{}P}\begin{math}\ast\end{math}{\it{}P}.{\it{}conjugate}()-{\it{}dn\_from\_arg\_mod}({\it{}Arg0}, {\it{}pow}({\it{}P}.{\it{}norm}(), 2))));

\nwused{\\{NWgqRer-nXe8t-4}}\nwidentdefs{\\{{\nwixident{test{\_}out}}{test:unout}}}\nwidentuses{\\{{\nwixident{Arg0}}{Arg0}}}\nwindexuse{\nwixident{Arg0}}{Arg0}{NWgqRer-22TNkn-5}\nwendcode{}\nwbegindocs{113}The reference point is unit under multiplication.
\nwenddocs{}\nwbegincode{114}\sublabel{NWgqRer-22TNkn-6}\nwmargintag{{\nwtagstyle{}\subpageref{NWgqRer-22TNkn-6}}}\moddef{Check identities~{\nwtagstyle{}\subpageref{NWgqRer-22TNkn-1}}}\plusendmoddef\Rm{}\nwstartdeflinemarkup\nwusesondefline{\\{NWgqRer-nXe8t-4}}\nwprevnextdefs{NWgqRer-22TNkn-5}{NWgqRer-22TNkn-7}\nwenddeflinemarkup
{\it{}test\_out}({\tt{}"Product {\char92}{\char92}((u, v)*(u\_0, v\_0){\char92}{\char92}) is {\char92}{\char92}((u, v){\char92}{\char92}): "}, {\it{}P}\begin{math}\ast\end{math}{\it{}P0}-{\it{}P});

\nwused{\\{NWgqRer-nXe8t-4}}\nwidentuses{\\{{\nwixident{test{\_}out}}{test:unout}}}\nwindexuse{\nwixident{test{\_}out}}{test:unout}{NWgqRer-22TNkn-6}\nwendcode{}\nwbegindocs{115}Addition is commutative, Prop.~\ref{item:add-is-comm-ass}.
\nwenddocs{}\nwbegincode{116}\sublabel{NWgqRer-22TNkn-7}\nwmargintag{{\nwtagstyle{}\subpageref{NWgqRer-22TNkn-7}}}\moddef{Check identities~{\nwtagstyle{}\subpageref{NWgqRer-22TNkn-1}}}\plusendmoddef\Rm{}\nwstartdeflinemarkup\nwusesondefline{\\{NWgqRer-nXe8t-4}}\nwprevnextdefs{NWgqRer-22TNkn-6}{NWgqRer-22TNkn-8}\nwenddeflinemarkup
{\it{}test\_out}({\tt{}"Add is commutative: "}, ({\it{}P}+{\it{}P1})-({\it{}P1}+{\it{}P}));

\nwused{\\{NWgqRer-nXe8t-4}}\nwidentuses{\\{{\nwixident{test{\_}out}}{test:unout}}}\nwindexuse{\nwixident{test{\_}out}}{test:unout}{NWgqRer-22TNkn-7}\nwendcode{}\nwbegindocs{117}Addition is associative, Prop.~\ref{item:add-is-comm-ass}.
\nwenddocs{}\nwbegincode{118}\sublabel{NWgqRer-22TNkn-8}\nwmargintag{{\nwtagstyle{}\subpageref{NWgqRer-22TNkn-8}}}\moddef{Check identities~{\nwtagstyle{}\subpageref{NWgqRer-22TNkn-1}}}\plusendmoddef\Rm{}\nwstartdeflinemarkup\nwusesondefline{\\{NWgqRer-nXe8t-4}}\nwprevnextdefs{NWgqRer-22TNkn-7}{NWgqRer-22TNkn-9}\nwenddeflinemarkup
{\it{}test\_out}({\tt{}"Add is associative: "}, (({\it{}P}+{\it{}P1})+ {\it{}P2})-({\it{}P}+({\it{}P1}+{\it{}P2})));

\nwused{\\{NWgqRer-nXe8t-4}}\nwidentuses{\\{{\nwixident{test{\_}out}}{test:unout}}}\nwindexuse{\nwixident{test{\_}out}}{test:unout}{NWgqRer-22TNkn-8}\nwendcode{}\nwbegindocs{119}Multiplication by a scalar is commutative.
\nwenddocs{}\nwbegincode{120}\sublabel{NWgqRer-22TNkn-9}\nwmargintag{{\nwtagstyle{}\subpageref{NWgqRer-22TNkn-9}}}\moddef{Check identities~{\nwtagstyle{}\subpageref{NWgqRer-22TNkn-1}}}\plusendmoddef\Rm{}\nwstartdeflinemarkup\nwusesondefline{\\{NWgqRer-nXe8t-4}}\nwprevnextdefs{NWgqRer-22TNkn-8}{NWgqRer-22TNkn-A}\nwenddeflinemarkup
{\it{}test\_out}({\tt{}"S-mult commutative: "}, {\it{}P}\begin{math}\ast\end{math}{\it{}a}-{\it{}a}\begin{math}\ast\end{math}{\it{}P});

\nwused{\\{NWgqRer-nXe8t-4}}\nwidentuses{\\{{\nwixident{test{\_}out}}{test:unout}}}\nwindexuse{\nwixident{test{\_}out}}{test:unout}{NWgqRer-22TNkn-9}\nwendcode{}\nwbegindocs{121}Multiplication by a scalar is associative.
\nwenddocs{}\nwbegincode{122}\sublabel{NWgqRer-22TNkn-A}\nwmargintag{{\nwtagstyle{}\subpageref{NWgqRer-22TNkn-A}}}\moddef{Check identities~{\nwtagstyle{}\subpageref{NWgqRer-22TNkn-1}}}\plusendmoddef\Rm{}\nwstartdeflinemarkup\nwusesondefline{\\{NWgqRer-nXe8t-4}}\nwprevnextdefs{NWgqRer-22TNkn-9}{NWgqRer-22TNkn-B}\nwenddeflinemarkup
{\it{}test\_out}({\tt{}"S-mult associative: "}, {\it{}b}\begin{math}\ast\end{math}{\it{}P}\begin{math}\ast\end{math}{\it{}a}-{\it{}a}\begin{math}\ast\end{math}{\it{}P}\begin{math}\ast\end{math}{\it{}b});

\nwused{\\{NWgqRer-nXe8t-4}}\nwidentuses{\\{{\nwixident{test{\_}out}}{test:unout}}}\nwindexuse{\nwixident{test{\_}out}}{test:unout}{NWgqRer-22TNkn-A}\nwendcode{}\nwbegindocs{123}Distributive law \(a(w_1+w_2)=aw_1+aw_2\), Prop.~\ref{item:distrib-scalar}.
\nwenddocs{}\nwbegincode{124}\sublabel{NWgqRer-22TNkn-B}\nwmargintag{{\nwtagstyle{}\subpageref{NWgqRer-22TNkn-B}}}\moddef{Check identities~{\nwtagstyle{}\subpageref{NWgqRer-22TNkn-1}}}\plusendmoddef\Rm{}\nwstartdeflinemarkup\nwusesondefline{\\{NWgqRer-nXe8t-4}}\nwprevnextdefs{NWgqRer-22TNkn-A}{NWgqRer-22TNkn-C}\nwenddeflinemarkup
{\it{}test\_out}({\tt{}"S-mult distributive 1: "}, {\it{}a}\begin{math}\ast\end{math}({\it{}P}+{\it{}P1})-({\it{}a}\begin{math}\ast\end{math}{\it{}P} +{\it{}a}\begin{math}\ast\end{math}{\it{}P1}));

\nwused{\\{NWgqRer-nXe8t-4}}\nwidentuses{\\{{\nwixident{test{\_}out}}{test:unout}}}\nwindexuse{\nwixident{test{\_}out}}{test:unout}{NWgqRer-22TNkn-B}\nwendcode{}\nwbegindocs{125}Distributive law \((a+b)w=aw+bw\), Prop.~\ref{item:distrib-scalar}.
\nwenddocs{}\nwbegincode{126}\sublabel{NWgqRer-22TNkn-C}\nwmargintag{{\nwtagstyle{}\subpageref{NWgqRer-22TNkn-C}}}\moddef{Check identities~{\nwtagstyle{}\subpageref{NWgqRer-22TNkn-1}}}\plusendmoddef\Rm{}\nwstartdeflinemarkup\nwusesondefline{\\{NWgqRer-nXe8t-4}}\nwprevnextdefs{NWgqRer-22TNkn-B}{NWgqRer-22TNkn-D}\nwenddeflinemarkup
{\it{}test\_out}({\tt{}"S-mult distributive 2: "}, {\it{}P}\begin{math}\ast\end{math}({\it{}a}+{\it{}b})-({\it{}P}\begin{math}\ast\end{math}{\it{}a} + {\it{}P}\begin{math}\ast\end{math}{\it{}b}));

\nwused{\\{NWgqRer-nXe8t-4}}\nwidentuses{\\{{\nwixident{test{\_}out}}{test:unout}}}\nwindexuse{\nwixident{test{\_}out}}{test:unout}{NWgqRer-22TNkn-C}\nwendcode{}\nwbegindocs{127}Product is commutative, Prop.~\ref{item:prod-comm-ass}.
\nwenddocs{}\nwbegincode{128}\sublabel{NWgqRer-22TNkn-D}\nwmargintag{{\nwtagstyle{}\subpageref{NWgqRer-22TNkn-D}}}\moddef{Check identities~{\nwtagstyle{}\subpageref{NWgqRer-22TNkn-1}}}\plusendmoddef\Rm{}\nwstartdeflinemarkup\nwusesondefline{\\{NWgqRer-nXe8t-4}}\nwprevnextdefs{NWgqRer-22TNkn-C}{NWgqRer-22TNkn-E}\nwenddeflinemarkup
{\it{}test\_out}({\tt{}"Product is symmetric (commutative): "}, {\it{}P}\begin{math}\ast\end{math}{\it{}P1}-{\it{}P1}\begin{math}\ast\end{math}{\it{}P});

\nwused{\\{NWgqRer-nXe8t-4}}\nwidentuses{\\{{\nwixident{test{\_}out}}{test:unout}}}\nwindexuse{\nwixident{test{\_}out}}{test:unout}{NWgqRer-22TNkn-D}\nwendcode{}\nwbegindocs{129}Product is associative, Prop.~\ref{item:prod-comm-ass}.
\nwenddocs{}\nwbegincode{130}\sublabel{NWgqRer-22TNkn-E}\nwmargintag{{\nwtagstyle{}\subpageref{NWgqRer-22TNkn-E}}}\moddef{Check identities~{\nwtagstyle{}\subpageref{NWgqRer-22TNkn-1}}}\plusendmoddef\Rm{}\nwstartdeflinemarkup\nwusesondefline{\\{NWgqRer-nXe8t-4}}\nwprevnextdefs{NWgqRer-22TNkn-D}{NWgqRer-22TNkn-F}\nwenddeflinemarkup
{\it{}test\_out}({\tt{}"Prod is associative: "}, ({\it{}P}\begin{math}\ast\end{math}{\it{}P1})\begin{math}\ast\end{math}{\it{}P2}-{\it{}P}\begin{math}\ast\end{math}({\it{}P1}\begin{math}\ast\end{math}{\it{}P2}));

\nwused{\\{NWgqRer-nXe8t-4}}\nwidentuses{\\{{\nwixident{test{\_}out}}{test:unout}}}\nwindexuse{\nwixident{test{\_}out}}{test:unout}{NWgqRer-22TNkn-E}\nwendcode{}\nwbegindocs{131}Product and addition are distributive, Prop.~\ref{item:distrib}.
\nwenddocs{}\nwbegincode{132}\sublabel{NWgqRer-22TNkn-F}\nwmargintag{{\nwtagstyle{}\subpageref{NWgqRer-22TNkn-F}}}\moddef{Check identities~{\nwtagstyle{}\subpageref{NWgqRer-22TNkn-1}}}\plusendmoddef\Rm{}\nwstartdeflinemarkup\nwusesondefline{\\{NWgqRer-nXe8t-4}}\nwprevnextdefs{NWgqRer-22TNkn-E}{\relax}\nwenddeflinemarkup
{\it{}test\_out}({\tt{}"Product is distributive: "}, ({\it{}P}+{\it{}P1})\begin{math}\ast\end{math}{\it{}P2}-({\it{}P}\begin{math}\ast\end{math}{\it{}P2}+{\it{}P1}\begin{math}\ast\end{math}{\it{}P2}));

\nwused{\\{NWgqRer-nXe8t-4}}\nwidentuses{\\{{\nwixident{test{\_}out}}{test:unout}}}\nwindexuse{\nwixident{test{\_}out}}{test:unout}{NWgqRer-22TNkn-F}\nwendcode{}\nwbegindocs{133}\nwdocspar
\subsection{Induced Representations}
\label{sec:induc-repr1}
 Here we calculate the basic formulae for Section~\ref{sec:induc-repr}.

\nwenddocs{}\nwbegindocs{134}\nwdocspar
\subsubsection{Encoded formulae}
\label{sec:encoded-formulae}

\nwenddocs{}\nwbegindocs{135}This routine encodes the map \(s: \Space{R}{2}\rightarrow
\SL\)~\eqref{eq:s-map}.  
\nwenddocs{}\nwbegincode{136}\sublabel{NWgqRer-11aTCz-1}\nwmargintag{{\nwtagstyle{}\subpageref{NWgqRer-11aTCz-1}}}\moddef{Induced representations routines~{\nwtagstyle{}\subpageref{NWgqRer-11aTCz-1}}}\endmoddef\Rm{}\nwstartdeflinemarkup\nwusesondefline{\\{NWgqRer-1p0Y9w-2}}\nwprevnextdefs{\relax}{NWgqRer-11aTCz-2}\nwenddeflinemarkup
{\bf{}ex} {\it{}s\_map}({\bf{}const} {\bf{}ex} & {\it{}u}, {\bf{}const} {\bf{}ex} & {\it{}v}) {\nwlbrace}
        {\bf{}return} {\bf{}matrix}(2, 2, {\bf{}lst}({\it{}v},{\it{}u},0,1));
{\nwrbrace}

{\bf{}ex} {\it{}s\_map}({\bf{}const} {\bf{}ex} & {\it{}P}) {\nwlbrace}
        {\bf{}if} ({\it{}P}.{\it{}nops}() \begin{math}\equiv\end{math} 2)
                {\bf{}return} {\it{}s\_map}({\it{}P}.{\it{}op}(0), {\it{}P}.{\it{}op}(1));
        {\it{}cerr} \begin{math}\ll\end{math} {\tt{}"s\_map() error: parameter should have two operands"} \begin{math}\ll\end{math} {\it{}endl};
        {\bf{}return} {\it{}s\_map}({\it{}P},1);
{\nwrbrace}

\nwalsodefined{\\{NWgqRer-11aTCz-2}\\{NWgqRer-11aTCz-3}\\{NWgqRer-11aTCz-4}}\nwused{\\{NWgqRer-1p0Y9w-2}}\nwendcode{}\nwbegindocs{137}This routine encodes the map \(r: \SL \rightarrow
H\)~\eqref{eq:r-map}. The first parameter is an element of \(\SL\), the
second---is a generic element of subgroup \(H\).
\nwenddocs{}\nwbegincode{138}\sublabel{NWgqRer-11aTCz-2}\nwmargintag{{\nwtagstyle{}\subpageref{NWgqRer-11aTCz-2}}}\moddef{Induced representations routines~{\nwtagstyle{}\subpageref{NWgqRer-11aTCz-1}}}\plusendmoddef\Rm{}\nwstartdeflinemarkup\nwusesondefline{\\{NWgqRer-1p0Y9w-2}}\nwprevnextdefs{NWgqRer-11aTCz-1}{NWgqRer-11aTCz-3}\nwenddeflinemarkup
{\bf{}ex} {\it{}r\_map}({\bf{}const} {\bf{}ex} & {\it{}M}, {\bf{}const} {\bf{}ex} & {\it{}K}) {\nwlbrace}
        {\bf{}ex} {\it{}K1}={\it{}K}.{\it{}evalm}(), {\it{}K2};
        {\bf{}lst} {\it{}vars} = ({\it{}is\_a}\begin{math}<\end{math}{\bf{}symbol}\begin{math}>\end{math}({\it{}K1}.{\it{}op}(2)) ? {\bf{}lst}({\it{}K1}.{\it{}op}(2)) :  {\bf{}lst}({\it{}K1}.{\it{}op}(1)));
        {\bf{}if} ({\it{}is\_a}\begin{math}<\end{math}{\bf{}symbol}\begin{math}>\end{math}({\it{}K1}.{\it{}op}(3))) {\nwlbrace}
                {\it{}vars} = {\it{}vars}.{\it{}append}({\it{}K1}.{\it{}op}(3));
                {\it{}K2} = {\it{}K1}.{\it{}subs}({\it{}lsolve}({\bf{}lst}(({\it{}M}\begin{math}\ast\end{math}{\it{}K1}).{\it{}evalm}().{\it{}op}(2)\begin{math}\equiv\end{math}0), {\it{}vars})).{\it{}subs}({\it{}K1}.{\it{}op}(3)\begin{math}\equiv\end{math}1);
        {\nwrbrace} {\bf{}else} 
                {\it{}K2} = {\it{}K1}.{\it{}subs}({\it{}lsolve}({\bf{}lst}(({\it{}M}\begin{math}\ast\end{math}{\it{}K1}).{\it{}evalm}().{\it{}op}(2)\begin{math}\equiv\end{math}0), {\it{}vars}));
        {\bf{}return} {\it{}pow}({\it{}K2}, -1).{\it{}evalm}();
{\nwrbrace}

\nwused{\\{NWgqRer-1p0Y9w-2}}\nwendcode{}\nwbegindocs{139}This is the inverse \(s^{-1}\) of the above map \(s\).
\nwenddocs{}\nwbegincode{140}\sublabel{NWgqRer-11aTCz-3}\nwmargintag{{\nwtagstyle{}\subpageref{NWgqRer-11aTCz-3}}}\moddef{Induced representations routines~{\nwtagstyle{}\subpageref{NWgqRer-11aTCz-1}}}\plusendmoddef\Rm{}\nwstartdeflinemarkup\nwusesondefline{\\{NWgqRer-1p0Y9w-2}}\nwprevnextdefs{NWgqRer-11aTCz-2}{NWgqRer-11aTCz-4}\nwenddeflinemarkup
{\bf{}ex} {\it{}s\_inv}({\bf{}const} {\bf{}ex} & {\it{}M}, {\bf{}const} {\bf{}ex} & {\it{}K}) {\nwlbrace}
        {\bf{}ex} {\it{}MK}=({\it{}M}\begin{math}\ast\end{math}{\it{}pow}({\it{}r\_map}({\it{}M},{\it{}K}),-1)).{\it{}evalm}();
        {\bf{}ex} {\it{}D}={\it{}MK}.{\it{}op}(3).{\it{}subs}({\it{}x}\begin{math}\equiv\end{math}1).{\it{}normal}();
        {\bf{}return} {\bf{}matrix}(1, 2, {\bf{}lst}(({\it{}MK}.{\it{}op}(1).{\it{}subs}({\it{}x}\begin{math}\equiv\end{math}1).{\it{}normal}()\begin{math}\div\end{math}{\it{}D}).{\it{}normal}(), 
                                                        ({\it{}MK}.{\it{}op}(0).{\it{}subs}({\it{}x}\begin{math}\equiv\end{math}1).{\it{}normal}()\begin{math}\div\end{math}{\it{}D}).{\it{}normal}()));
{\nwrbrace}

\nwused{\\{NWgqRer-1p0Y9w-2}}\nwendcode{}\nwbegindocs{141}This is a matrix form of the above inverse map {\Tt{}\Rm{}{\it{}s\_inv}()\nwendquote}.
\nwenddocs{}\nwbegincode{142}\sublabel{NWgqRer-11aTCz-4}\nwmargintag{{\nwtagstyle{}\subpageref{NWgqRer-11aTCz-4}}}\moddef{Induced representations routines~{\nwtagstyle{}\subpageref{NWgqRer-11aTCz-1}}}\plusendmoddef\Rm{}\nwstartdeflinemarkup\nwusesondefline{\\{NWgqRer-1p0Y9w-2}}\nwprevnextdefs{NWgqRer-11aTCz-3}{\relax}\nwenddeflinemarkup
{\bf{}ex} {\it{}s\_inv\_m}({\bf{}const} {\bf{}ex} & {\it{}M}, {\bf{}const} {\bf{}ex} & {\it{}K}) {\nwlbrace}
        {\bf{}return} ({\it{}M}\begin{math}\ast\end{math}{\it{}pow}({\it{}r\_map}({\it{}M},{\it{}K}),-1)).{\it{}evalm}();
{\nwrbrace}

\nwused{\\{NWgqRer-1p0Y9w-2}}\nwendcode{}\nwbegindocs{143}\nwdocspar
\subsubsection{Caculation of induced representation formulae}
\label{sec:cacul-induc-repr}

\nwenddocs{}\nwbegindocs{144} Firstly we define a generic element {\Tt{}\Rm{}{\it{}M}\nwendquote} of \(\SL\).
\nwenddocs{}\nwbegincode{145}\sublabel{NWgqRer-1h89LW-1}\nwmargintag{{\nwtagstyle{}\subpageref{NWgqRer-1h89LW-1}}}\moddef{Induced representations~{\nwtagstyle{}\subpageref{NWgqRer-1h89LW-1}}}\endmoddef\Rm{}\nwstartdeflinemarkup\nwusesondefline{\\{NWgqRer-3sxui-4}}\nwprevnextdefs{\relax}{NWgqRer-1h89LW-2}\nwenddeflinemarkup
{\bf{}ex} {\it{}M}={\bf{}matrix}(2,2, {\bf{}lst}({\it{}a},{\it{}b},{\it{}c},{\it{}d})), {\it{}H};

\nwalsodefined{\\{NWgqRer-1h89LW-2}\\{NWgqRer-1h89LW-3}\\{NWgqRer-1h89LW-4}\\{NWgqRer-1h89LW-5}\\{NWgqRer-1h89LW-6}\\{NWgqRer-1h89LW-7}\\{NWgqRer-1h89LW-8}}\nwused{\\{NWgqRer-3sxui-4}}\nwendcode{}\nwbegindocs{146}We consider the three cases.
\nwenddocs{}\nwbegincode{147}\sublabel{NWgqRer-1h89LW-2}\nwmargintag{{\nwtagstyle{}\subpageref{NWgqRer-1h89LW-2}}}\moddef{Induced representations~{\nwtagstyle{}\subpageref{NWgqRer-1h89LW-1}}}\plusendmoddef\Rm{}\nwstartdeflinemarkup\nwusesondefline{\\{NWgqRer-3sxui-4}}\nwprevnextdefs{NWgqRer-1h89LW-1}{NWgqRer-1h89LW-3}\nwenddeflinemarkup
{\bf{}char}\begin{math}\ast\end{math} {\it{}cases}[]={\nwlbrace}{\tt{}"Elliptic"}, {\tt{}"Parabolic ({\char92}{\char92}(N^{\char92}{\char92}prime{\char92}{\char92}))"}, {\tt{}"Hyperbolic"}{\nwrbrace};\nwindexdefn{\nwixident{cases}}{cases}{NWgqRer-1h89LW-2}

\nwused{\\{NWgqRer-3sxui-4}}\nwidentdefs{\\{{\nwixident{cases}}{cases}}}\nwendcode{}\nwbegindocs{148}In the those cases {\Tt{}\Rm{}{\it{}subgroups}\nwendquote} holds a generic element of a
subgroup \(H\), see~\eqref{eq:matrix-exp-eh}
and~\eqref{eq:parab-prime-exp-matrix}.  
\nwenddocs{}\nwbegincode{149}\sublabel{NWgqRer-1h89LW-3}\nwmargintag{{\nwtagstyle{}\subpageref{NWgqRer-1h89LW-3}}}\moddef{Induced representations~{\nwtagstyle{}\subpageref{NWgqRer-1h89LW-1}}}\plusendmoddef\Rm{}\nwstartdeflinemarkup\nwusesondefline{\\{NWgqRer-3sxui-4}}\nwprevnextdefs{NWgqRer-1h89LW-2}{NWgqRer-1h89LW-4}\nwenddeflinemarkup
{\bf{}ex} {\it{}subgroups}={\bf{}lst}({\bf{}matrix}(2, 2, {\bf{}lst}({\it{}x},-{\it{}y},{\it{}y},{\it{}x})),
                                   {\bf{}matrix}(2, 2, {\bf{}lst}(1,0,{\it{}y},1)),
                                   {\bf{}matrix}(2, 2, {\bf{}lst}({\it{}x},{\it{}y},{\it{}y},{\it{}x})));

\nwused{\\{NWgqRer-3sxui-4}}\nwendcode{}\nwbegindocs{150}Now we run a cycle over the three cases\ldots
\nwenddocs{}\nwbegincode{151}\sublabel{NWgqRer-1h89LW-4}\nwmargintag{{\nwtagstyle{}\subpageref{NWgqRer-1h89LW-4}}}\moddef{Induced representations~{\nwtagstyle{}\subpageref{NWgqRer-1h89LW-1}}}\plusendmoddef\Rm{}\nwstartdeflinemarkup\nwusesondefline{\\{NWgqRer-3sxui-4}}\nwprevnextdefs{NWgqRer-1h89LW-3}{NWgqRer-1h89LW-5}\nwenddeflinemarkup
{\bf{}for}({\bf{}int} {\it{}i}=0; {\it{}i}\begin{math}<\end{math}3; {\it{}i}\protect\PP) {\nwlbrace}
        {\it{}H}={\it{}subgroups}[{\it{}i}];
        {\it{}cout} \begin{math}\ll\end{math} {\it{}cases}[{\it{}i}] \begin{math}\ll\end{math} {\tt{}" case of induced representations{\char92}{\char92}{\char92}{\char92}"} \begin{math}\ll\end{math} {\it{}endl};
        //formula\_out("M*H: ", (M*H).evalm());

\nwused{\\{NWgqRer-3sxui-4}}\nwidentuses{\\{{\nwixident{cases}}{cases}}\\{{\nwixident{formula{\_}out}}{formula:unout}}}\nwindexuse{\nwixident{cases}}{cases}{NWgqRer-1h89LW-4}\nwindexuse{\nwixident{formula{\_}out}}{formula:unout}{NWgqRer-1h89LW-4}\nwendcode{}\nwbegindocs{152}\ldots and output expression of \(r\)~\eqref{eq:r-map},\ldots 
\nwenddocs{}\nwbegincode{153}\sublabel{NWgqRer-1h89LW-5}\nwmargintag{{\nwtagstyle{}\subpageref{NWgqRer-1h89LW-5}}}\moddef{Induced representations~{\nwtagstyle{}\subpageref{NWgqRer-1h89LW-1}}}\plusendmoddef\Rm{}\nwstartdeflinemarkup\nwusesondefline{\\{NWgqRer-3sxui-4}}\nwprevnextdefs{NWgqRer-1h89LW-4}{NWgqRer-1h89LW-6}\nwenddeflinemarkup
        {\it{}formula\_out}({\tt{}"map {\char92}{\char92}(r(M){\char92}{\char92}): "}, {\it{}r\_map}({\it{}M},{\it{}H}));

\nwused{\\{NWgqRer-3sxui-4}}\nwidentuses{\\{{\nwixident{formula{\_}out}}{formula:unout}}}\nwindexuse{\nwixident{formula{\_}out}}{formula:unout}{NWgqRer-1h89LW-5}\nwendcode{}\nwbegindocs{154}\ldots matrix form of the inverse \(s^{-1}\)~\eqref{eq:r-map},\ldots 
\nwenddocs{}\nwbegincode{155}\sublabel{NWgqRer-1h89LW-6}\nwmargintag{{\nwtagstyle{}\subpageref{NWgqRer-1h89LW-6}}}\moddef{Induced representations~{\nwtagstyle{}\subpageref{NWgqRer-1h89LW-1}}}\plusendmoddef\Rm{}\nwstartdeflinemarkup\nwusesondefline{\\{NWgqRer-3sxui-4}}\nwprevnextdefs{NWgqRer-1h89LW-5}{NWgqRer-1h89LW-7}\nwenddeflinemarkup
        {\it{}formula\_out}({\tt{}"map {\char92}{\char92}(s^{\char123}-1{\char125}(M){\char92}{\char92}): "}, {\it{}s\_inv\_m}({\it{}M},{\it{}H}).{\it{}subs}({\it{}a}\begin{math}\equiv\end{math}(1+{\it{}b}\begin{math}\ast\end{math}{\it{}c})\begin{math}\div\end{math}{\it{}d}).{\it{}normal}());

\nwused{\\{NWgqRer-3sxui-4}}\nwidentuses{\\{{\nwixident{formula{\_}out}}{formula:unout}}}\nwindexuse{\nwixident{formula{\_}out}}{formula:unout}{NWgqRer-1h89LW-6}\nwendcode{}\nwbegindocs{156}\ldots expression for the argument of the character in~\eqref{eq:def-ind},\ldots 
\nwenddocs{}\nwbegincode{157}\sublabel{NWgqRer-1h89LW-7}\nwmargintag{{\nwtagstyle{}\subpageref{NWgqRer-1h89LW-7}}}\moddef{Induced representations~{\nwtagstyle{}\subpageref{NWgqRer-1h89LW-1}}}\plusendmoddef\Rm{}\nwstartdeflinemarkup\nwusesondefline{\\{NWgqRer-3sxui-4}}\nwprevnextdefs{NWgqRer-1h89LW-6}{NWgqRer-1h89LW-8}\nwenddeflinemarkup
        {\it{}formula\_out}({\tt{}"character: "}, {\it{}r\_map}({\it{}M}\begin{math}\ast\end{math}{\it{}s\_map}({\it{}P}),{\it{}H}));

\nwused{\\{NWgqRer-3sxui-4}}\nwidentuses{\\{{\nwixident{formula{\_}out}}{formula:unout}}}\nwindexuse{\nwixident{formula{\_}out}}{formula:unout}{NWgqRer-1h89LW-7}\nwendcode{}\nwbegindocs{158}\ldots and finally the action~\eqref{eq:g-action} of \(\SL\) on the
homogeneous space.
\nwenddocs{}\nwbegincode{159}\sublabel{NWgqRer-1h89LW-8}\nwmargintag{{\nwtagstyle{}\subpageref{NWgqRer-1h89LW-8}}}\moddef{Induced representations~{\nwtagstyle{}\subpageref{NWgqRer-1h89LW-1}}}\plusendmoddef\Rm{}\nwstartdeflinemarkup\nwusesondefline{\\{NWgqRer-3sxui-4}}\nwprevnextdefs{NWgqRer-1h89LW-7}{\relax}\nwenddeflinemarkup
        {\it{}formula\_out}({\tt{}"Moebius map: "}, {\it{}s\_inv}({\it{}M}\begin{math}\ast\end{math}{\it{}s\_map}({\it{}P}.{\it{}to\_matrix}()),{\it{}H}));
        {\it{}test\_out}({\tt{}"Moebius map is given by the imaginary unit: "}, {\it{}s\_inv}({\it{}M}\begin{math}\ast\end{math}{\it{}s\_map}({\it{}P}),{\it{}H}) - 
                         {\it{}clifford\_moebius\_map}({\it{}a}\begin{math}\ast\end{math}{\it{}one}, {\it{}b}\begin{math}\ast\end{math}{\it{}e0}, -{\it{}c}\begin{math}\ast\end{math}{\it{}e0}, {\it{}d}\begin{math}\ast\end{math}{\it{}one},{\it{}P}.{\it{}to\_matrix}(),{\it{}e}).{\it{}subs}({\it{}sign}\begin{math}\equiv\end{math}{\it{}i}-1));
        {\it{}cout} \begin{math}\ll\end{math}  ({\it{}latexout} ? {\tt{}"{\char92}{\char92}vspace{\char123}2mm{\char125}{\char92}{\char92}hrule"} : 
                          {\tt{}"----------------------------------------"} ) \begin{math}\ll\end{math} {\it{}endl};
{\nwrbrace}

\nwused{\\{NWgqRer-3sxui-4}}\nwidentuses{\\{{\nwixident{formula{\_}out}}{formula:unout}}\\{{\nwixident{test{\_}out}}{test:unout}}}\nwindexuse{\nwixident{formula{\_}out}}{formula:unout}{NWgqRer-1h89LW-8}\nwindexuse{\nwixident{test{\_}out}}{test:unout}{NWgqRer-1h89LW-8}\nwendcode{}\nwbegindocs{160}\nwdocspar
\subsection{Program Outline}
\label{sec:program-outline}
Here is the outline how we use the above parts.

\nwenddocs{}\nwbegindocs{161}Routines for {\Tt{}\Rm{}{\bf{}dual\_number}\nwendquote} are collected in a separate library. We
start from the definition {\Tt{}\Rm{}{\bf{}dual\_number}\nwendquote} class in the header file.
\nwenddocs{}\nwbegincode{162}\sublabel{NWgqRer-36Ytqo-1}\nwmargintag{{\nwtagstyle{}\subpageref{NWgqRer-36Ytqo-1}}}\moddef{dualnum.h~{\nwtagstyle{}\subpageref{NWgqRer-36Ytqo-1}}}\endmoddef\Rm{}\nwstartdeflinemarkup\nwenddeflinemarkup
 \LA{}Initialisation~{\nwtagstyle{}\subpageref{NWgqRer-40D9Pp-1}}\RA{}
 \LA{}N-Nprime separation~{\nwtagstyle{}\subpageref{NWgqRer-4EWpKB-1}}\RA{}
 \LA{}Dual number class declaration~{\nwtagstyle{}\subpageref{NWgqRer-4GrAIY-1}}\RA{}
 \LA{}Additional routines declarations~{\nwtagstyle{}\subpageref{NWgqRer-NQ5ML-1}}\RA{}

\nwnotused{dualnum.h}\nwendcode{}\nwbegindocs{163}Here is the file with the implementation.
\nwenddocs{}\nwbegincode{164}\sublabel{NWgqRer-3HhVai-1}\nwmargintag{{\nwtagstyle{}\subpageref{NWgqRer-3HhVai-1}}}\moddef{dualnum.cpp~{\nwtagstyle{}\subpageref{NWgqRer-3HhVai-1}}}\endmoddef\Rm{}\nwstartdeflinemarkup\nwenddeflinemarkup
{\bf{}\char35{}include}{\tt{} \begin{math}<\end{math}dualnum.h\begin{math}>\end{math}}
 \LA{}Algebraic procedures~{\nwtagstyle{}\subpageref{NWgqRer-3jzoaG-1}}\RA{}
 \LA{}Dual number class implementation~{\nwtagstyle{}\subpageref{NWgqRer-3uuIiL-1}}\RA{}
 \LA{}Dual number class further implementation~{\nwtagstyle{}\subpageref{NWgqRer-4dxCds-1}}\RA{}
 \LA{}Output routines~{\nwtagstyle{}\subpageref{NWgqRer-4MrGQl-1}}\RA{}

\nwnotused{dualnum.cpp}\nwendcode{}\nwbegindocs{165}\nwdocspar
\subsubsection{Test program outline}
\label{sec:test-program-outline}

Firstly we load {\Tt{}\Rm{}{\bf{}dual\_number}\nwendquote} support.
\nwenddocs{}\nwbegincode{166}\sublabel{NWgqRer-1p0Y9w-1}\nwmargintag{{\nwtagstyle{}\subpageref{NWgqRer-1p0Y9w-1}}}\moddef{*~{\nwtagstyle{}\subpageref{NWgqRer-1p0Y9w-1}}}\endmoddef\Rm{}\nwstartdeflinemarkup\nwprevnextdefs{\relax}{NWgqRer-1p0Y9w-2}\nwenddeflinemarkup
{\bf{}\char35{}include}{\tt{} \begin{math}<\end{math}dualnum.h\begin{math}>\end{math}}

\nwalsodefined{\\{NWgqRer-1p0Y9w-2}}\nwnotused{*}\nwendcode{}\nwbegindocs{167}The rest of the program makes all checks.
\nwenddocs{}\nwbegincode{168}\sublabel{NWgqRer-1p0Y9w-2}\nwmargintag{{\nwtagstyle{}\subpageref{NWgqRer-1p0Y9w-2}}}\moddef{*~{\nwtagstyle{}\subpageref{NWgqRer-1p0Y9w-1}}}\plusendmoddef\Rm{}\nwstartdeflinemarkup\nwprevnextdefs{NWgqRer-1p0Y9w-1}{\relax}\nwenddeflinemarkup
 \LA{}Definition of variables~{\nwtagstyle{}\subpageref{NWgqRer-2iImO3-1}}\RA{}
 \LA{}Test routine~{\nwtagstyle{}\subpageref{NWgqRer-nXe8t-1}}\RA{}
 \LA{}Induced representations routines~{\nwtagstyle{}\subpageref{NWgqRer-11aTCz-1}}\RA{}
 \LA{}Main procedure~{\nwtagstyle{}\subpageref{NWgqRer-3sxui-1}}\RA{}

\nwendcode{}\nwbegindocs{169}This is the initialisation part
\nwenddocs{}\nwbegincode{170}\sublabel{NWgqRer-40D9Pp-1}\nwmargintag{{\nwtagstyle{}\subpageref{NWgqRer-40D9Pp-1}}}\moddef{Initialisation~{\nwtagstyle{}\subpageref{NWgqRer-40D9Pp-1}}}\endmoddef\Rm{}\nwstartdeflinemarkup\nwusesondefline{\\{NWgqRer-36Ytqo-1}}\nwenddeflinemarkup
{\bf{}\char35{}include}{\tt{} \begin{math}<\end{math}cycle.h\begin{math}>\end{math}}
{\bf{}\char35{}include}{\tt{} \begin{math}<\end{math}fstream\begin{math}>\end{math}}
{\bf{}using} {\bf{}namespace} {\it{}std};
{\bf{}using} {\bf{}namespace} {\it{}GiNaC};

{\bf{}\char35{}include}{\tt{} "ginac-utils.h"}

\nwused{\\{NWgqRer-36Ytqo-1}}\nwendcode{}\nwbegindocs{171}\nwdocspar
\subsubsection{Variables}
\label{sec:variables}

\nwenddocs{}\nwbegindocs{172}These {\Tt{}\Rm{}{\bf{}realsymbol}\nwendquote}s are used in our calculations.
\nwenddocs{}\nwbegincode{173}\sublabel{NWgqRer-2iImO3-1}\nwmargintag{{\nwtagstyle{}\subpageref{NWgqRer-2iImO3-1}}}\moddef{Definition of variables~{\nwtagstyle{}\subpageref{NWgqRer-2iImO3-1}}}\endmoddef\Rm{}\nwstartdeflinemarkup\nwusesondefline{\\{NWgqRer-1p0Y9w-2}}\nwprevnextdefs{\relax}{NWgqRer-2iImO3-2}\nwenddeflinemarkup
{\bf{}realsymbol} {\it{}u}({\tt{}"u"}), {\it{}v}({\tt{}"v"}), {\it{}u1}({\tt{}"u'"}), {\it{}v1}({\tt{}"v'"}), {\it{}u2}({\tt{}"u''"}), {\it{}v2}({\tt{}"v''"}),
        {\it{}a}({\tt{}"a"}), {\it{}b}({\tt{}"b"}), {\it{}c}({\tt{}"c"}), {\it{}d}({\tt{}"d"}), {\it{}x}({\tt{}"x"}), {\it{}y}({\tt{}"y"}),

\nwalsodefined{\\{NWgqRer-2iImO3-2}\\{NWgqRer-2iImO3-3}\\{NWgqRer-2iImO3-4}\\{NWgqRer-2iImO3-5}}\nwused{\\{NWgqRer-1p0Y9w-2}}\nwendcode{}\nwbegindocs{174}Finally this variable keeps the signature of the metric space. 
\nwenddocs{}\nwbegincode{175}\sublabel{NWgqRer-2iImO3-2}\nwmargintag{{\nwtagstyle{}\subpageref{NWgqRer-2iImO3-2}}}\moddef{Definition of variables~{\nwtagstyle{}\subpageref{NWgqRer-2iImO3-1}}}\plusendmoddef\Rm{}\nwstartdeflinemarkup\nwusesondefline{\\{NWgqRer-1p0Y9w-2}}\nwprevnextdefs{NWgqRer-2iImO3-1}{NWgqRer-2iImO3-3}\nwenddeflinemarkup
        {\it{}sign}({\tt{}"s"}, {\tt{}"{\char92}{\char92}sigma"});

\nwused{\\{NWgqRer-1p0Y9w-2}}\nwendcode{}\nwbegindocs{176}This an index used for the definition of Clifford units.
\nwenddocs{}\nwbegincode{177}\sublabel{NWgqRer-2iImO3-3}\nwmargintag{{\nwtagstyle{}\subpageref{NWgqRer-2iImO3-3}}}\moddef{Definition of variables~{\nwtagstyle{}\subpageref{NWgqRer-2iImO3-1}}}\plusendmoddef\Rm{}\nwstartdeflinemarkup\nwusesondefline{\\{NWgqRer-1p0Y9w-2}}\nwprevnextdefs{NWgqRer-2iImO3-2}{NWgqRer-2iImO3-4}\nwenddeflinemarkup
{\bf{}varidx} {\it{}mu}({\bf{}symbol}({\tt{}"mu"}, {\tt{}"{\char92}{\char92}mu"}), 2);

\nwused{\\{NWgqRer-1p0Y9w-2}}\nwendcode{}\nwbegindocs{178}Three generic points which are used in calculations.
\nwenddocs{}\nwbegincode{179}\sublabel{NWgqRer-2iImO3-4}\nwmargintag{{\nwtagstyle{}\subpageref{NWgqRer-2iImO3-4}}}\moddef{Definition of variables~{\nwtagstyle{}\subpageref{NWgqRer-2iImO3-1}}}\plusendmoddef\Rm{}\nwstartdeflinemarkup\nwusesondefline{\\{NWgqRer-1p0Y9w-2}}\nwprevnextdefs{NWgqRer-2iImO3-3}{NWgqRer-2iImO3-5}\nwenddeflinemarkup
{\bf{}dual\_number} {\it{}P}({\it{}u}, {\it{}v}), {\it{}P1}({\it{}u1}, {\it{}v1}), {\it{}P2}({\it{}u2}, {\it{}v2});

\nwused{\\{NWgqRer-1p0Y9w-2}}\nwidentuses{\\{{\nwixident{dual{\_}number}}{dual:unnumber}}}\nwindexuse{\nwixident{dual{\_}number}}{dual:unnumber}{NWgqRer-2iImO3-4}\nwendcode{}\nwbegindocs{180}Here we define a parabolic Clifford units {\Tt{}\Rm{}{\it{}e0}\nwendquote}, {\Tt{}\Rm{}{\it{}e1}\nwendquote}.
\nwenddocs{}\nwbegincode{181}\sublabel{NWgqRer-2iImO3-5}\nwmargintag{{\nwtagstyle{}\subpageref{NWgqRer-2iImO3-5}}}\moddef{Definition of variables~{\nwtagstyle{}\subpageref{NWgqRer-2iImO3-1}}}\plusendmoddef\Rm{}\nwstartdeflinemarkup\nwusesondefline{\\{NWgqRer-1p0Y9w-2}}\nwprevnextdefs{NWgqRer-2iImO3-4}{\relax}\nwenddeflinemarkup
{\bf{}ex} {\it{}e} = {\it{}clifford\_unit}({\it{}mu}, {\it{}diag\_matrix}({\bf{}lst}(-1, {\it{}sign}))),
        {\it{}e0} = {\it{}e}.{\it{}subs}({\it{}mu}\begin{math}\equiv\end{math}0),
        {\it{}e1} = {\it{}e}.{\it{}subs}({\it{}mu}\begin{math}\equiv\end{math}1),
        {\it{}one} = {\it{}dirac\_ONE}();

\nwused{\\{NWgqRer-1p0Y9w-2}}\nwendcode{}\nwbegindocs{182}\nwdocspar
\subsubsection{Test routine}
\label{sec:test-routine}
This routine make the same sequence of checks for both cases of
subgroups \(N\) and \(N^\prime\).

\nwenddocs{}\nwbegindocs{183}First we define the reference point {\Tt{}\Rm{}({\it{}u0},{\it{}v0})\nwendquote}.
\nwenddocs{}\nwbegincode{184}\sublabel{NWgqRer-nXe8t-1}\nwmargintag{{\nwtagstyle{}\subpageref{NWgqRer-nXe8t-1}}}\moddef{Test routine~{\nwtagstyle{}\subpageref{NWgqRer-nXe8t-1}}}\endmoddef\Rm{}\nwstartdeflinemarkup\nwusesondefline{\\{NWgqRer-1p0Y9w-2}}\nwprevnextdefs{\relax}{NWgqRer-nXe8t-2}\nwenddeflinemarkup
{\bf{}void} {\it{}parab\_rot\_sub}() {\nwlbrace}\nwindexdefn{\nwixident{parab{\_}rot{\_}sub}}{parab:unrot:unsub}{NWgqRer-nXe8t-1}
        {\bf{}ex} {\it{}X},
                {\it{}u0}={\it{}Arg0}, {\it{}v0}={\it{}v\_from\_norm}({\it{}u0}, 1),
                {\it{}P0}={\bf{}matrix}(1, 2, {\bf{}lst}({\it{}u0}, {\it{}v0})),

\nwalsodefined{\\{NWgqRer-nXe8t-2}\\{NWgqRer-nXe8t-3}\\{NWgqRer-nXe8t-4}}\nwused{\\{NWgqRer-1p0Y9w-2}}\nwidentdefs{\\{{\nwixident{parab{\_}rot{\_}sub}}{parab:unrot:unsub}}}\nwidentuses{\\{{\nwixident{Arg0}}{Arg0}}}\nwindexuse{\nwixident{Arg0}}{Arg0}{NWgqRer-nXe8t-1}\nwendcode{}\nwbegindocs{185}These two matrices define the Cayley transform and its inverse.
\nwenddocs{}\nwbegincode{186}\sublabel{NWgqRer-nXe8t-2}\nwmargintag{{\nwtagstyle{}\subpageref{NWgqRer-nXe8t-2}}}\moddef{Test routine~{\nwtagstyle{}\subpageref{NWgqRer-nXe8t-1}}}\plusendmoddef\Rm{}\nwstartdeflinemarkup\nwusesondefline{\\{NWgqRer-1p0Y9w-2}}\nwprevnextdefs{NWgqRer-nXe8t-1}{NWgqRer-nXe8t-3}\nwenddeflinemarkup
                {\it{}TC}={\bf{}matrix}(2, 2, {\bf{}lst}({\it{}one}, -{\it{}e1}, -{\it{}e1}, {\it{}one})),
                {\it{}TCI}={\bf{}matrix}(2, 2, {\bf{}lst}({\it{}one}, {\it{}e1}, {\it{}e1}, {\it{}one}));

\nwused{\\{NWgqRer-1p0Y9w-2}}\nwendcode{}\nwbegindocs{187}For the subgroup \(N\) we consider upper-triangular matrices, for
\(N^\prime\)---lower-triangular. 
\nwenddocs{}\nwbegincode{188}\sublabel{NWgqRer-nXe8t-3}\nwmargintag{{\nwtagstyle{}\subpageref{NWgqRer-nXe8t-3}}}\moddef{Test routine~{\nwtagstyle{}\subpageref{NWgqRer-nXe8t-1}}}\plusendmoddef\Rm{}\nwstartdeflinemarkup\nwusesondefline{\\{NWgqRer-1p0Y9w-2}}\nwprevnextdefs{NWgqRer-nXe8t-2}{NWgqRer-nXe8t-4}\nwenddeflinemarkup
 {\bf{}if} ({\it{}is\_subgroup\_N})
         {\it{}X}={\bf{}matrix}(2, 2, {\bf{}lst}({\it{}one}, {\it{}e0}\begin{math}\ast\end{math}{\it{}x}, 0, {\it{}one}));
 {\bf{}else}
         {\it{}X}={\bf{}matrix}(2, 2, {\bf{}lst}({\it{}one}, 0, {\it{}e0}\begin{math}\ast\end{math}{\it{}x}, {\it{}one}));

\nwused{\\{NWgqRer-1p0Y9w-2}}\nwendcode{}\nwbegindocs{189}Common part of test routine.
\nwenddocs{}\nwbegincode{190}\sublabel{NWgqRer-nXe8t-4}\nwmargintag{{\nwtagstyle{}\subpageref{NWgqRer-nXe8t-4}}}\moddef{Test routine~{\nwtagstyle{}\subpageref{NWgqRer-nXe8t-1}}}\plusendmoddef\Rm{}\nwstartdeflinemarkup\nwusesondefline{\\{NWgqRer-1p0Y9w-2}}\nwprevnextdefs{NWgqRer-nXe8t-3}{\relax}\nwenddeflinemarkup
 \LA{}Show expressions~{\nwtagstyle{}\subpageref{NWgqRer-1Pc9Jw-1}}\RA{}
 \LA{}Check identities~{\nwtagstyle{}\subpageref{NWgqRer-22TNkn-1}}\RA{}
 {\it{}cout} \begin{math}\ll\end{math}  ({\it{}latexout} ? {\tt{}"{\char92}{\char92}vspace{\char123}2mm{\char125}{\char92}{\char92}hrule"} : 
                   {\tt{}"----------------------------------------"} ) \begin{math}\ll\end{math} {\it{}endl};
{\nwrbrace}

\nwused{\\{NWgqRer-1p0Y9w-2}}\nwendcode{}\nwbegindocs{191}\nwdocspar
\subsubsection{Main procedure}
\label{sec:main-procedure}
It just calls the test routine, calculates the induced representation and
draws a few pictures.

\nwenddocs{}\nwbegindocs{192} We output formulae in \LaTeX\ mode.
\nwenddocs{}\nwbegincode{193}\sublabel{NWgqRer-3sxui-1}\nwmargintag{{\nwtagstyle{}\subpageref{NWgqRer-3sxui-1}}}\moddef{Main procedure~{\nwtagstyle{}\subpageref{NWgqRer-3sxui-1}}}\endmoddef\Rm{}\nwstartdeflinemarkup\nwusesondefline{\\{NWgqRer-1p0Y9w-2}}\nwprevnextdefs{\relax}{NWgqRer-3sxui-2}\nwenddeflinemarkup
{\bf{}int} {\it{}main}(){\nwlbrace}\nwindexdefn{\nwixident{main}}{main}{NWgqRer-3sxui-1}
        {\it{}latexout}={\bf{}true};

\nwalsodefined{\\{NWgqRer-3sxui-2}\\{NWgqRer-3sxui-3}\\{NWgqRer-3sxui-4}\\{NWgqRer-3sxui-5}}\nwused{\\{NWgqRer-1p0Y9w-2}}\nwidentdefs{\\{{\nwixident{main}}{main}}}\nwendcode{}\nwbegindocs{194}Preparation of output stream.
\nwenddocs{}\nwbegincode{195}\sublabel{NWgqRer-3sxui-2}\nwmargintag{{\nwtagstyle{}\subpageref{NWgqRer-3sxui-2}}}\moddef{Main procedure~{\nwtagstyle{}\subpageref{NWgqRer-3sxui-1}}}\plusendmoddef\Rm{}\nwstartdeflinemarkup\nwusesondefline{\\{NWgqRer-1p0Y9w-2}}\nwprevnextdefs{NWgqRer-3sxui-1}{NWgqRer-3sxui-3}\nwenddeflinemarkup
        {\it{}cout} \begin{math}\ll\end{math} {\it{}boolalpha};
        {\bf{}if} ({\it{}latexout})
                {\it{}cout} \begin{math}\ll\end{math} {\it{}latex};

\nwused{\\{NWgqRer-1p0Y9w-2}}\nwendcode{}\nwbegindocs{196}Now we call the test routine for both \(N\) and \(N^\prime\) subgroups.
\nwenddocs{}\nwbegincode{197}\sublabel{NWgqRer-3sxui-3}\nwmargintag{{\nwtagstyle{}\subpageref{NWgqRer-3sxui-3}}}\moddef{Main procedure~{\nwtagstyle{}\subpageref{NWgqRer-3sxui-1}}}\plusendmoddef\Rm{}\nwstartdeflinemarkup\nwusesondefline{\\{NWgqRer-1p0Y9w-2}}\nwprevnextdefs{NWgqRer-3sxui-2}{NWgqRer-3sxui-4}\nwenddeflinemarkup
 {\it{}is\_subgroup\_N} = {\bf{}true};
 {\it{}parab\_rot\_sub}();
 {\it{}is\_subgroup\_N} = {\bf{}false};
 {\it{}parab\_rot\_sub}();

\nwused{\\{NWgqRer-1p0Y9w-2}}\nwidentuses{\\{{\nwixident{parab{\_}rot{\_}sub}}{parab:unrot:unsub}}}\nwindexuse{\nwixident{parab{\_}rot{\_}sub}}{parab:unrot:unsub}{NWgqRer-3sxui-3}\nwendcode{}\nwbegindocs{198} Calculation of induced representations formulae.
\nwenddocs{}\nwbegincode{199}\sublabel{NWgqRer-3sxui-4}\nwmargintag{{\nwtagstyle{}\subpageref{NWgqRer-3sxui-4}}}\moddef{Main procedure~{\nwtagstyle{}\subpageref{NWgqRer-3sxui-1}}}\plusendmoddef\Rm{}\nwstartdeflinemarkup\nwusesondefline{\\{NWgqRer-1p0Y9w-2}}\nwprevnextdefs{NWgqRer-3sxui-3}{NWgqRer-3sxui-5}\nwenddeflinemarkup
 \LA{}Induced representations~{\nwtagstyle{}\subpageref{NWgqRer-1h89LW-1}}\RA{}

\nwused{\\{NWgqRer-1p0Y9w-2}}\nwendcode{}\nwbegindocs{200}And we finishing by drawing several pictures for
Figs.~\ref{fig:rotations} and~\ref{fig:p-rotations}.
\nwenddocs{}\nwbegincode{201}\sublabel{NWgqRer-3sxui-5}\nwmargintag{{\nwtagstyle{}\subpageref{NWgqRer-3sxui-5}}}\moddef{Main procedure~{\nwtagstyle{}\subpageref{NWgqRer-3sxui-1}}}\plusendmoddef\Rm{}\nwstartdeflinemarkup\nwusesondefline{\\{NWgqRer-1p0Y9w-2}}\nwprevnextdefs{NWgqRer-3sxui-4}{\relax}\nwenddeflinemarkup
 \LA{}Drawing pictures~{\nwtagstyle{}\subpageref{NWgqRer-Ghke2-1}}\RA{}
{\nwrbrace}

\nwused{\\{NWgqRer-1p0Y9w-2}}\nwendcode{}\nwbegindocs{202}\nwdocspar
\subsection{Drawing Orbits}
\label{sec:drawing-orbits}
To draw cycles we use {\Tt{}\Rm{}{\bf{}cycle}\nwendquote} library~\cite{Kisil05b}.

\nwenddocs{}\nwbegindocs{203}Elliptic orbits (circles).
\nwenddocs{}\nwbegincode{204}\sublabel{NWgqRer-Ghke2-1}\nwmargintag{{\nwtagstyle{}\subpageref{NWgqRer-Ghke2-1}}}\moddef{Drawing pictures~{\nwtagstyle{}\subpageref{NWgqRer-Ghke2-1}}}\endmoddef\Rm{}\nwstartdeflinemarkup\nwusesondefline{\\{NWgqRer-3sxui-5}}\nwprevnextdefs{\relax}{NWgqRer-Ghke2-2}\nwenddeflinemarkup
{\it{}ofstream} {\it{}asymptote}({\tt{}"parab-rot-data.asy"});
{\it{}asymptote} \begin{math}\ll\end{math} {\tt{}"path[] K="};
{\bf{}for}({\bf{}int} {\it{}i}=0; {\it{}i}\begin{math}<\end{math}6; {\it{}i}\protect\PP)
        {\bf{}cycle2D}({\bf{}lst}(0,0),{\it{}e}.{\it{}subs}({\it{}sign}\begin{math}\equiv\end{math}-1),{\it{}i}\begin{math}\ast\end{math}{\it{}i}\begin{math}\ast\end{math}.04)
                .{\it{}asy\_path}({\it{}asymptote}, -1.5, 1.5, -2, 2, 0, ({\it{}i}\begin{math}>\end{math}0));
{\it{}asymptote} \begin{math}\ll\end{math} {\tt{}";"} \begin{math}\ll\end{math} {\it{}endl};

\nwalsodefined{\\{NWgqRer-Ghke2-2}\\{NWgqRer-Ghke2-3}\\{NWgqRer-Ghke2-4}}\nwused{\\{NWgqRer-3sxui-5}}\nwendcode{}\nwbegindocs{205}Hyperbolic orbits. 
\nwenddocs{}\nwbegincode{206}\sublabel{NWgqRer-Ghke2-2}\nwmargintag{{\nwtagstyle{}\subpageref{NWgqRer-Ghke2-2}}}\moddef{Drawing pictures~{\nwtagstyle{}\subpageref{NWgqRer-Ghke2-1}}}\plusendmoddef\Rm{}\nwstartdeflinemarkup\nwusesondefline{\\{NWgqRer-3sxui-5}}\nwprevnextdefs{NWgqRer-Ghke2-1}{NWgqRer-Ghke2-3}\nwenddeflinemarkup
{\it{}asymptote} \begin{math}\ll\end{math} {\tt{}"path[] A="};
{\bf{}for}({\bf{}int} {\it{}i}=0; {\it{}i}\begin{math}<\end{math}6; {\it{}i}\protect\PP)
        {\bf{}cycle2D}({\bf{}lst}(0,0),{\it{}e}.{\it{}subs}({\it{}sign}\begin{math}\equiv\end{math}1),-{\it{}i}\begin{math}\ast\end{math}{\it{}i}\begin{math}\ast\end{math}.04)
                .{\it{}asy\_path}({\it{}asymptote}, -1.5, 1.5, -1.5, 2, 0, ({\it{}i}\begin{math}>\end{math}0));
{\it{}asymptote} \begin{math}\ll\end{math} {\tt{}";"} \begin{math}\ll\end{math} {\it{}endl};

\nwused{\\{NWgqRer-3sxui-5}}\nwendcode{}\nwbegindocs{207}Parabolic orbits, subgroup \(N\).
\nwenddocs{}\nwbegincode{208}\sublabel{NWgqRer-Ghke2-3}\nwmargintag{{\nwtagstyle{}\subpageref{NWgqRer-Ghke2-3}}}\moddef{Drawing pictures~{\nwtagstyle{}\subpageref{NWgqRer-Ghke2-1}}}\plusendmoddef\Rm{}\nwstartdeflinemarkup\nwusesondefline{\\{NWgqRer-3sxui-5}}\nwprevnextdefs{NWgqRer-Ghke2-2}{NWgqRer-Ghke2-4}\nwenddeflinemarkup
{\it{}asymptote} \begin{math}\ll\end{math} {\tt{}"path[] N="};
{\bf{}for}({\bf{}int} {\it{}i}=0; {\it{}i}\begin{math}<\end{math}6; {\it{}i}\protect\PP)
        {\bf{}cycle2D}(1,{\bf{}lst}(0,{\bf{}numeric}(1,2)),{\bf{}numeric}({\it{}i},2)-1,{\it{}e}.{\it{}subs}({\it{}sign}\begin{math}\equiv\end{math}0))
                .{\it{}asy\_path}({\it{}asymptote}, -1.5, 1.5, -2, 2, 0, ({\it{}i}\begin{math}>\end{math}0));
{\it{}asymptote} \begin{math}\ll\end{math} {\tt{}";"} \begin{math}\ll\end{math} {\it{}endl};

\nwused{\\{NWgqRer-3sxui-5}}\nwendcode{}\nwbegindocs{209} Parabolic orbits, subgroup \(N^\prime\). 
\nwenddocs{}\nwbegincode{210}\sublabel{NWgqRer-Ghke2-4}\nwmargintag{{\nwtagstyle{}\subpageref{NWgqRer-Ghke2-4}}}\moddef{Drawing pictures~{\nwtagstyle{}\subpageref{NWgqRer-Ghke2-1}}}\plusendmoddef\Rm{}\nwstartdeflinemarkup\nwusesondefline{\\{NWgqRer-3sxui-5}}\nwprevnextdefs{NWgqRer-Ghke2-3}{\relax}\nwenddeflinemarkup
{\it{}asymptote} \begin{math}\ll\end{math} {\tt{}"path[] N1="};
{\bf{}for}({\bf{}int} {\it{}i}=0; {\it{}i}\begin{math}<\end{math}5; {\it{}i}\protect\PP)
        {\bf{}cycle2D}(.5\begin{math}\ast\end{math}{\it{}i}\begin{math}\ast\end{math}{\it{}i}\begin{math}\ast\end{math}{\it{}i}+1,{\bf{}lst}(0,{\bf{}numeric}(1,2)),-1,{\it{}e}.{\it{}subs}({\it{}sign}\begin{math}\equiv\end{math}0))
                .{\it{}asy\_path}({\it{}asymptote}, -1.5, 1.5, -1.5, 2, 0, ({\it{}i}\begin{math}>\end{math}0));
{\it{}asymptote} \begin{math}\ll\end{math} {\tt{}";"} \begin{math}\ll\end{math} {\it{}endl};

{\it{}asymptote}.{\it{}close}();

\nwused{\\{NWgqRer-3sxui-5}}\nwendcode{}\nwbegindocs{211}\nwdocspar
\subsubsection{Output routines}
\label{sec:output-routines}

We use standardised routines to output results of calculations.
\nwenddocs{}\nwbegincode{212}\sublabel{NWgqRer-4MrGQl-1}\nwmargintag{{\nwtagstyle{}\subpageref{NWgqRer-4MrGQl-1}}}\moddef{Output routines~{\nwtagstyle{}\subpageref{NWgqRer-4MrGQl-1}}}\endmoddef\Rm{}\nwstartdeflinemarkup\nwusesondefline{\\{NWgqRer-3HhVai-1}}\nwprevnextdefs{\relax}{NWgqRer-4MrGQl-2}\nwenddeflinemarkup
{\bf{}void} {\it{}formula\_out}({\bf{}char}\begin{math}\ast\end{math} {\it{}S}, {\bf{}ex} {\it{}F}) {\nwlbrace}\nwindexdefn{\nwixident{formula{\_}out}}{formula:unout}{NWgqRer-4MrGQl-1}
        {\it{}cout} \begin{math}\ll\end{math} {\it{}S} \begin{math}\ll\end{math} ({\it{}latexout} ? {\tt{}"{\char92}{\char92}("} : {\tt{}""} ) \begin{math}\ll\end{math} {\it{}F} \begin{math}\ll\end{math} ({\it{}latexout} ? {\tt{}"{\char92}{\char92}){\char92}{\char92}{\char92}{\char92}"} : {\tt{}""} )
                 \begin{math}\ll\end{math} {\it{}endl};
{\nwrbrace}

\nwalsodefined{\\{NWgqRer-4MrGQl-2}}\nwused{\\{NWgqRer-3HhVai-1}}\nwidentdefs{\\{{\nwixident{formula{\_}out}}{formula:unout}}}\nwendcode{}\nwbegindocs{213}This routine is used to check identities.
\nwenddocs{}\nwbegincode{214}\sublabel{NWgqRer-4MrGQl-2}\nwmargintag{{\nwtagstyle{}\subpageref{NWgqRer-4MrGQl-2}}}\moddef{Output routines~{\nwtagstyle{}\subpageref{NWgqRer-4MrGQl-1}}}\plusendmoddef\Rm{}\nwstartdeflinemarkup\nwusesondefline{\\{NWgqRer-3HhVai-1}}\nwprevnextdefs{NWgqRer-4MrGQl-1}{\relax}\nwenddeflinemarkup
{\bf{}void} {\it{}test\_out}({\bf{}char}\begin{math}\ast\end{math} {\it{}S}, {\bf{}ex} {\it{}T}) {\nwlbrace}\nwindexdefn{\nwixident{test{\_}out}}{test:unout}{NWgqRer-4MrGQl-2}
        {\it{}cout} \begin{math}\ll\end{math} {\it{}S} \begin{math}\ll\end{math} ({\it{}latexout} ? {\tt{}"{\char92}{\char92}textbf{\char123}"} : {\tt{}"*"} ) 
                 \begin{math}\ll\end{math} ({\it{}is\_a}\begin{math}<\end{math}{\bf{}dual\_number}\begin{math}>\end{math}({\it{}T}) ? {\it{}ex\_to}\begin{math}<\end{math}{\bf{}dual\_number}\begin{math}>\end{math}({\it{}T}).{\it{}normal}().{\it{}is\_zero}() :
                         {\it{}T}.{\it{}evalm}().{\it{}normal}().{\it{}is\_zero\_matrix}()) \begin{math}\ll\end{math} ({\it{}latexout} ? {\tt{}"{\char125}{\char92}{\char92}{\char92}{\char92}"} : {\tt{}"*"} )
                 \begin{math}\ll\end{math} {\it{}endl};
{\nwrbrace}

\nwused{\\{NWgqRer-3HhVai-1}}\nwidentdefs{\\{{\nwixident{test{\_}out}}{test:unout}}}\nwidentuses{\\{{\nwixident{dual{\_}number}}{dual:unnumber}}}\nwindexuse{\nwixident{dual{\_}number}}{dual:unnumber}{NWgqRer-4MrGQl-2}\nwendcode{}\nwbegindocs{215}Here is declarations of additional routines for the header file.
\nwenddocs{}\nwbegincode{216}\sublabel{NWgqRer-NQ5ML-1}\nwmargintag{{\nwtagstyle{}\subpageref{NWgqRer-NQ5ML-1}}}\moddef{Additional routines declarations~{\nwtagstyle{}\subpageref{NWgqRer-NQ5ML-1}}}\endmoddef\Rm{}\nwstartdeflinemarkup\nwusesondefline{\\{NWgqRer-36Ytqo-1}}\nwenddeflinemarkup
{\bf{}bool} {\it{}latexout};
{\bf{}dual\_number} {\it{}dn\_from\_arg\_mod}({\bf{}const} {\bf{}ex} & {\it{}a}, {\bf{}const} {\bf{}ex} & {\it{}n});
{\bf{}ex} {\it{}v\_from\_norm}({\bf{}const} {\bf{}ex} & {\it{}u}, {\bf{}const} {\bf{}ex} & {\it{}n});
{\bf{}void} {\it{}test\_out}({\bf{}char}\begin{math}\ast\end{math} {\it{}S}, {\bf{}ex} {\it{}T});
{\bf{}void} {\it{}formula\_out}({\bf{}char}\begin{math}\ast\end{math} {\it{}S}, {\bf{}ex} {\it{}F});

\nwused{\\{NWgqRer-36Ytqo-1}}\nwidentuses{\\{{\nwixident{dual{\_}number}}{dual:unnumber}}\\{{\nwixident{formula{\_}out}}{formula:unout}}\\{{\nwixident{test{\_}out}}{test:unout}}}\nwindexuse{\nwixident{dual{\_}number}}{dual:unnumber}{NWgqRer-NQ5ML-1}\nwindexuse{\nwixident{formula{\_}out}}{formula:unout}{NWgqRer-NQ5ML-1}\nwindexuse{\nwixident{test{\_}out}}{test:unout}{NWgqRer-NQ5ML-1}\nwendcode{}\nwbegindocs{217}\nwdocspar
\subsection{Header and Implementation of the {\Tt{}\Rm{}{\bf{}dual\_number}\nwendquote} Class}
\label{sec:impl-dual_n-class}

\nwenddocs{}\nwbegindocs{218}\nwdocspar

\subsubsection{Header File for {\Tt{}\Rm{}{\bf{}dual\_number}\nwendquote}}
\label{sec:head-file-dual_n}

\nwenddocs{}\nwbegindocs{219}We use the standard \GiNaC\ machinery do define {\Tt{}\Rm{}{\bf{}dual\_number}\nwendquote}s as
derived of the class {\Tt{}\Rm{}{\bf{}basic}\nwendquote}.
\nwenddocs{}\nwbegincode{220}\sublabel{NWgqRer-4GrAIY-1}\nwmargintag{{\nwtagstyle{}\subpageref{NWgqRer-4GrAIY-1}}}\moddef{Dual number class declaration~{\nwtagstyle{}\subpageref{NWgqRer-4GrAIY-1}}}\endmoddef\Rm{}\nwstartdeflinemarkup\nwusesondefline{\\{NWgqRer-36Ytqo-1}}\nwprevnextdefs{\relax}{NWgqRer-4GrAIY-2}\nwenddeflinemarkup
{\bf{}class} {\bf{}dual\_number} : {\bf{}public} {\bf{}basic}
{\nwlbrace}
        {\it{}GINAC\_DECLARE\_REGISTERED\_CLASS}({\bf{}dual\_number}, {\bf{}basic})
        {\bf{}static} {\bf{}const} {\it{}tinfo\_static\_t} {\it{}return\_type\_tinfo\_static}[256];

\nwalsodefined{\\{NWgqRer-4GrAIY-2}\\{NWgqRer-4GrAIY-3}\\{NWgqRer-4GrAIY-4}\\{NWgqRer-4GrAIY-5}\\{NWgqRer-4GrAIY-6}}\nwused{\\{NWgqRer-36Ytqo-1}}\nwidentuses{\\{{\nwixident{dual{\_}number}}{dual:unnumber}}\\{{\nwixident{tinfo{\_}static{\_}t}}{tinfo:unstatic:unt}}}\nwindexuse{\nwixident{dual{\_}number}}{dual:unnumber}{NWgqRer-4GrAIY-1}\nwindexuse{\nwixident{tinfo{\_}static{\_}t}}{tinfo:unstatic:unt}{NWgqRer-4GrAIY-1}\nwendcode{}\nwbegindocs{221}Public methods (constructors, algebraic, etc.)
\nwenddocs{}\nwbegincode{222}\sublabel{NWgqRer-4GrAIY-2}\nwmargintag{{\nwtagstyle{}\subpageref{NWgqRer-4GrAIY-2}}}\moddef{Dual number class declaration~{\nwtagstyle{}\subpageref{NWgqRer-4GrAIY-1}}}\plusendmoddef\Rm{}\nwstartdeflinemarkup\nwusesondefline{\\{NWgqRer-36Ytqo-1}}\nwprevnextdefs{NWgqRer-4GrAIY-1}{NWgqRer-4GrAIY-3}\nwenddeflinemarkup
{\bf{}public}:
        \LA{}Public methods~{\nwtagstyle{}\subpageref{NWgqRer-1eKCCy-1}}\RA{}
        \LA{}Technical methods~{\nwtagstyle{}\subpageref{NWgqRer-2XWNbs-1}}\RA{}

\nwused{\\{NWgqRer-36Ytqo-1}}\nwendcode{}\nwbegindocs{223}We redefine protected methods for printing only.
\nwenddocs{}\nwbegincode{224}\sublabel{NWgqRer-4GrAIY-3}\nwmargintag{{\nwtagstyle{}\subpageref{NWgqRer-4GrAIY-3}}}\moddef{Dual number class declaration~{\nwtagstyle{}\subpageref{NWgqRer-4GrAIY-1}}}\plusendmoddef\Rm{}\nwstartdeflinemarkup\nwusesondefline{\\{NWgqRer-36Ytqo-1}}\nwprevnextdefs{NWgqRer-4GrAIY-2}{NWgqRer-4GrAIY-4}\nwenddeflinemarkup
{\bf{}protected}:
        {\bf{}void} {\it{}do\_print}({\bf{}const} {\it{}print\_context} & {\it{}c}, {\bf{}unsigned} {\it{}level}) {\bf{}const};
        {\bf{}void} {\it{}do\_print\_latex}({\bf{}const} {\it{}print\_latex} & {\it{}c}, {\bf{}unsigned} {\it{}level}) {\bf{}const};

\nwused{\\{NWgqRer-36Ytqo-1}}\nwendcode{}\nwbegindocs{225}Private members: two components of a {\Tt{}\Rm{}{\bf{}dual\_number}\nwendquote}.
\nwenddocs{}\nwbegincode{226}\sublabel{NWgqRer-4GrAIY-4}\nwmargintag{{\nwtagstyle{}\subpageref{NWgqRer-4GrAIY-4}}}\moddef{Dual number class declaration~{\nwtagstyle{}\subpageref{NWgqRer-4GrAIY-1}}}\plusendmoddef\Rm{}\nwstartdeflinemarkup\nwusesondefline{\\{NWgqRer-36Ytqo-1}}\nwprevnextdefs{NWgqRer-4GrAIY-3}{NWgqRer-4GrAIY-5}\nwenddeflinemarkup
{\bf{}protected}: 
        {\bf{}ex} {\it{}u\_comp};
        {\bf{}ex} {\it{}v\_comp};
{\nwrbrace};

\nwused{\\{NWgqRer-36Ytqo-1}}\nwendcode{}\nwbegindocs{227}The following methods are needed for \GiNaC\ to work properly.
\nwenddocs{}\nwbegincode{228}\sublabel{NWgqRer-2XWNbs-1}\nwmargintag{{\nwtagstyle{}\subpageref{NWgqRer-2XWNbs-1}}}\moddef{Technical methods~{\nwtagstyle{}\subpageref{NWgqRer-2XWNbs-1}}}\endmoddef\Rm{}\nwstartdeflinemarkup\nwusesondefline{\\{NWgqRer-4GrAIY-2}}\nwenddeflinemarkup
        {\bf{}dual\_number} {\it{}normal}() {\bf{}const} {\nwlbrace} {\bf{}return} {\bf{}dual\_number}({\it{}u\_comp}.{\it{}normal}(), {\it{}v\_comp}.{\it{}normal}()); {\nwrbrace}
        {\bf{}dual\_number} {\it{}subs}({\bf{}const} {\bf{}ex} & {\it{}e}, {\bf{}unsigned} {\it{}options} = 0) {\bf{}const};
        {\bf{}bool} {\it{}is\_zero}() {\bf{}const};
        {\bf{}bool} {\it{}is\_equal}({\bf{}const} {\bf{}ex} & {\it{}other}) {\bf{}const};
        {\it{}size\_t} {\it{}nops}() {\bf{}const} {\nwlbrace} {\bf{}return} 2; {\nwrbrace}
        {\bf{}ex} {\it{}op}({\it{}size\_t} {\it{}i}) {\bf{}const};
        {\bf{}ex} & {\it{}let\_op}({\it{}size\_t} {\it{}i});

\nwused{\\{NWgqRer-4GrAIY-2}}\nwidentuses{\\{{\nwixident{dual{\_}number}}{dual:unnumber}}}\nwindexuse{\nwixident{dual{\_}number}}{dual:unnumber}{NWgqRer-2XWNbs-1}\nwendcode{}\nwbegindocs{229} We overload standard algebraic operations for {\Tt{}\Rm{}{\bf{}dual\_number}\nwendquote}.
\nwenddocs{}\nwbegincode{230}\sublabel{NWgqRer-4GrAIY-5}\nwmargintag{{\nwtagstyle{}\subpageref{NWgqRer-4GrAIY-5}}}\moddef{Dual number class declaration~{\nwtagstyle{}\subpageref{NWgqRer-4GrAIY-1}}}\plusendmoddef\Rm{}\nwstartdeflinemarkup\nwusesondefline{\\{NWgqRer-36Ytqo-1}}\nwprevnextdefs{NWgqRer-4GrAIY-4}{NWgqRer-4GrAIY-6}\nwenddeflinemarkup
{\bf{}const} {\bf{}dual\_number} {\bf{}operator}+({\bf{}const} {\bf{}dual\_number} & {\it{}lh}, {\bf{}const} {\bf{}dual\_number} & {\it{}rh});\nwindexdefn{\nwixident{dual{\_}number}}{dual:unnumber}{NWgqRer-4GrAIY-5}
{\bf{}const} {\bf{}dual\_number} {\bf{}operator}-({\bf{}const} {\bf{}dual\_number} & {\it{}lh}, {\bf{}const} {\bf{}dual\_number} & {\it{}rh});\nwindexdefn{\nwixident{dual{\_}number}}{dual:unnumber}{NWgqRer-4GrAIY-5}
{\bf{}const} {\bf{}dual\_number} {\bf{}operator}\begin{math}\ast\end{math}({\bf{}const} {\bf{}dual\_number} & {\it{}lh}, {\bf{}const} {\bf{}dual\_number} & {\it{}rh});\nwindexdefn{\nwixident{dual{\_}number}}{dual:unnumber}{NWgqRer-4GrAIY-5}
{\bf{}const} {\bf{}dual\_number} {\bf{}operator}\begin{math}\ast\end{math}({\bf{}const} {\bf{}dual\_number} & {\it{}lh}, {\bf{}const} {\bf{}ex} & {\it{}rh});\nwindexdefn{\nwixident{dual{\_}number}}{dual:unnumber}{NWgqRer-4GrAIY-5}
{\bf{}const} {\bf{}dual\_number} {\bf{}operator}\begin{math}\ast\end{math}({\bf{}const} {\bf{}ex} & {\it{}lh}, {\bf{}const} {\bf{}dual\_number} & {\it{}rh});\nwindexdefn{\nwixident{dual{\_}number}}{dual:unnumber}{NWgqRer-4GrAIY-5}
{\bf{}const} {\bf{}dual\_number} {\bf{}operator}\begin{math}\div\end{math}({\bf{}const} {\bf{}dual\_number} & {\it{}lh}, {\bf{}const} {\bf{}ex} & {\it{}rh});\nwindexdefn{\nwixident{dual{\_}number}}{dual:unnumber}{NWgqRer-4GrAIY-5}
{\bf{}const} {\bf{}dual\_number} {\bf{}operator}\begin{math}\div\end{math}({\bf{}const} {\bf{}ex} & {\it{}lh}, {\bf{}const} {\bf{}dual\_number} & {\it{}rh});\nwindexdefn{\nwixident{dual{\_}number}}{dual:unnumber}{NWgqRer-4GrAIY-5}
{\bf{}const} {\bf{}dual\_number} {\bf{}operator}\begin{math}\div\end{math}({\bf{}const} {\bf{}dual\_number} & {\it{}lh}, {\bf{}const} {\bf{}ex} & {\it{}rh});\nwindexdefn{\nwixident{dual{\_}number}}{dual:unnumber}{NWgqRer-4GrAIY-5}

\nwused{\\{NWgqRer-36Ytqo-1}}\nwidentdefs{\\{{\nwixident{dual{\_}number}}{dual:unnumber}}}\nwendcode{}\nwbegincode{231}\sublabel{NWgqRer-4GrAIY-6}\nwmargintag{{\nwtagstyle{}\subpageref{NWgqRer-4GrAIY-6}}}\moddef{Dual number class declaration~{\nwtagstyle{}\subpageref{NWgqRer-4GrAIY-1}}}\plusendmoddef\Rm{}\nwstartdeflinemarkup\nwusesondefline{\\{NWgqRer-36Ytqo-1}}\nwprevnextdefs{NWgqRer-4GrAIY-5}{\relax}\nwenddeflinemarkup
{\bf{}dual\_number} {\it{}dn\_from\_arg\_mod}({\bf{}const} {\bf{}ex} & {\it{}a}, {\bf{}const} {\bf{}ex} & {\it{}n});
// End of "header"

\nwused{\\{NWgqRer-36Ytqo-1}}\nwidentuses{\\{{\nwixident{dual{\_}number}}{dual:unnumber}}}\nwindexuse{\nwixident{dual{\_}number}}{dual:unnumber}{NWgqRer-4GrAIY-6}\nwendcode{}\nwbegindocs{232}\nwdocspar
\subsubsection{Standard Implementation Part}
\label{sec:stand-impl-part}

\nwenddocs{}\nwbegindocs{233}The implementation uses standard \GiNaC\ technique.
\nwenddocs{}\nwbegincode{234}\sublabel{NWgqRer-3uuIiL-1}\nwmargintag{{\nwtagstyle{}\subpageref{NWgqRer-3uuIiL-1}}}\moddef{Dual number class implementation~{\nwtagstyle{}\subpageref{NWgqRer-3uuIiL-1}}}\endmoddef\Rm{}\nwstartdeflinemarkup\nwusesondefline{\\{NWgqRer-3HhVai-1}}\nwprevnextdefs{\relax}{NWgqRer-3uuIiL-2}\nwenddeflinemarkup

{\it{}GINAC\_IMPLEMENT\_REGISTERED\_CLASS\_OPT}({\bf{}dual\_number}, {\bf{}basic}, 
                                                  {\it{}print\_func}\begin{math}<\end{math}{\it{}print\_context}\begin{math}>\end{math}(&{\bf{}dual\_number}::{\it{}do\_print}).
                                                  {\it{}print\_func}\begin{math}<\end{math}{\it{}print\_latex}\begin{math}>\end{math}(&{\bf{}dual\_number}::{\it{}do\_print\_latex}))

{\bf{}const} {\it{}tinfo\_static\_t} {\bf{}dual\_number}::{\it{}return\_type\_tinfo\_static}[256] = {\nwlbrace}{\nwlbrace}{\nwrbrace}{\nwrbrace};\nwindexdefn{\nwixident{tinfo{\_}static{\_}t}}{tinfo:unstatic:unt}{NWgqRer-3uuIiL-1}

{\it{}DEFAULT\_ARCHIVING}({\bf{}dual\_number})\nwindexdefn{\nwixident{dual{\_}number}}{dual:unnumber}{NWgqRer-3uuIiL-1}

\nwalsodefined{\\{NWgqRer-3uuIiL-2}\\{NWgqRer-3uuIiL-3}\\{NWgqRer-3uuIiL-4}\\{NWgqRer-3uuIiL-5}\\{NWgqRer-3uuIiL-6}\\{NWgqRer-3uuIiL-7}\\{NWgqRer-3uuIiL-8}\\{NWgqRer-3uuIiL-9}\\{NWgqRer-3uuIiL-A}\\{NWgqRer-3uuIiL-B}\\{NWgqRer-3uuIiL-C}\\{NWgqRer-3uuIiL-D}\\{NWgqRer-3uuIiL-E}\\{NWgqRer-3uuIiL-F}\\{NWgqRer-3uuIiL-G}}\nwused{\\{NWgqRer-3HhVai-1}}\nwidentdefs{\\{{\nwixident{dual{\_}number}}{dual:unnumber}}\\{{\nwixident{tinfo{\_}static{\_}t}}{tinfo:unstatic:unt}}}\nwendcode{}\nwbegindocs{235}\nwdocspar
\subsubsection{Implementation of Constructors}
\label{sec:impl-constr}
Default constructor.
\nwenddocs{}\nwbegincode{236}\sublabel{NWgqRer-3uuIiL-2}\nwmargintag{{\nwtagstyle{}\subpageref{NWgqRer-3uuIiL-2}}}\moddef{Dual number class implementation~{\nwtagstyle{}\subpageref{NWgqRer-3uuIiL-1}}}\plusendmoddef\Rm{}\nwstartdeflinemarkup\nwusesondefline{\\{NWgqRer-3HhVai-1}}\nwprevnextdefs{NWgqRer-3uuIiL-1}{NWgqRer-3uuIiL-3}\nwenddeflinemarkup
{\bf{}dual\_number}::{\bf{}dual\_number}() : {\it{}inherited}(&{\bf{}dual\_number}::{\it{}tinfo\_static}), {\it{}u\_comp}(0), {\it{}v\_comp}(0)
{\nwlbrace}
        {\it{}setflag}({\it{}status\_flags}::{\it{}not\_shareable});
{\nwrbrace}

\nwused{\\{NWgqRer-3HhVai-1}}\nwidentuses{\\{{\nwixident{dual{\_}number}}{dual:unnumber}}}\nwindexuse{\nwixident{dual{\_}number}}{dual:unnumber}{NWgqRer-3uuIiL-2}\nwendcode{}\nwbegindocs{237}Constructor from two components.
\nwenddocs{}\nwbegincode{238}\sublabel{NWgqRer-3uuIiL-3}\nwmargintag{{\nwtagstyle{}\subpageref{NWgqRer-3uuIiL-3}}}\moddef{Dual number class implementation~{\nwtagstyle{}\subpageref{NWgqRer-3uuIiL-1}}}\plusendmoddef\Rm{}\nwstartdeflinemarkup\nwusesondefline{\\{NWgqRer-3HhVai-1}}\nwprevnextdefs{NWgqRer-3uuIiL-2}{NWgqRer-3uuIiL-4}\nwenddeflinemarkup
{\bf{}dual\_number}::{\bf{}dual\_number}({\bf{}const} {\bf{}ex} & {\it{}a}, {\bf{}const} {\bf{}ex} & {\it{}b}) : {\it{}inherited}(&{\bf{}dual\_number}::{\it{}tinfo\_static}), 
{\it{}u\_comp}({\it{}a}), {\it{}v\_comp}({\it{}b})
{\nwlbrace}
{\nwrbrace}

\nwused{\\{NWgqRer-3HhVai-1}}\nwidentuses{\\{{\nwixident{dual{\_}number}}{dual:unnumber}}}\nwindexuse{\nwixident{dual{\_}number}}{dual:unnumber}{NWgqRer-3uuIiL-3}\nwendcode{}\nwbegindocs{239}Constructor from a single expression. It may contain two components\ldots
\nwenddocs{}\nwbegincode{240}\sublabel{NWgqRer-3uuIiL-4}\nwmargintag{{\nwtagstyle{}\subpageref{NWgqRer-3uuIiL-4}}}\moddef{Dual number class implementation~{\nwtagstyle{}\subpageref{NWgqRer-3uuIiL-1}}}\plusendmoddef\Rm{}\nwstartdeflinemarkup\nwusesondefline{\\{NWgqRer-3HhVai-1}}\nwprevnextdefs{NWgqRer-3uuIiL-3}{NWgqRer-3uuIiL-5}\nwenddeflinemarkup
{\bf{}dual\_number}::{\bf{}dual\_number}({\bf{}const} {\bf{}ex} & {\it{}P}) : {\it{}inherited}(&{\bf{}dual\_number}::{\it{}tinfo\_static})
{\nwlbrace}
        {\bf{}if} ({\it{}is\_a}\begin{math}<\end{math}{\bf{}lst}\begin{math}>\end{math}({\it{}P}) \begin{math}\vee\end{math} {\it{}is\_a}\begin{math}<\end{math}{\bf{}matrix}\begin{math}>\end{math}({\it{}P}) \begin{math}\vee\end{math} {\it{}is\_a}\begin{math}<\end{math}{\bf{}dual\_number}\begin{math}>\end{math}({\it{}P})) {\nwlbrace}
                {\it{}u\_comp} = {\it{}P}.{\it{}op}(0);
                {\it{}v\_comp} = {\it{}P}.{\it{}op}(1);

\nwused{\\{NWgqRer-3HhVai-1}}\nwidentuses{\\{{\nwixident{dual{\_}number}}{dual:unnumber}}}\nwindexuse{\nwixident{dual{\_}number}}{dual:unnumber}{NWgqRer-3uuIiL-4}\nwendcode{}\nwbegindocs{241}\ldots if it is a real expression we embed it into {\Tt{}\Rm{}{\bf{}dual\_number}\nwendquote}\ldots
\nwenddocs{}\nwbegincode{242}\sublabel{NWgqRer-3uuIiL-5}\nwmargintag{{\nwtagstyle{}\subpageref{NWgqRer-3uuIiL-5}}}\moddef{Dual number class implementation~{\nwtagstyle{}\subpageref{NWgqRer-3uuIiL-1}}}\plusendmoddef\Rm{}\nwstartdeflinemarkup\nwusesondefline{\\{NWgqRer-3HhVai-1}}\nwprevnextdefs{NWgqRer-3uuIiL-4}{NWgqRer-3uuIiL-6}\nwenddeflinemarkup
        {\nwrbrace} {\bf{}else} {\bf{}if} ({\it{}P}.{\it{}imag\_part}().{\it{}is\_zero}()) {\nwlbrace}
                {\it{}u\_comp} = {\it{}Arg0};
                {\it{}v\_comp} = {\it{}v\_from\_norm}({\it{}Arg0}, {\it{}P});

\nwused{\\{NWgqRer-3HhVai-1}}\nwidentuses{\\{{\nwixident{Arg0}}{Arg0}}}\nwindexuse{\nwixident{Arg0}}{Arg0}{NWgqRer-3uuIiL-5}\nwendcode{}\nwbegindocs{243}\ldots or if its a complex expression we decompose it into the real and
imaginary parts.
\nwenddocs{}\nwbegincode{244}\sublabel{NWgqRer-3uuIiL-6}\nwmargintag{{\nwtagstyle{}\subpageref{NWgqRer-3uuIiL-6}}}\moddef{Dual number class implementation~{\nwtagstyle{}\subpageref{NWgqRer-3uuIiL-1}}}\plusendmoddef\Rm{}\nwstartdeflinemarkup\nwusesondefline{\\{NWgqRer-3HhVai-1}}\nwprevnextdefs{NWgqRer-3uuIiL-5}{NWgqRer-3uuIiL-7}\nwenddeflinemarkup
        {\nwrbrace} {\bf{}else} {\nwlbrace}
                {\it{}u\_comp} = {\it{}P}.{\it{}real\_part}();
                {\it{}v\_comp} = {\it{}P}.{\it{}imag\_part}();
        {\nwrbrace}
{\nwrbrace}

\nwused{\\{NWgqRer-3HhVai-1}}\nwendcode{}\nwbegindocs{245}Comparison routine.
\nwenddocs{}\nwbegincode{246}\sublabel{NWgqRer-3uuIiL-7}\nwmargintag{{\nwtagstyle{}\subpageref{NWgqRer-3uuIiL-7}}}\moddef{Dual number class implementation~{\nwtagstyle{}\subpageref{NWgqRer-3uuIiL-1}}}\plusendmoddef\Rm{}\nwstartdeflinemarkup\nwusesondefline{\\{NWgqRer-3HhVai-1}}\nwprevnextdefs{NWgqRer-3uuIiL-6}{NWgqRer-3uuIiL-8}\nwenddeflinemarkup
{\bf{}int} {\bf{}dual\_number}::{\it{}compare\_same\_type}({\bf{}const} {\bf{}basic} & {\it{}other}) {\bf{}const}\nwindexdefn{\nwixident{dual{\_}number}}{dual:unnumber}{NWgqRer-3uuIiL-7}
{\nwlbrace}
        {\it{}GINAC\_ASSERT}({\it{}is\_a}\begin{math}<\end{math}{\bf{}dual\_number}\begin{math}>\end{math}({\it{}other}));
        {\bf{}const} {\bf{}dual\_number} &{\it{}o} = {\bf{}static\_cast}\begin{math}<\end{math}{\bf{}const} {\bf{}dual\_number} &\begin{math}>\end{math}({\it{}other});
        
        {\bf{}int} {\it{}cmpval} = {\it{}u\_comp}.{\it{}compare}({\it{}o}.{\it{}op}(0));
        {\bf{}if} ({\it{}cmpval}\begin{math}\neq\end{math}0) {\bf{}return} {\it{}cmpval};

        {\bf{}return} {\it{}v\_comp}.{\it{}compare}({\it{}o}.{\it{}op}(1));
{\nwrbrace}

\nwused{\\{NWgqRer-3HhVai-1}}\nwidentdefs{\\{{\nwixident{dual{\_}number}}{dual:unnumber}}}\nwendcode{}\nwbegindocs{247}Equality of two dual numbers.
\nwenddocs{}\nwbegincode{248}\sublabel{NWgqRer-3uuIiL-8}\nwmargintag{{\nwtagstyle{}\subpageref{NWgqRer-3uuIiL-8}}}\moddef{Dual number class implementation~{\nwtagstyle{}\subpageref{NWgqRer-3uuIiL-1}}}\plusendmoddef\Rm{}\nwstartdeflinemarkup\nwusesondefline{\\{NWgqRer-3HhVai-1}}\nwprevnextdefs{NWgqRer-3uuIiL-7}{NWgqRer-3uuIiL-9}\nwenddeflinemarkup
{\bf{}bool} {\bf{}dual\_number}::{\it{}is\_equal}({\bf{}const} {\bf{}ex} & {\it{}other}) {\bf{}const}
{\nwlbrace}
        {\it{}GINAC\_ASSERT}({\it{}is\_a}\begin{math}<\end{math}{\bf{}dual\_number}\begin{math}>\end{math}({\it{}other}));
        {\bf{}const} {\bf{}dual\_number} &{\it{}o} = {\bf{}static\_cast}\begin{math}<\end{math}{\bf{}const} {\bf{}dual\_number} &\begin{math}>\end{math}({\it{}other});
        
        {\bf{}return} {\it{}u\_comp}.{\it{}is\_equal}({\it{}o}.{\it{}op}(0)) \begin{math}\wedge\end{math} {\it{}v\_comp}.{\it{}is\_equal}({\it{}o}.{\it{}op}(1));
{\nwrbrace}

\nwused{\\{NWgqRer-3HhVai-1}}\nwidentuses{\\{{\nwixident{dual{\_}number}}{dual:unnumber}}}\nwindexuse{\nwixident{dual{\_}number}}{dual:unnumber}{NWgqRer-3uuIiL-8}\nwendcode{}\nwbegindocs{249}\nwdocspar
\subsubsection{Printing}
\label{sec:printing}
Default printing.
\nwenddocs{}\nwbegincode{250}\sublabel{NWgqRer-3uuIiL-9}\nwmargintag{{\nwtagstyle{}\subpageref{NWgqRer-3uuIiL-9}}}\moddef{Dual number class implementation~{\nwtagstyle{}\subpageref{NWgqRer-3uuIiL-1}}}\plusendmoddef\Rm{}\nwstartdeflinemarkup\nwusesondefline{\\{NWgqRer-3HhVai-1}}\nwprevnextdefs{NWgqRer-3uuIiL-8}{NWgqRer-3uuIiL-A}\nwenddeflinemarkup
{\bf{}void} {\bf{}dual\_number}::{\it{}do\_print}({\bf{}const} {\it{}print\_context} & {\it{}c}, {\bf{}unsigned} {\it{}level}) {\bf{}const}\nwindexdefn{\nwixident{dual{\_}number}}{dual:unnumber}{NWgqRer-3uuIiL-9}
{\nwlbrace}
        {\it{}c}.{\it{}s} \begin{math}\ll\end{math} {\tt{}"["};
        {\it{}u\_comp}.{\it{}print}({\it{}c});
        {\it{}c}.{\it{}s} \begin{math}\ll\end{math} {\tt{}","};
        {\it{}v\_comp}.{\it{}print}({\it{}c});
        {\it{}c}.{\it{}s} \begin{math}\ll\end{math} {\tt{}"]"};
{\nwrbrace}

\nwused{\\{NWgqRer-3HhVai-1}}\nwidentdefs{\\{{\nwixident{dual{\_}number}}{dual:unnumber}}}\nwendcode{}\nwbegindocs{251}\LaTeX\ printing.
\nwenddocs{}\nwbegincode{252}\sublabel{NWgqRer-3uuIiL-A}\nwmargintag{{\nwtagstyle{}\subpageref{NWgqRer-3uuIiL-A}}}\moddef{Dual number class implementation~{\nwtagstyle{}\subpageref{NWgqRer-3uuIiL-1}}}\plusendmoddef\Rm{}\nwstartdeflinemarkup\nwusesondefline{\\{NWgqRer-3HhVai-1}}\nwprevnextdefs{NWgqRer-3uuIiL-9}{NWgqRer-3uuIiL-B}\nwenddeflinemarkup
{\bf{}void} {\bf{}dual\_number}::{\it{}do\_print\_latex}({\bf{}const} {\it{}print\_latex} & {\it{}c}, {\bf{}unsigned} {\it{}level}) {\bf{}const}\nwindexdefn{\nwixident{dual{\_}number}}{dual:unnumber}{NWgqRer-3uuIiL-A}
{\nwlbrace}
        {\it{}c}.{\it{}s} \begin{math}\ll\end{math} {\tt{}"{\char92}{\char92}left({\char92}{\char92}begin{\char123}array{\char125}{\char123}cc{\char125}"};
        {\it{}u\_comp}.{\it{}print}({\it{}c});
        {\it{}c}.{\it{}s} \begin{math}\ll\end{math} {\tt{}"&"};
        {\it{}v\_comp}.{\it{}print}({\it{}c});
        {\it{}c}.{\it{}s} \begin{math}\ll\end{math} {\tt{}"{\char92}{\char92}end{\char123}array{\char125}{\char92}{\char92}right)"};
{\nwrbrace}

\nwused{\\{NWgqRer-3HhVai-1}}\nwidentdefs{\\{{\nwixident{dual{\_}number}}{dual:unnumber}}}\nwendcode{}\nwbegindocs{253}\nwdocspar
\subsubsection{Overloading algebraic operations}
\label{sec:overl-algebr-oper}
Addition.
\nwenddocs{}\nwbegincode{254}\sublabel{NWgqRer-3uuIiL-B}\nwmargintag{{\nwtagstyle{}\subpageref{NWgqRer-3uuIiL-B}}}\moddef{Dual number class implementation~{\nwtagstyle{}\subpageref{NWgqRer-3uuIiL-1}}}\plusendmoddef\Rm{}\nwstartdeflinemarkup\nwusesondefline{\\{NWgqRer-3HhVai-1}}\nwprevnextdefs{NWgqRer-3uuIiL-A}{NWgqRer-3uuIiL-C}\nwenddeflinemarkup
{\bf{}const} {\bf{}dual\_number} {\bf{}operator}+({\bf{}const} {\bf{}dual\_number} & {\it{}lh}, {\bf{}const} {\bf{}dual\_number} & {\it{}rh}) \nwindexdefn{\nwixident{dual{\_}number}}{dual:unnumber}{NWgqRer-3uuIiL-B}
{\nwlbrace}
        {\bf{}return} {\it{}lh}.{\it{}add}({\it{}rh});
{\nwrbrace}

{\bf{}const} {\bf{}dual\_number} {\bf{}operator}-({\bf{}const} {\bf{}dual\_number} & {\it{}lh}, {\bf{}const} {\bf{}dual\_number} & {\it{}rh})\nwindexdefn{\nwixident{dual{\_}number}}{dual:unnumber}{NWgqRer-3uuIiL-B}
{\nwlbrace}
        {\bf{}return} {\it{}lh}.{\it{}sub}({\it{}rh});
{\nwrbrace}

\nwused{\\{NWgqRer-3HhVai-1}}\nwidentdefs{\\{{\nwixident{dual{\_}number}}{dual:unnumber}}}\nwendcode{}\nwbegindocs{255}Multiplication.
\nwenddocs{}\nwbegincode{256}\sublabel{NWgqRer-3uuIiL-C}\nwmargintag{{\nwtagstyle{}\subpageref{NWgqRer-3uuIiL-C}}}\moddef{Dual number class implementation~{\nwtagstyle{}\subpageref{NWgqRer-3uuIiL-1}}}\plusendmoddef\Rm{}\nwstartdeflinemarkup\nwusesondefline{\\{NWgqRer-3HhVai-1}}\nwprevnextdefs{NWgqRer-3uuIiL-B}{NWgqRer-3uuIiL-D}\nwenddeflinemarkup
{\bf{}const} {\bf{}dual\_number} {\bf{}operator}\begin{math}\ast\end{math}({\bf{}const} {\bf{}dual\_number} & {\it{}lh}, {\bf{}const} {\bf{}ex} & {\it{}rh})\nwindexdefn{\nwixident{dual{\_}number}}{dual:unnumber}{NWgqRer-3uuIiL-C}
{\nwlbrace}
        {\bf{}return} {\it{}lh}.{\it{}mul}({\it{}rh});
{\nwrbrace}

{\bf{}const} {\bf{}dual\_number} {\bf{}operator}\begin{math}\ast\end{math}({\bf{}const} {\bf{}ex} & {\it{}lh}, {\bf{}const} {\bf{}dual\_number} & {\it{}rh})\nwindexdefn{\nwixident{dual{\_}number}}{dual:unnumber}{NWgqRer-3uuIiL-C}
{\nwlbrace}
        {\bf{}return} {\it{}rh}.{\it{}mul}({\it{}lh});
{\nwrbrace}

{\bf{}const} {\bf{}dual\_number} {\bf{}operator}\begin{math}\ast\end{math}({\bf{}const} {\bf{}dual\_number} & {\it{}lh}, {\bf{}const} {\bf{}dual\_number} & {\it{}rh})\nwindexdefn{\nwixident{dual{\_}number}}{dual:unnumber}{NWgqRer-3uuIiL-C}
{\nwlbrace}
        {\bf{}return} {\it{}lh}.{\it{}mul}({\it{}rh});
{\nwrbrace}

\nwused{\\{NWgqRer-3HhVai-1}}\nwidentdefs{\\{{\nwixident{dual{\_}number}}{dual:unnumber}}}\nwendcode{}\nwbegindocs{257}Division.
\nwenddocs{}\nwbegincode{258}\sublabel{NWgqRer-3uuIiL-D}\nwmargintag{{\nwtagstyle{}\subpageref{NWgqRer-3uuIiL-D}}}\moddef{Dual number class implementation~{\nwtagstyle{}\subpageref{NWgqRer-3uuIiL-1}}}\plusendmoddef\Rm{}\nwstartdeflinemarkup\nwusesondefline{\\{NWgqRer-3HhVai-1}}\nwprevnextdefs{NWgqRer-3uuIiL-C}{NWgqRer-3uuIiL-E}\nwenddeflinemarkup
{\bf{}const} {\bf{}dual\_number} {\bf{}operator}\begin{math}\div\end{math}({\bf{}const} {\bf{}dual\_number} & {\it{}lh}, {\bf{}const} {\bf{}dual\_number} & {\it{}rh})\nwindexdefn{\nwixident{dual{\_}number}}{dual:unnumber}{NWgqRer-3uuIiL-D}
{\nwlbrace}
        {\bf{}return} {\it{}lh}.{\it{}mul}({\it{}rh}.{\it{}power}(-1));
{\nwrbrace}

{\bf{}const} {\bf{}dual\_number} {\bf{}operator}\begin{math}\div\end{math}({\bf{}const} {\bf{}ex} & {\it{}lh}, {\bf{}const} {\bf{}dual\_number} & {\it{}rh})\nwindexdefn{\nwixident{dual{\_}number}}{dual:unnumber}{NWgqRer-3uuIiL-D}
{\nwlbrace}
        {\bf{}return} {\it{}rh}.{\it{}power}(-1)\begin{math}\ast\end{math}{\it{}lh};
{\nwrbrace}

{\bf{}const} {\bf{}dual\_number} {\bf{}operator}\begin{math}\div\end{math}({\bf{}const} {\bf{}dual\_number} & {\it{}lh}, {\bf{}const} {\bf{}ex} & {\it{}rh})\nwindexdefn{\nwixident{dual{\_}number}}{dual:unnumber}{NWgqRer-3uuIiL-D}
{\nwlbrace}
        {\bf{}return} {\it{}lh}.{\it{}mul}({\it{}pow}({\it{}rh}, -1));
{\nwrbrace}

\nwused{\\{NWgqRer-3HhVai-1}}\nwidentdefs{\\{{\nwixident{dual{\_}number}}{dual:unnumber}}}\nwendcode{}\nwbegindocs{259}\nwdocspar

\subsubsection{Component-related functions}
\label{sec:comp-relat-funct}

\nwenddocs{}\nwbegincode{260}\sublabel{NWgqRer-3uuIiL-E}\nwmargintag{{\nwtagstyle{}\subpageref{NWgqRer-3uuIiL-E}}}\moddef{Dual number class implementation~{\nwtagstyle{}\subpageref{NWgqRer-3uuIiL-1}}}\plusendmoddef\Rm{}\nwstartdeflinemarkup\nwusesondefline{\\{NWgqRer-3HhVai-1}}\nwprevnextdefs{NWgqRer-3uuIiL-D}{NWgqRer-3uuIiL-F}\nwenddeflinemarkup
{\bf{}ex} {\bf{}dual\_number}::{\it{}op}({\it{}size\_t} {\it{}i}) {\bf{}const}
{\nwlbrace}
        {\it{}GINAC\_ASSERT}({\it{}i}\begin{math}<\end{math}{\it{}nops}());
        
        {\bf{}switch} ({\it{}i}) {\nwlbrace}
        {\bf{}case} 0:
                {\bf{}return} {\it{}u\_comp};
        {\bf{}case} 1:
                {\bf{}return} {\it{}v\_comp};
        {\bf{}default}:
                {\bf{}throw}({\it{}std}::{\it{}invalid\_argument}({\tt{}"dual\_number::op(): requested"} 
                                                                        {\tt{}" operand out of the range (2)"}));
        {\nwrbrace}
{\nwrbrace}

\nwused{\\{NWgqRer-3HhVai-1}}\nwidentuses{\\{{\nwixident{dual{\_}number}}{dual:unnumber}}}\nwindexuse{\nwixident{dual{\_}number}}{dual:unnumber}{NWgqRer-3uuIiL-E}\nwendcode{}\nwbegindocs{261}\nwdocspar
\nwenddocs{}\nwbegincode{262}\sublabel{NWgqRer-3uuIiL-F}\nwmargintag{{\nwtagstyle{}\subpageref{NWgqRer-3uuIiL-F}}}\moddef{Dual number class implementation~{\nwtagstyle{}\subpageref{NWgqRer-3uuIiL-1}}}\plusendmoddef\Rm{}\nwstartdeflinemarkup\nwusesondefline{\\{NWgqRer-3HhVai-1}}\nwprevnextdefs{NWgqRer-3uuIiL-E}{NWgqRer-3uuIiL-G}\nwenddeflinemarkup
{\bf{}dual\_number} {\bf{}dual\_number}::{\it{}subs}({\bf{}const} {\bf{}ex} & {\it{}e}, {\bf{}unsigned} {\it{}options}) {\bf{}const}
{\nwlbrace}
        {\it{}exmap} {\it{}m};
        {\bf{}if} ({\it{}e}.{\it{}info}({\it{}info\_flags}::{\it{}list})) {\nwlbrace}
                {\bf{}lst} {\it{}l} = {\it{}ex\_to}\begin{math}<\end{math}{\bf{}lst}\begin{math}>\end{math}({\it{}e});
                {\bf{}for} ({\bf{}lst}::{\it{}const\_iterator} {\it{}i} = {\it{}l}.{\it{}begin}(); {\it{}i} \begin{math}\neq\end{math} {\it{}l}.{\it{}end}(); \protect\PP{\it{}i})
                        {\it{}m}.{\it{}insert}({\it{}std}::{\it{}make\_pair}({\it{}i}\begin{math}\rightarrow\end{math}{\it{}op}(0), {\it{}i}\begin{math}\rightarrow\end{math}{\it{}op}(1)));
        {\nwrbrace} {\bf{}else} {\bf{}if} ({\it{}is\_a}\begin{math}<\end{math}{\bf{}relational}\begin{math}>\end{math}({\it{}e})) {\nwlbrace}
                {\it{}m}.{\it{}insert}({\it{}std}::{\it{}make\_pair}({\it{}e}.{\it{}op}(0), {\it{}e}.{\it{}op}(1)));
        {\nwrbrace} {\bf{}else}
                {\bf{}throw}({\it{}std}::{\it{}invalid\_argument}({\tt{}"dual\_number::subs(): the parameter"}
                                                                        {\tt{}" should be a relational or a lst"}));
        
 {\bf{}return} {\it{}ex\_to}\begin{math}<\end{math}{\bf{}dual\_number}\begin{math}>\end{math}({\it{}inherited}::{\it{}subs}({\it{}m}, {\it{}options}));
{\nwrbrace}

\nwused{\\{NWgqRer-3HhVai-1}}\nwidentuses{\\{{\nwixident{dual{\_}number}}{dual:unnumber}}}\nwindexuse{\nwixident{dual{\_}number}}{dual:unnumber}{NWgqRer-3uuIiL-F}\nwendcode{}\nwbegindocs{263}\nwdocspar
\nwenddocs{}\nwbegincode{264}\sublabel{NWgqRer-3uuIiL-G}\nwmargintag{{\nwtagstyle{}\subpageref{NWgqRer-3uuIiL-G}}}\moddef{Dual number class implementation~{\nwtagstyle{}\subpageref{NWgqRer-3uuIiL-1}}}\plusendmoddef\Rm{}\nwstartdeflinemarkup\nwusesondefline{\\{NWgqRer-3HhVai-1}}\nwprevnextdefs{NWgqRer-3uuIiL-F}{\relax}\nwenddeflinemarkup
{\bf{}ex} & {\bf{}dual\_number}::{\it{}let\_op}({\it{}size\_t} {\it{}i})
{\nwlbrace}
        {\it{}GINAC\_ASSERT}({\it{}i}\begin{math}<\end{math}{\it{}nops}());
        
        {\it{}ensure\_if\_modifiable}();
        {\bf{}switch} ({\it{}i}) {\nwlbrace}
        {\bf{}case} 0:
                {\bf{}return} {\it{}u\_comp};
        {\bf{}case} 1:
                {\bf{}return} {\it{}v\_comp};
        {\bf{}default}:
                {\bf{}throw}({\it{}std}::{\it{}invalid\_argument}({\tt{}"dual\_number::let\_op(): requested operand"}
                                                                        {\tt{}" out of the range (2)"}));
        {\nwrbrace}
{\nwrbrace}

\nwused{\\{NWgqRer-3HhVai-1}}\nwidentuses{\\{{\nwixident{dual{\_}number}}{dual:unnumber}}}\nwindexuse{\nwixident{dual{\_}number}}{dual:unnumber}{NWgqRer-3uuIiL-G}\nwendcode{}

\nwixlogsorted{c}{{*}{NWgqRer-1p0Y9w-1}{\nwixd{NWgqRer-1p0Y9w-1}\nwixd{NWgqRer-1p0Y9w-2}}}%
\nwixlogsorted{c}{{Additional routines declarations}{NWgqRer-NQ5ML-1}{\nwixu{NWgqRer-36Ytqo-1}\nwixd{NWgqRer-NQ5ML-1}}}%
\nwixlogsorted{c}{{Algebraic procedures}{NWgqRer-3jzoaG-1}{\nwixd{NWgqRer-3jzoaG-1}\nwixd{NWgqRer-3jzoaG-2}\nwixd{NWgqRer-3jzoaG-3}\nwixd{NWgqRer-3jzoaG-4}\nwixd{NWgqRer-3jzoaG-5}\nwixd{NWgqRer-3jzoaG-6}\nwixd{NWgqRer-3jzoaG-7}\nwixu{NWgqRer-3HhVai-1}}}%
\nwixlogsorted{c}{{Check identities}{NWgqRer-22TNkn-1}{\nwixd{NWgqRer-22TNkn-1}\nwixd{NWgqRer-22TNkn-2}\nwixd{NWgqRer-22TNkn-3}\nwixd{NWgqRer-22TNkn-4}\nwixd{NWgqRer-22TNkn-5}\nwixd{NWgqRer-22TNkn-6}\nwixd{NWgqRer-22TNkn-7}\nwixd{NWgqRer-22TNkn-8}\nwixd{NWgqRer-22TNkn-9}\nwixd{NWgqRer-22TNkn-A}\nwixd{NWgqRer-22TNkn-B}\nwixd{NWgqRer-22TNkn-C}\nwixd{NWgqRer-22TNkn-D}\nwixd{NWgqRer-22TNkn-E}\nwixd{NWgqRer-22TNkn-F}\nwixu{NWgqRer-nXe8t-4}}}%
\nwixlogsorted{c}{{Definition of variables}{NWgqRer-2iImO3-1}{\nwixu{NWgqRer-1p0Y9w-2}\nwixd{NWgqRer-2iImO3-1}\nwixd{NWgqRer-2iImO3-2}\nwixd{NWgqRer-2iImO3-3}\nwixd{NWgqRer-2iImO3-4}\nwixd{NWgqRer-2iImO3-5}}}%
\nwixlogsorted{c}{{Drawing pictures}{NWgqRer-Ghke2-1}{\nwixu{NWgqRer-3sxui-5}\nwixd{NWgqRer-Ghke2-1}\nwixd{NWgqRer-Ghke2-2}\nwixd{NWgqRer-Ghke2-3}\nwixd{NWgqRer-Ghke2-4}}}%
\nwixlogsorted{c}{{Dual number class declaration}{NWgqRer-4GrAIY-1}{\nwixu{NWgqRer-36Ytqo-1}\nwixd{NWgqRer-4GrAIY-1}\nwixd{NWgqRer-4GrAIY-2}\nwixd{NWgqRer-4GrAIY-3}\nwixd{NWgqRer-4GrAIY-4}\nwixd{NWgqRer-4GrAIY-5}\nwixd{NWgqRer-4GrAIY-6}}}%
\nwixlogsorted{c}{{Dual number class further implementation}{NWgqRer-4dxCds-1}{\nwixd{NWgqRer-4dxCds-1}\nwixd{NWgqRer-4dxCds-2}\nwixd{NWgqRer-4dxCds-3}\nwixd{NWgqRer-4dxCds-4}\nwixd{NWgqRer-4dxCds-5}\nwixd{NWgqRer-4dxCds-6}\nwixd{NWgqRer-4dxCds-7}\nwixd{NWgqRer-4dxCds-8}\nwixd{NWgqRer-4dxCds-9}\nwixu{NWgqRer-3HhVai-1}}}%
\nwixlogsorted{c}{{Dual number class implementation}{NWgqRer-3uuIiL-1}{\nwixu{NWgqRer-3HhVai-1}\nwixd{NWgqRer-3uuIiL-1}\nwixd{NWgqRer-3uuIiL-2}\nwixd{NWgqRer-3uuIiL-3}\nwixd{NWgqRer-3uuIiL-4}\nwixd{NWgqRer-3uuIiL-5}\nwixd{NWgqRer-3uuIiL-6}\nwixd{NWgqRer-3uuIiL-7}\nwixd{NWgqRer-3uuIiL-8}\nwixd{NWgqRer-3uuIiL-9}\nwixd{NWgqRer-3uuIiL-A}\nwixd{NWgqRer-3uuIiL-B}\nwixd{NWgqRer-3uuIiL-C}\nwixd{NWgqRer-3uuIiL-D}\nwixd{NWgqRer-3uuIiL-E}\nwixd{NWgqRer-3uuIiL-F}\nwixd{NWgqRer-3uuIiL-G}}}%
\nwixlogsorted{c}{{dualnum.cpp}{NWgqRer-3HhVai-1}{\nwixd{NWgqRer-3HhVai-1}}}%
\nwixlogsorted{c}{{dualnum.h}{NWgqRer-36Ytqo-1}{\nwixd{NWgqRer-36Ytqo-1}}}%
\nwixlogsorted{c}{{Induced representations}{NWgqRer-1h89LW-1}{\nwixd{NWgqRer-1h89LW-1}\nwixd{NWgqRer-1h89LW-2}\nwixd{NWgqRer-1h89LW-3}\nwixd{NWgqRer-1h89LW-4}\nwixd{NWgqRer-1h89LW-5}\nwixd{NWgqRer-1h89LW-6}\nwixd{NWgqRer-1h89LW-7}\nwixd{NWgqRer-1h89LW-8}\nwixu{NWgqRer-3sxui-4}}}%
\nwixlogsorted{c}{{Induced representations routines}{NWgqRer-11aTCz-1}{\nwixd{NWgqRer-11aTCz-1}\nwixd{NWgqRer-11aTCz-2}\nwixd{NWgqRer-11aTCz-3}\nwixd{NWgqRer-11aTCz-4}\nwixu{NWgqRer-1p0Y9w-2}}}%
\nwixlogsorted{c}{{Initialisation}{NWgqRer-40D9Pp-1}{\nwixu{NWgqRer-36Ytqo-1}\nwixd{NWgqRer-40D9Pp-1}}}%
\nwixlogsorted{c}{{Main procedure}{NWgqRer-3sxui-1}{\nwixu{NWgqRer-1p0Y9w-2}\nwixd{NWgqRer-3sxui-1}\nwixd{NWgqRer-3sxui-2}\nwixd{NWgqRer-3sxui-3}\nwixd{NWgqRer-3sxui-4}\nwixd{NWgqRer-3sxui-5}}}%
\nwixlogsorted{c}{{N-Nprime separation}{NWgqRer-4EWpKB-1}{\nwixd{NWgqRer-4EWpKB-1}\nwixd{NWgqRer-4EWpKB-2}\nwixd{NWgqRer-4EWpKB-3}\nwixu{NWgqRer-36Ytqo-1}}}%
\nwixlogsorted{c}{{Output routines}{NWgqRer-4MrGQl-1}{\nwixu{NWgqRer-3HhVai-1}\nwixd{NWgqRer-4MrGQl-1}\nwixd{NWgqRer-4MrGQl-2}}}%
\nwixlogsorted{c}{{Public methods}{NWgqRer-1eKCCy-1}{\nwixd{NWgqRer-1eKCCy-1}\nwixd{NWgqRer-1eKCCy-2}\nwixd{NWgqRer-1eKCCy-3}\nwixd{NWgqRer-1eKCCy-4}\nwixd{NWgqRer-1eKCCy-5}\nwixd{NWgqRer-1eKCCy-6}\nwixd{NWgqRer-1eKCCy-7}\nwixd{NWgqRer-1eKCCy-8}\nwixd{NWgqRer-1eKCCy-9}\nwixu{NWgqRer-4GrAIY-2}}}%
\nwixlogsorted{c}{{Show expressions}{NWgqRer-1Pc9Jw-1}{\nwixd{NWgqRer-1Pc9Jw-1}\nwixd{NWgqRer-1Pc9Jw-2}\nwixd{NWgqRer-1Pc9Jw-3}\nwixd{NWgqRer-1Pc9Jw-4}\nwixd{NWgqRer-1Pc9Jw-5}\nwixd{NWgqRer-1Pc9Jw-6}\nwixd{NWgqRer-1Pc9Jw-7}\nwixd{NWgqRer-1Pc9Jw-8}\nwixd{NWgqRer-1Pc9Jw-9}\nwixd{NWgqRer-1Pc9Jw-A}\nwixu{NWgqRer-nXe8t-4}}}%
\nwixlogsorted{c}{{Technical methods}{NWgqRer-2XWNbs-1}{\nwixu{NWgqRer-4GrAIY-2}\nwixd{NWgqRer-2XWNbs-1}}}%
\nwixlogsorted{c}{{Test routine}{NWgqRer-nXe8t-1}{\nwixu{NWgqRer-1p0Y9w-2}\nwixd{NWgqRer-nXe8t-1}\nwixd{NWgqRer-nXe8t-2}\nwixd{NWgqRer-nXe8t-3}\nwixd{NWgqRer-nXe8t-4}}}%
\nwixlogsorted{i}{{\nwixident{Arg0}}{Arg0}}%
\nwixlogsorted{i}{{\nwixident{cases}}{cases}}%
\nwixlogsorted{i}{{\nwixident{dual{\_}number}}{dual:unnumber}}%
\nwixlogsorted{i}{{\nwixident{formula{\_}out}}{formula:unout}}%
\nwixlogsorted{i}{{\nwixident{main}}{main}}%
\nwixlogsorted{i}{{\nwixident{parab{\_}rot{\_}sub}}{parab:unrot:unsub}}%
\nwixlogsorted{i}{{\nwixident{test{\_}out}}{test:unout}}%
\nwixlogsorted{i}{{\nwixident{tinfo{\_}static{\_}t}}{tinfo:unstatic:unt}}%
\nwbegindocs{265}This is the end of {\Tt{}\Rm{}{\bf{}dual\_number}\nwendquote} implementation. 
\nwenddocs{}
}{}

\end{document}